\documentclass{amsart}

\usepackage[T1]{fontenc}
\usepackage[utf8]{inputenc}
\usepackage[american]{babel}

\usepackage{amsmath}
\usepackage{amsthm}
\usepackage{amssymb}
\usepackage[only,rrbracket,llbracket]{stmaryrd}

\theoremstyle{plain}
\newtheorem{thm}{Theorem}[section]
\newtheorem{lem}[thm]{Lemma}
\newtheorem{cor}[thm]{Corollary}
\theoremstyle{remark}
\newtheorem{rem}[thm]{Remark}

\usepackage[draft]{hyperref}

\usepackage{lmodern}
\usepackage[scr=rsfs]{mathalfa}

\newcount\afterconjunctionpenalty
\makeatletter
\@namedef{u8:\detokenize{⋅}}{\penalty\afterconjunctionpenalty\ }
\makeatother

\newcommand*\hairspace{\kern 0.08333em}
\newcommand*\thinglue{\nobreak\hspace{.16667em plus .08333em}}

\let\op\mathcal

\let\wt\tilde
\let\wb\bar
\let\wh\hat

\usepackage{mathtools}

\renewcommand{\Re}{\operatorname{Re}}
\renewcommand{\Im}{\operatorname{Im}}
\providecommand{\rbr}[1]{(#1)}
\providecommand{\bigrbr}[1]{\bigl(#1\bigr)}
\providecommand{\Bigrbr}[1]{\Bigl(#1\Bigr)}
\providecommand{\biggrbr}[1]{\biggl(#1\biggr)}
\providecommand{\sbr}[1]{[#1]}
\providecommand{\cbr}[1]{\{#1\}}
\providecommand{\bigcbr}[1]{\bigl\{#1\bigr\}}
\providecommand{\abr}[1]{\langle#1\rangle}
\providecommand{\abs}[1]{\lvert#1\rvert}
\providecommand{\bigabs}[1]{\bigl\lvert#1\bigr\rvert}
\providecommand{\norm}[1]{\lVert#1\rVert}
\providecommand{\Norm}[1]{\mathopen{|\mkern-1.5mu|\mkern-1.5mu|}#1\mathclose{|\mkern-1.5mu|\mkern-1.5mu|}}
\providecommand{\ceil}[1]{\lceil#1\rceil}
\let\seminorm\sbr
\providecommand{\Seminorm}[1]{\llbracket#1\rrbracket}
\newcommand\restrict[1]{\raisebox{-0.1ex}{$\rvert_{#1}$}}
\DeclareMathOperator{\diam}{diam}
\DeclareMathOperator{\dist}{dist}
\DeclareMathOperator{\Div}{div}
\DeclareMathOperator{\rank}{rank}
\DeclareMathOperator{\supp}{supp}
\newcommand*{\findiff}{\varDelta}
\newcommand*{\join}{{\vee}}
\newcommand*{\meet}{{\wedge}}

\let\set\mathscr
\def\B{\mathbf{B}}
\def\C{\mathbb{C}}
\def\N{\mathbb{N}}
\def\R{\mathbb{R}}
\def\T{\mathbb{T}}
\def\Z{\mathbb{Z}}

\DeclareMathOperator{\BMO}{BMO}
\DeclareMathOperator{\VMO}{VMO}

\newcommand*{\osc}{\omega}

\newcommand*{\wc}{{\mkern 2mu\cdot\mkern 2mu}}

\newcommand*{\dd}{\mathop{}\!d}

\newcommand*{\dual}{{+}}
\newcommand*{\dualref}[1]{\hyperref[#1]{\ref*{#1}\textsuperscript{+}}}
\newcommand*{\cell}{Q}
\newcommand*{\cellradius}{r_\cell}
\newcommand*{\chf}{\chi}

\def\Xint#1{\,\mathchoice
{\XXint\displaystyle\textstyle{#1}}%
{\XXint\textstyle\scriptstyle{#1}}%
{\XXint\scriptstyle\scriptscriptstyle{#1}}%
{\XXint\scriptscriptstyle\scriptscriptstyle{#1}}%
\!\int}
\def\XXint#1#2#3{{\setbox0=\hbox{$#1{#2#3}{\int}$ }
\vcenter{\hbox{$#2#3$ }}\kern-.6\wd0}}
\def\dashint{\Xint-}

\begin{document}
\frenchspacing

\title[Homogenization for locally periodic elliptic problems]
      {Homogenization for\\ locally periodic elliptic problems on a domain}
\author{Nikita N. Senik}
\address{St.~Petersburg State University, Universitetskaya nab.~7/9, St.~Petersburg 199034, Russia}
\email{nnsenik@gmail.com}
\thanks{The~research was supported by Russian Science Foundation grant~17-11-01069.}


\date{\today.}


\begin{abstract}
Let $\Omega$ be a Lipschitz domain in $\mathbb R^d$, and let
$\mathcal A^\varepsilon=-\operatorname{div}A(x,x/\varepsilon)\nabla$ be
a strongly elliptic operator on $\Omega$. We suppose that $\varepsilon$ is
small and the function~$A$ is Lipschitz in the first variable and periodic
in the second, so the coefficients of $\mathcal A^\varepsilon$ are locally
periodic and rapidly oscillate. Given $\mu$ in the resolvent set, we are
interested in finding the rates of approximations, as $\varepsilon\to0$,
for $(\mathcal A^\varepsilon-\mu)^{-1}$ and
$\nabla(\mathcal A^\varepsilon-\mu)^{-1}$ in the operator topology on $L_p$
for suitable~$p$. It is well-known that the rates depend on regularity of
the effective operator~$\mathcal A^0$. We prove that if $(\mathcal A^0-\mu)^{-1}$
and its adjoint are bounded from $L_p(\Omega)^n$ to the Lipschitz--Besov space~%
$\Lambda_p^{1+s}(\Omega)^n$ with $s\in(0,1]$, then the rates are, respectively,
$\varepsilon^s$ and $\varepsilon^{s/p}$. The results are applied to the
Dirichlet, Neumann and mixed Dirichlet--Neumann problems for strongly elliptic
operators with uniformly bounded and $\operatorname{VMO}$ coefficients.
\end{abstract}

\maketitle

\section{Introduction}

Let⋅$\Omega$ be a bounded domain, and let $A\colon\Omega\times\R^{d}\to\C^{d\times d}$
be a uniformly elliptic function which is smooth in the first variable
and periodic  in the second. A⋅classical result in homogenization
theory tells us that, for any⋅$f$ in⋅$H^{-1}\rbr{\Omega}$, the
dual of the Sobolev space~$\mathring{H}^{1}\rbr{\Omega}$, the solution~$u_{\varepsilon}\in\mathring{H}^{1}\rbr{\Omega}$
of the Dirichlet problem
\begin{equation}
\begin{aligned}-\Div A\rbr{x,\varepsilon^{-1}x}\hairspace\nabla u_{\varepsilon}=f & \qquad\text{in }\Omega,\\
u_{\varepsilon}=\mathrlap{0}\hphantom{f} & \qquad\text{on }\partial\Omega,
\end{aligned}
\label{eq: Intro | Heterogeneous equation}
\end{equation}
converges, as⋅$\varepsilon\to0$,  to the solution~$u_{0}$ of
a similar problem
\begin{equation}
\begin{aligned}-\Div A^{0}\rbr x\hairspace\nabla u_{0}=f & \qquad\text{in }\Omega,\\
u_{0}=\mathrlap{0}\hphantom{f} & \qquad\text{on }\partial\Omega,
\end{aligned}
\label{eq: Intro | Homogeneous equation}
\end{equation}
where $A^{0}\colon\Omega\to\C^{d\times d}$ is a smooth function.
  In⋅applications, this usually is interpreted as approximation
of a highly heterogeneous medium, described by the rapidly oscillating
locally periodic function~$x\mapsto A\rbr{x,\varepsilon^{-1}x}$,
with a homogeneous one, described by the  slowly varying function~$x\mapsto A^{0}\rbr x$.

There are various ways to prove the convergence. Among the first were
the method of asymptotic expansions, using powerful tools of asymptotic
analysis (see⋅\cite{BLP:1978} or~\cite{BP:1984}), and the energy
method, based on the compensated compactness phenomenon (see~\cite{MT:1997}
and also⋅\cite{Tartar:2010}). Another way of dealing with the problem~(\ref{eq: Intro | Heterogeneous equation})
is to use the two-scale convergence technique, which was developed
later (see, e.g.,~\cite{Al:1992}). In⋅any case, one finds that
$u_{\varepsilon}$ converges to⋅$u_{0}$ weakly in the Sobolev space~$\mathring{H}^{1}\rbr{\Omega}$,
and therefore strongly in the Lebesgue space~$L_{2}\rbr{\Omega}$.
The⋅latter  can be phrased as saying that the resolvent of the
operator⋅$-\Div A\rbr{x,\varepsilon^{-1}x}\hairspace\nabla$ converges
to the resolvent of the operator⋅$-\Div A^{0}\rbr x\hairspace\nabla$
in the respective strong operator topology. A⋅simple argument, see~\cite{AlaireConca:1998},
 using a compact embedding theorem then shows that the resolvent
converges in the uniform operator topology on⋅$L_{2}\rbr{\Omega}$,
 the strongest operator topology on⋅$L_{2}\rbr{\Omega}$. However,
this says nothing about the rate of convergence, nor does it apply
to the case of unbounded⋅$\Omega$ (or⋅quasi-bounded, to be precise;
see~\cite[Section~6.9]{AF:2003}).

A⋅sharp-order bound on the rate  was found in the pioneering paper~\cite{BSu:2001}
(see⋅also~\cite{BSu:2003}) for a purely periodic problem (when
the coefficients depend on $x/\varepsilon$ only) on~$\R^{d}$. Uniform
operator approximations in homogenization theory have attracted considerable
attention since then, with a number of interesting results~– see
\cite{Gri:2004}, \cite{Gri:2006}, \cite{Zh:2005}, \cite{ZhPas:2005},
\cite{Bor:2008}, \cite{KLS:2012}, \cite{Su:2013-1}, \cite{Su:2013-2},
\cite{ChC:2016} and~\cite{ZhPas:2016}, to name just a few.

As⋅each weakly convergent sequence of operators is bounded, one may
ask whether a sequence of the resolvents of⋅$-\Div A\rbr{x,\varepsilon^{-1}x}\hairspace\nabla$
converge in the uniform operator topology on⋅$L_{p}\rbr{\Omega}$
provided that it is bounded in the operator norm from⋅$W_{p}^{-1}\rbr{\Omega}$
to⋅$\mathring{W}_{p}^{1}\rbr{\Omega}$ for⋅$p$ other than~$2$.
Another question that naturally arises in this context is which domains
and boundary conditions are allowed to still yield the convergence
of the resolvent in the uniform operator topology on~$L_{p}\rbr{\Omega}$.
The⋅answer we give in this paper is somewhat implicit, for it is formulated
in terms of, e.g., boundary regularity results, but we provide some
 examples as well.

Let $\Omega$ be a uniformly  Lipschitz domain (possibly unbounded).
For⋅fixed⋅$p\in\rbr{1,\infty}$, let $\set W_{p}^{1}\rbr{\Omega;\C^{n}}$
be a subspace of⋅$W_{p}^{1}\rbr{\Omega}^{n}$ that contains⋅$\mathring{W}_{p}^{1}\rbr{\Omega}^{n}$,
and let $\set W_{\smash[t]{\cramped{p^{\dual}}}}^{1}\rbr{\Omega;\C^{n}}$ be defined similarly
for the exponent~$p^{\dual}$ conjugate to~$p$. Let $A_{kl}$ be
$\C^{n\times n}$\nobreakdash-matrix-valued mappings on⋅$\wb{\Omega}\times\R^{d}$
that are Lipschitz in the first variable and periodic (with⋅respect
to~$\Z^{d}$) in the second and set~$A=\cbr{A_{kl}}_{k,l=1}^{d}$.
We will study the matrix operator
\[
\op A^{\varepsilon}=-\Div A\rbr{x,\varepsilon^{-1}x}\hairspace\nabla
\]
acting between⋅$\set W_{p}^{1}\rbr{\Omega;\C^{n}}$ and⋅$\set W_{p}^{-1}\rbr{\Omega;\C^{n}}$,
the dual of⋅$\set W_{\smash[t]{\cramped{p^{\dual}}}}^{1}\rbr{\Omega;\C^{n}}$. We point
out that a function in⋅$\set W_{p}^{1}\rbr{\Omega;\C^{n}}$ may satisfy
mixed boundary conditions and even different components of this function
may satisfy different boundary conditions.

Suppose that, for some⋅$\mu\in\C$ and all sufficiently small⋅$\varepsilon$,
the operators~$\op A^{\varepsilon}-\mu$ are isomorphisms with uniformly
bounded (in~$\varepsilon$) inverses. This condition is obviously
necessary for the sequence~$\rbr{\op A^{\varepsilon}-\mu}^{-1}$
to have a limit even in the weak operator topology and thus  is not
much related to homogenization; the next two definitely are. We assume
that, for each⋅$x\in\Omega$, the cell problem
\[
-\Div A\rbr{x,\wc}\hairspace\rbr{\nabla N\rbr{x,\wc}+I}=0
\]
has a unique solution in⋅$W_{p}^{1}\rbr{\T^{d}}$ which is Lipschitz
in~$x$, and  the resolvent~$\rbr{\op A^{0}-\mu}^{-1}$ of the
effective operator~$\op A^{0}$ is continuous from⋅$L_{p}\rbr{\Omega}^{n}$
to⋅$\set W_{p}^{1}\rbr{\Omega}^{n}\cap\Lambda_{p}^{1+s}\rbr{\Omega}^{n}$
for some~$s\in(0,1]$, where $\Lambda_{p}^{1+s}\rbr{\Omega}$ is
the Lipschitz–Besov space~$B_{p,\infty}^{1+s}\rbr{\Omega}$.

The⋅basic examples are the Dirichlet and the Neumann problems for
strongly elliptic operators~$\op A^{\varepsilon}$ on a bounded $C^{1,1}$
domain. In⋅this case, there is a sector~$\set S$ in the complex
plain and an open neighborhood~$\set P$ of the exponent~$2$ such
that our assumptions hold for any⋅$\mu\notin\set S$ and~$p\in\set P$.
Moreover, $\set P=\rbr{1,\infty}$ as long as the function~$A$
belongs to the $\VMO$ space in the ``periodic'' variable. See⋅Section~\ref{sec: Examples}
for details.

In⋅this paper we  prove that
\begin{gather}
\norm{\rbr{\mathcal{A}^{\varepsilon}-\mu}^{-1}-\rbr{\mathcal{A}^{0}-\mu}^{-1}}_{L_{p}\rbr{\Omega}^{n}\to L_{p}\rbr{\Omega}^{n}}\le C\varepsilon^{s/p},\label{est: Intro | Convergence in Lp}\\
\norm{\nabla\rbr{\mathcal{A}^{\varepsilon}-\mu}^{-1}-\nabla\rbr{\mathcal{A}^{0}-\mu}^{-1}-\varepsilon\nabla\op K_{\mu}^{\varepsilon}}_{L_{p}\rbr{\Omega}^{n}\to L_{p}\rbr{\Omega}^{n}}\le C\varepsilon^{s/p},\label{est: Intro | Approximation in Wp1}
\end{gather}
where $\op K_{\mu}^{\varepsilon}$ is a so-called corrector, see Theorem~\ref{thm: Convergence and Approximation with 1st corrector}.
If, in addition, the adjoint of⋅$\op A^{\varepsilon}$ satisfies similar
assumptions as⋅$\op A^{\varepsilon}$ with some exponent~$s^{\dual}\in(0,s]$,
then
\begin{equation}
\norm{\rbr{\mathcal{A}^{\varepsilon}-\mu}^{-1}-\rbr{\mathcal{A}^{0}-\mu}^{-1}}_{L_{p}\rbr{\Omega}^{n}\to L_{p}\rbr{\Omega}^{n}}\le C\varepsilon^{s/p+s^{\dual}{\mkern-4mu }/p^{\dual}},\label{est: Intro | Convergence in Lp with e}
\end{equation}
see⋅Theorem~\ref{thm: Convergence with e}. For⋅$s=s^{\dual}=1$,
the convergence rate in⋅(\ref{est: Intro | Convergence in Lp with e})
is the same as in the whole space case, which is known to be sharp
with respect to the order. If⋅we have a uniform Caccioppoli-type inequality
for⋅$\op A^{\varepsilon}$, then the estimate~(\ref{est: Intro | Approximation in Wp1})
can be improved as well, but only away from the boundary. Thus, for
 a domain~$\Sigma\subset\Omega$ separated from the boundary
of⋅$\Omega$ by a positive distance,
\begin{equation}
\norm{\nabla\rbr{\mathcal{A}^{\varepsilon}-\mu}^{-1}-\nabla\rbr{\mathcal{A}^{0}-\mu}^{-1}-\varepsilon\nabla\op K_{\mu}^{\varepsilon}}_{L_{p}\rbr{\Omega}^{n}\to L_{p}\rbr{\Sigma}^{n}}\le C\varepsilon^{s/p+s^{\dual}{\mkern-4mu }/p^{\dual}},\label{est: Intro | Approximation in Wp1 with e}
\end{equation}
see⋅Corollary~\ref{cor: Interior approximation with 1st corrector}.
We mention that, in the whole space case, one can also find the second
term in the approximation~(\ref{est: Intro | Convergence in Lp with e})
so that the rate becomes of order~$\varepsilon^{2}$, see~\cite{BSu:2005}
and⋅\cite{Se:2017-1}, where the case⋅$p=2$ was handled.

Purely periodic homogenization problems on bounded domains have been
heavily investigated. By⋅using the unfolding method, Griso~\cite{Gri:2004},⋅\cite{Gri:2006}
 studied scalar problems  with Dirichlet, Neumann and mixed boundary
conditions within the Hilbert-space framework (i.e., when~$p=2$).
For⋅the Dirichlet and Neumann problems on $C^{1,1}$⋅domains, he obtained
the  estimates~(\ref{est: Intro | Convergence in Lp})–(\ref{est: Intro | Approximation in Wp1 with e})
with⋅$s=s^{\dual}=1$. For⋅mixed Dirichlet–Neumann problems on $C^{0,1}$⋅domains,
the approximations~(\ref{est: Intro | Convergence in Lp}),⋅(\ref{est: Intro | Approximation in Wp1})
were obtained with error of order~$\varepsilon^{\gamma}$ for some~$\gamma\in(0,1/3]$;
furthermore, it was shown for self-adjoint operators on polygonal
domains in⋅$\R^{2}$ or⋅$\R^{3}$ that (\ref{est: Intro | Convergence in Lp})–(\ref{est: Intro | Approximation in Wp1 with e})
hold for $s=s^{\dual}$ as long as $\rbr{\op A^{0}-\mu}^{-1}$ maps⋅$L_{2}\rbr{\Omega}$
to⋅$W_{2}^{1+s}\rbr{\Omega}$ continuously (the⋅proof may be carried
over to more general domains).  In⋅the⋅case when⋅$s=s^{\dual}=1$
and⋅$p=2$, Zhikov~\cite{Zh:2005} and Zhikov with⋅Pastukhova~\cite{ZhPas:2005}
(see also the survey paper~\cite{ZhPas:2016} and the references
therein) proved (\ref{est: Intro | Convergence in Lp})–(\ref{est: Intro | Approximation in Wp1})
 for scalar problems and the linear elasticity system on sufficiently
smooth domains with Dirichlet or Neumann boundary conditions. Fairly
general self-adjoint strongly elliptic systems on $C^{1,1}$⋅domains
with Dirichlet or Neumann boundary conditions were studied by Suslina⋅\cite{PSu:2012},⋅\cite{Su:2013-1},⋅\cite{Su:2013-2},
where the estimates~(\ref{est: Intro | Convergence in Lp})–(\ref{est: Intro | Approximation in Wp1 with e})
were proved for⋅$s=s^{\dual}=1$ and~$p=2$.  In⋅\cite{KLS:2012},
the authors considered self-adjoint Dirichlet and Neumann problems
on $C^{0,1}$⋅domains with H\"{o}lder continuous coefficients and,
for⋅$p=2$, obtained the approximation~(\ref{est: Intro | Convergence in Lp with e})
with error of order~$\varepsilon\abs{\ln\varepsilon}^{\sigma}$ for
any~$\sigma>1/2$. They also improved the rate to⋅$\varepsilon$
if⋅$s=s^{\dual}=1$.  Later, Shen~\cite{Shen:2018} proved that,
without any continuity assumptions on the coefficients,
\[
\norm{\rbr{\mathcal{A}^{\varepsilon}-\mu}^{-1}f-\rbr{\mathcal{A}^{0}-\mu}^{-1}f}_{L_{q^{\dual}}\rbr{\Omega}^{n}}\le C\varepsilon\norm{\rbr{\mathcal{A}^{0}-\mu}^{-1}f}_{W_{q}^{2}\rbr{\Omega}^{n}}
\]
with⋅$q=2d/\rbr{d+1}$. Moreover, if $\Omega$ is a $C^{1,1}$⋅domain
and the coefficients are H\"{o}lder continuous, then he showed that
\[
\norm{\rbr{\mathcal{A}^{\varepsilon}-\mu}^{-1}-\rbr{\mathcal{A}^{0}-\mu}^{-1}}_{L_{p}\rbr{\Omega}^{n}\to L_{q}\rbr{\Omega}^{n}}\le C\varkappa\rbr{\varepsilon},
\]
where⋅$q=pd/\rbr{p-d}$ for⋅$p\in\rbr{1,d}$, $q=\infty$ for⋅$p>d$
and $\varkappa\rbr{\varepsilon}=\varepsilon$ for the Dirichlet problem,
$\varkappa\rbr{\varepsilon}=\varepsilon\abs{\ln\varepsilon}$ for
the Neumann problem. We also mention the paper~\cite{Nazarov:2006}
of Nazarov, who treats, among other things, interior estimates for
the scalar self-adjoint Dirichlet problem with sufficiently smooth
coefficients. In⋅particular, if $\Omega$ is a $C^{5}$⋅domain and
$\Sigma\subset\Omega$ is separated from $\partial\Omega$ by a positive
distance, then
\[
\norm{\rbr{\mathcal{A}^{\varepsilon}-\mu}^{-1}-\rbr{\mathcal{A}^{0}-\mu}^{-1}}_{L_{p}\rbr{\Omega}\to L_{q}\rbr{\Sigma}}\le C\varepsilon,
\]
where⋅$q<pd/\rbr{p-d}$ for⋅$p\in\rbr{1,d}$ and⋅$q=\infty$ for⋅$p>d$.

To⋅prove our results, we study a first-order approximation, involving
the resolvents of the original and the effective operators and the
corrector. First-order approximations are well-known in homogenization
theory, see, e.g.,⋅\cite{BLP:1978} or⋅\cite{ZhKO:1993}. The⋅one
we use here differs from the classical one in that the corrector is
now regularized. The⋅idea of using a smoothing to regularize the
classical corrector is due to Cioranescu, Damlamian and Griso, see~\cite{CDG:2002}.
Besides the standard mollification, we employ the Steklov smoothing
operator, which is the most simple and which had already proved to
be quite useful for both linear and non-linear problems; see⋅\cite{Zh:2005}
and⋅\cite{ZhPas:2005}, where that smoothing first appeared in the
context of homogenization, as well as \cite{PasT:2007}, \cite{Su:2013-1}
and~\cite{Su:2013-2}. We adopt the technique related to the Steklov
smoothing operator from these papers. 

For⋅the⋅first-order approximation, we derive an operator representation
that splits the problem into interior and boundary parts, see~(\ref{eq: The identity}).
The⋅interior part is treated in the same way as for the whole space
case, cf.~\cite{Se:2017-1}. On⋅the other hand, the boundary part,
being supported in  a small neighborhood of the boundary, is small
 as well, no matter what the boundary conditions.

We note that, once the estimates~(\ref{est: Intro | Convergence in Lp})–(\ref{est: Intro | Approximation in Wp1 with e})
are obtained, a limiting argument will lead to similar results for
locally periodic operators whose coefficients are H\"{o}lder continuous
in the first variable, see⋅\cite{Se:2017-3} for some details. We
also mention the paper~\cite{Se:2020}, where the elliptic bounds~(\ref{est: Intro | Approximation in Wp1})–(\ref{est: Intro | Convergence in Lp with e})
for the Dirichlet problem with⋅$s=1$ and⋅$p=2$ were carried over
to the parabolic semigroup by keeping track of the rate dependence
on both the small parameter~$\varepsilon$ and spectral parameter~$\mu$.

\section{\label{sec: Notation}Notation}

The⋅symbol~$\norm{\wc}_{U}$ will stand for the norm on a normed
space~$U$. If⋅$U$ and⋅$V$ are Banach spaces, then $\B\rbr{U,V}$
is the Banach space of bounded linear operators from⋅$U$ to~$V$.
When⋅$U=V$, the space~$\B\rbr U=\B\rbr{U,U}$ becomes a Banach algebra
with the identity map~$\op I$. The⋅norm and the inner product on⋅$\C^{n}$
are denoted by⋅$\abs{\wc}$ and⋅$\abr{\wc,\wc}$. We shall often
identify $\B\rbr{\C^{n},\C^{m}}$ and~$\C^{m\times n}$.

Let⋅$\Sigma$ be a domain in⋅$\R^{d}$ and⋅$U$ a Banach space. Then
$L_{0}\rbr{\Sigma;U}$ is the vector space of all strongly measurable
functions on⋅$\Sigma$ with values in~$U$. The⋅symbol~$L_{p}\rbr{\Sigma;U}$,
$p\in\sbr{1,\infty}$, stands for the⋅$L_{p}$\nobreakdash-space
of  $L_{0}\rbr{\Sigma;U}$\nobreakdash-functions.  For⋅$p<\infty$
and⋅$s\ge0$, we let $W_{p}^{s}\rbr{\Sigma;U}$ denote the usual Sobolev
space or Sobolev–Slobodetskii space of $L_{0}\rbr{\Sigma;U}$\nobreakdash-functions
with norm
\[
\norm u_{W_{p}^{s}\rbr{\Sigma;U}}=\biggrbr{\sum_{\abs{\alpha}=0}^{m}\norm{D^{\alpha}u}_{L_{p}\rbr{\Sigma;U}}^{p}}^{1/p}
\]
if⋅$s=m\in\N_{0}$~and
\[
\norm u_{W_{p}^{s}\rbr{\Sigma;U}}=\biggrbr{\sum_{\abs{\alpha}=m}\Seminorm{D^{\alpha}u}_{W_{p}^{r}\rbr{\Sigma;U}}^{p}+\norm u_{W_{p}^{m}\rbr{\Sigma;U}}^{p}}^{1/p}
\]
if⋅$s=m+r$ with $m\in\N_{0}$ and~$r\in\rbr{0,1}$. Here $D=-i\nabla$~and
\[
\Seminorm u_{W_{p}^{r}\rbr{\Sigma;U}}=\biggrbr{\int_{0}^{\infty}t^{-1-rp}\osc_{p}\rbr{u;t}^{p}\dd t}^{1/p},
\]
where $\osc_{p}\rbr{u;\wc}$ is the⋅$L_{p}$\nobreakdash-modulus
of continuity of⋅$u$, given⋅by
\[
\osc_{p}\rbr{u;t}=\sup_{\abs h\le t}\norm{u\rbr{\wc+h}-u\rbr{\wc}}_{L_{p}\rbr{\Sigma\cap\rbr{-h+\Sigma};U}}.
\]
In⋅the⋅case⋅$U=\C^{n}$, we write⋅$\norm{\wc}_{p,\Sigma}$ and⋅$\norm{\wc}_{s,p,\Sigma}$
for the norms on⋅$L_{p}\rbr{\Sigma}^{n}=L_{p}\rbr{\Sigma;U}$ and⋅$W_{p}^{s}\rbr{\Sigma}^{n}=W_{p}^{s}\rbr{\Sigma;U}$,
respectively, and⋅$\rbr{\wc,\wc}_{\Sigma}$ for the inner product
on~$L_{2}\rbr{\Sigma}^{n}$. The⋅dual space of⋅$W_{p}^{s}\rbr{\Sigma}$
under the pairing~$\rbr{\wc,\wc}_{\Sigma}$ is denoted by⋅$W_{p}^{s}\rbr{\Sigma}^{*}$,
with $\norm{\wc}_{-s,p^{\dual},\Sigma}^{*}$ standing for the norm.
Here $p^{\dual}$ is the exponent conjugate to⋅$p$, that is,⋅$1/p^{\dual}=1-1/p$.
The⋅closure of⋅$C_{c}^{\infty}\rbr{\Sigma}$ in⋅$W_{p}^{s}\rbr{\Sigma}$
is⋅$\mathring{W}_{p}^{s}\rbr{\Sigma}$, and $W_{\smash[t]{\cramped{p^{\dual}}}}^{-s}\rbr{\Sigma}$
is its dual, with norm~$\norm{\wc}_{-s,p^{\dual},\Sigma}$. As⋅usual,
$H^{s}=W_{2}^{s}$, $H^{-s}=W_{2}^{-s}$,~etc.

Let⋅$\cell$ be the closed cube in⋅$\R^{d}$ with center~$0$ and
side length~$1$, sides being parallel to the axes. Then  $\wt W_{p}^{m}\rbr{\cell}$
denotes the completion of⋅$\wt C^{m}\rbr{\cell}$ in the $W_{p}^{m}$\nobreakdash-norm.
Here $\wt C^{m}\rbr{\cell}$ is the class of $m$\nobreakdash-times
continuously differentiable functions  on⋅$\cell$ whose periodic
extension to⋅$\R^{d}$ enjoys the same smoothness. Notice that $\wt L_{p}\rbr{\cell}$
can be identified with the space of all periodic functions in~$L_{p,\text{loc}}\rbr{\R^{d}}$.
In⋅a⋅similar fashion, we define⋅$\wt W_{p}^{m}\rbr{\R^{d}\times\cell}$.
The⋅dual of⋅$\wt W_{p}^{m}$ is denoted by~$\wt W_{\smash[t]{\cramped{p^{\dual}}}}^{-m}$.

For⋅$p\in\sbr{1,\infty}$ and⋅$s=m+r$ with $m\in\N_{0}$ and⋅$r\in(0,1]$,
we also introduce the space~$\Lambda_{p}^{s}\rbr{\Sigma;U}$ of functions~$u\in W_{p}^{m}\rbr{\Sigma;U}$
with finite norm
\[
\abs u_{\Lambda_{p}^{s}\rbr{\Sigma;U}}=\biggrbr{\sum_{\abs{\alpha}=m}\seminorm{D^{\alpha}u}_{\Lambda_{p}^{r}\rbr{\Sigma;U}}^{p}+\norm u_{W_{p}^{m}\rbr{\Sigma;U}}^{p}}^{1/p},
\]
where
\[
\seminorm u_{\Lambda_{p}^{r}\rbr{\Sigma;U}}=\sup_{t\in\R_{+}}\abs t^{-r}\osc_{p}\rbr{u;t}.
\]
We use $\abs{\wc}_{s,p,\Sigma}$ and⋅$\seminorm{\wc}_{s,p,\Sigma}$
as shorthand for the norm and the seminorm on⋅$\Lambda_{p}^{s}\rbr{\Sigma}^{n}$.

Let $B$ be the open unit ball in⋅$\R^{d}$ centered at the origin,
and let $B_{+}$ be the open unit half-ball with~$x_{d}\in\rbr{0,1}$.
We say that a domain~$\Sigma$ in⋅$\R^{d}$ satisfies the uniform
weak Lipschitz condition if it is a $d$\nobreakdash-dimensional
Lipschitz manifold with boundary embedded into~$\R^{d}$. More precisely,
there is a uniformly locally-finite open covering~$\cbr{W_{k}}$
of⋅$\partial\Sigma$ and a sequence of bi-Lipschitz transformations~$\omega_{k}\colon W_{k}\to B$
so that (1)~$\omega_{k}\rbr{W_{k}\cap\Sigma}=B_{+}$ and $\omega_{k}\rbr{W_{k}\cap\partial\Sigma}=\partial B_{+}\!\setminus\partial B$;
(2)~$\sup_{k}\seminorm{\omega_{k}}_{1,\infty,W_{k}}$ and⋅$\sup_{k}\seminorm{\omega_{k}^{-1}}_{1,\infty,B}$
are finite; and (3)~for some⋅$\delta>0$, any open ball~$B_{\delta}\rbr x$
with center~$x\in\partial\Sigma$ and radius~$\delta$ is contained
in some coordinate patch~$W_{k}$. The⋅last
two conditions are automatically satisfied provided that the boundary
of⋅$\Sigma$ is compact. Notice that the domain~$\R^{d}\setminus\wb{\Sigma}$
is uniformly weakly Lipschitz whenever $\Sigma$~is.

For⋅such~$\Sigma$, there exists a $C^{\infty}$\nobreakdash-partition
of unity~$\cbr{\varphi_{k}}$ subordinate to⋅$\cbr{W_{k}}$ with
the property that $\sup_{k}\norm{D^{\alpha}\varphi_{k}}_{\infty,W_{k}}$
is finite for any~$\alpha$, see~\cite[Chapter~6, Section~3]{St:1970}.
Then there is a linear strong $m$\nobreakdash-extension operator~$\op E_{m}$,
which maps $W_{p}^{k}\rbr{\Sigma}$ and⋅$C^{k}\rbr{\wb{\Sigma}}$
into, respectively, $W_{p}^{k}\rbr{\R^{d}}$ and⋅$C^{k}\rbr{\wb{\R}^{d}}$
for every integer~$k\in\sbr{0,m}$ and every⋅$p\in[1,\infty)$, see~\cite{AF:2003}.
It⋅follows that standard density and embedding results which hold
for⋅$\R^{d}$ must also hold for~$\Sigma$. In⋅particular, the
Sobolev  theorem states that $W_{p}^{1}\rbr{\Sigma}$ is continuously
embedded in⋅$L_{q}\rbr{\Sigma}$ for any⋅$q\in\sbr{p,p^{*}}$. Here
$p^{*}$ is the Sobolev conjugate to $p$ given by⋅$1/p^{*}=1/p-1/d$
if⋅$p<d$; $p^{*}$ is any finite number greater than or equal to⋅$p$
if⋅$p=d$; and $p^{*}=\infty$ if~$p>d$. By⋅$p_{*}$ we denote an
exponent such that⋅$W_{p_{*}}^{1}\rbr{\Sigma}$ is embedded in⋅$L_{p}\rbr{\Sigma}$,
that is, $p_{*}=1$ if⋅$p\in[1,d^{\dual})$, $p_{*}$ satisfies⋅$1/p_{*}=1/p+1/d$
if⋅$p\in[d^{\dual},\infty)$ and⋅$p_{*}$ is any number greater than⋅$d$
if~$p=\infty$. We also note that $\Lambda_{p}^{s}\rbr{\Sigma;U}$
coincides with the Sobolev space~$W_{p}^{s}\rbr{\Sigma;U}$ if⋅$s\in\N$
and with the Lipschitz–Besov space~$B_{p,\infty}^{s}\rbr{\Sigma;U}$
if $s\notin\N$ (see, e.g.,~\cite{AF:2003} and⋅\cite{Agranovich:2013}
for the definition); in particular, $\Lambda_{\infty}^{s}\rbr{\Sigma;U}$,
with $s=m+r$ as above, is the space~$C^{m,r}\rbr{\wb{\Sigma};U}$
of H\"{o}lder continuous functions. By⋅interpolation, $\op E_{m}$
is then an extension operator for⋅$\Lambda_{p}^{s}\rbr{\Sigma;U}$
with any real~$s\in\sbr{0,m}$ and any⋅$p\in\sbr{1,\infty}$.

We shall also need the $\BMO\rbr{\R^{d}}$ and⋅$\VMO\rbr{\R^{d}}$
spaces. The⋅former consists of all $u\in L_{1,\text{loc}}\rbr{\R^{d}}$
such that
\[
\norm u_{\BMO}=\sup_{B_{R}}\dashint_{B_{R}}\abs{u\rbr x-m_{B_{R}}\rbr u}\dd x<\infty,
\]
where $B_{R}\subset\R^{d}$ is a ball of radius~$R$ and $m_{B_{R}}\rbr u=\dashint_{B_{R}}u\rbr y\dd y$
is the mean value of $u$ over~$B_{R}$. The⋅latter is the subspace
in⋅$\BMO\rbr{\R^{d}}$ of all functions~$u$ for which the $\VMO$\nobreakdash-modulus,
given~by
\[
\eta\rbr{u;r}=\sup_{B_{R}\colon R<r}\dashint_{B_{R}}\abs{u\rbr x-m_{B_{R}}\rbr u}\dd x,
\]
tends to zero as~$r\to0$. We refer the reader to⋅\cite{Gra:2014-2}
and⋅\cite{Gar:2007} for more on this matter.

For⋅a⋅set~$\Sigma\subset\R^{d}$, we let $\Sigma_{\delta}$ denote
a neighborhood of~$\Sigma$:
\[
\Sigma_{\delta}=\bigcbr{x\in\R^{d}\colon\dist\rbr{x,\Sigma}<\cellradius\delta},
\]
where $2\cellradius=\diam\cell=d^{1/2}$. Thus, $\Sigma+\delta\cell\subset\Sigma_{\delta}$
for any~$\delta\ge0$.

We will often use the notation~$\alpha\lesssim\beta$  (which is
the same as saying that~$\beta\gtrsim\alpha$) to mean that there
is a positive constant~$C$ depending only on some fixed parameters
(listed in Theorem~\ref{thm: Convergence and Approximation with 1st corrector}–Corollary~\ref{cor: Interior convergence for fractional derivative})
such that~$\alpha\le C\beta$. Finally, $\alpha\meet\beta$ and⋅$\alpha\join\beta$
are, respectively, the smaller and the larger of $\alpha$ and~$\beta$.

\section{\label{sec: Original operator}Original operator}

Let⋅$\Omega\subset\R^{d}$ be a (possibly unbounded) domain  satisfying
the uniform weak Lipschitz condition. Define the  operation~$\tau^{\varepsilon}$,⋅$\varepsilon>0$,
that takes a function~$u\colon\Omega\times\R^{d}\to L_{0}\rbr{\cell}$
to the function~$\tau^{\varepsilon}u\colon\Omega\to L_{0}\rbr{\cell}$
given by~$\tau^{\varepsilon}u\rbr x=u\rbr{x,\varepsilon^{-1}x}$.
We do not distinguish between⋅$\rbr{u\rbr{x,y}}\rbr z$ and⋅$u\rbr{x,y,z}$,
etc., and so for all⋅$x\in\Omega$ and⋅$z\in\cell$
\begin{equation}
\tau^{\varepsilon}u\rbr{x,z}=u\rbr{x,\varepsilon^{-1}x,z}.\label{def: Original operator | =0003C4=001D4B}
\end{equation}
Obviously, $\tau^{\varepsilon}$ is a homomorphism of the respective
algebras; in other words, for any two functions~$u$ and⋅$v$ from⋅$\Omega\times\R^{d}$
to~$L_{0}\rbr{\cell}$
\begin{equation}
\tau^{\varepsilon}\rbr{u+v}=\tau^{\varepsilon}u+\tau^{\varepsilon}v,\qquad\tau^{\varepsilon}uv=\tau^{\varepsilon}u\cdot\tau^{\varepsilon}v.\label{eq: Original operator | =0003C4=001D4B is multiplicative}
\end{equation}
Here ``$\cdot$'' denotes the pointwise product of functions, and
we adopt the convention that an operator acts on everything to its
right, until the end of the expression or the ``$\cdot$'' sign
is reached, i.e., $\tau^{\varepsilon}uv=\tau^{\varepsilon}\rbr{uv}$
and $\tau^{\varepsilon}u\cdot v=\rbr{\tau^{\varepsilon}u}v$. We will
use the notation~$u^{\varepsilon}=\tau^{\varepsilon}u$.

Let⋅$A_{kl}$, with $1\le k,l\le d$, be functions in~$C^{0,1}\rbr{\wb{\Omega};\wt L_{\infty}\rbr{\cell}}^{n\times n}$.
Then $A=\cbr{A_{kl}}$ can be thought of as a bounded mapping⋅$A\colon\wb{\Omega}\times\R^{d}\to\B\rbr{\C^{d\times n}}$
which is Lipschitz in the first variable and periodic in the second.
It⋅follows that $A$ satisfies a Carath\'{e}odory-type condition,
i.e., $A\rbr{\wc,y}$ is  continuous on⋅$\wb{\Omega}$ for almost
every⋅$y\in\cell$ uniformly with respect to⋅$y$ and $A\rbr{x,\wc}$
is  measurable on⋅$\R^{d}$ for each⋅$x\in\wb{\Omega}$ (see the
proof of⋅Lemma~5.6 in~\cite{Al:1992}). Therefore, $A^{\varepsilon}$
is measurable and uniformly bounded.

Fix⋅$p\in\rbr{1,\infty}$. Let $\set W_{p}^{1}\rbr{\Omega;\C^{n}}$
and⋅$\set W_{\smash[t]{\cramped{p^{\dual}}}}^{1}\rbr{\Omega;\C^{n}}$ be subspaces of,
respectively, $W_{p}^{1}\rbr{\Omega}^{n}$ and⋅$W_{\smash[t]{\cramped{p^{\dual}}}}^{1}\rbr{\Omega}^{n}$
that contain all functions in⋅$C_{c}^{\infty}\rbr{\Omega}^{n}$; thus
\begin{equation}
\mathring{W}_{p}^{1}\rbr{\Omega}^{n}\subseteq\set W_{p}^{1}\rbr{\Omega;\C^{n}}\subseteq W_{p}^{1}\rbr{\Omega}^{n}.\label{def: Original operator | Wp1}
\end{equation}
By⋅$\set W_{p}^{-1}\rbr{\Omega;\C^{n}}$ and⋅$\Norm{\wc}_{-1,p,\Omega}$
we denote the dual of⋅$\set W_{\smash[t]{\cramped{p^{\dual}}}}^{1}\rbr{\Omega;\C^{n}}$
(under the $L_{2}$\nobreakdash-pair\-ing) and the associated norm.
Since $\set W_{p}^{-1}\rbr{\Omega;\C^{n}}$ is isometrically isomorphic
to the quotient space⋅$\rbr{W_{\smash[t]{\cramped{p^{\dual}}}}^{1}\rbr{\Omega}^{n}}^{*}\!/\rbr{\set W_{\smash[t]{\cramped{p^{\dual}}}}^{1}\rbr{\Omega;\C^{n}}}^{\bot}$,
where $\rbr{\set W_{\smash[t]{\cramped{p^{\dual}}}}^{1}\rbr{\Omega;\C^{n}}}^{\bot}$ is
the subspace of all functionals in⋅$\rbr{W_{\smash[t]{\cramped{p^{\dual}}}}^{1}\rbr{\Omega}^{n}}^{*}$
vanishing on⋅$\set W_{\smash[t]{\cramped{p^{\dual}}}}^{1}\rbr{\Omega;\C^{n}}$, the natural
projection
\begin{equation}
\pi\colon f\mapsto f+\rbr{\set W_{\smash[t]{\cramped{p^{\dual}}}}^{1}\rbr{\Omega;\C^{n}}}^{\bot}\label{def: Quotient map}
\end{equation}
can be thought of as  an epimorphism of⋅$\rbr{W_{\smash[t]{\cramped{p^{\dual}}}}^{1}\rbr{\Omega}^{n}}^{*}$
onto⋅$\set W_{p}^{-1}\rbr{\Omega;\C^{n}}$:
\begin{equation}
\Norm{\pi f}_{-1,p,\Omega}\le\norm f_{-1,p,\Omega}^{*}.\label{est: Norm of pi}
\end{equation}

Consider the matrix operator~$\op A^{\varepsilon}\colon\set W_{p}^{1}\rbr{\Omega;\C^{n}}\to\set W_{p}^{-1}\rbr{\Omega;\C^{n}}$
given⋅by
\begin{equation}
\op A^{\varepsilon}=D^{*}A^{\varepsilon}D,\label{def: Original operator | A=001D4B}
\end{equation}
that is, $\op A^{\varepsilon}$ sends each⋅$u\in\set W_{p}^{1}\rbr{\Omega;\C^{n}}$
to the functional~$v\mapsto\rbr{A^{\varepsilon}Du,Dv}_{\Omega}$
belonging to~$\set W_{p}^{-1}\rbr{\Omega;\C^{n}}$. It⋅is plain that
$\op A^{\varepsilon}$ is bounded uniformly with respect to~$\varepsilon$:
\begin{equation}
\Norm{\op A^{\varepsilon}u}_{-1,p,\Omega}\le\norm A_{L_{\infty}}\norm{Du}_{p,\Omega},\qquad u\in\set W_{p}^{1}\rbr{\Omega;\C^{n}}.\label{est: Original operator | A=001D4B is bounded}
\end{equation}
We further assume that, for some⋅$\mu\in\C$, there is $\varepsilon_{\mu}\in(0,1]$
so that the operators~$\smash[b]{\op A_{\mu}^{\varepsilon}}=\op A^{\varepsilon}-\mu$
are isomorphisms for any⋅$\varepsilon\in\set E_{\mu}=(0,\varepsilon_{\mu}]$
and, moreover, have uniformly bounded inverses:
\begin{equation}
\norm{\rbr{\op A_{\mu}^{\varepsilon}}^{-1}f}_{1,p,\Omega}\lesssim\Norm f_{-1,p,\Omega},\qquad f\in\set W_{p}^{-1}\rbr{\Omega;\C^{n}}.\label{est: Original operator | Norm of (A=001D4B-=0003BC)=00207B=0000B9}
\end{equation}

Let $\smash[b]{\rbr{\op A_{\mu}^{\varepsilon}}^{\dual}}$ be the adjont
of~$\smash[b]{\op A_{\mu}^{\varepsilon}}$. The⋅corresponding objects
and results related to⋅$\smash[b]{\rbr{\op A_{\mu}^{\varepsilon}}^{\dual}}$,
will be marked with~``$\dual$'' too. Notice that  $\smash[b]{\rbr{\op A_{\mu}^{\varepsilon}}^{\dual}}$
obeys (\dualref{est: Original operator | A=001D4B is bounded})
and~(\dualref{est: Original operator | Norm of (A=001D4B-=0003BC)=00207B=0000B9}),
with the same constants, in fact.
\begin{rem}
Basic examples to keep in mind are the extreme cases where the space~$\set W_{p}^{1}\rbr{\Omega;\C^{n}}$
coincides with either $\mathring{W}_{p}^{1}\rbr{\Omega}^{n}$ or~$W_{p}^{1}\rbr{\Omega}^{n}$.
The⋅first case corresponds to the homogeneous Dirichlet problem, and
the second, to the homogeneous Neumann problem. Notice that components
of a function in⋅$\set W_{p}^{1}\rbr{\Omega;\C^{n}}$ may satisfy
different boundary conditions. We also point out that the space~$\set W_{p\meet p^{\dual}}^{1}\rbr{\Omega;\C^{n}}$
need not be locally embedded into⋅$\set W_{p\join p^{\dual}}^{1}\rbr{\Omega;\C^{n}}$,
although it usually does in applications.
\end{rem}
The⋅following result will be useful for interior estimates.
\begin{lem}
\label{lem: Original operator | Interior energy estimate}Suppose
that the inverse of⋅$\op A_{\mu}^{\varepsilon}$ is also bounded from⋅$\rbr{W_{q^{\dual}}^{1}\rbr{\Omega}^{n}}^{*}$
to⋅$W_{q}^{1}\rbr{\Omega}^{n}$ for some⋅$q\in[1,p)$. Assume further
that, given any⋅$\chi\in C_{c}^{0,1}\rbr{\Omega}$, there is $\chi^{\prime}\in C_{c}^{0,1}\rbr{\Omega}$,
with $\supp\chi\subset\supp\chi^{\prime}$, so that
\begin{equation}
\norm{D\chi u}_{q,\supp\chi}\lesssim\norm u_{q,\supp\chi^{\prime}}+\norm{\chi^{\prime}\op A_{\mu}^{\varepsilon}u}_{-1,q,\Omega}^{*}\label{est: Original operator | Interior energy estimate}
\end{equation}
for⋅all⋅$u\in C_{c}^{1}\rbr{\Omega}^{n}$. Then the same conclusion
holds with⋅$p$ in place of~$q$.
\end{lem}
\begin{proof}
Fix $\chi\in C_{c}^{0,1}\rbr{\Omega}$ and choose a sequence of cutoff
functions~$\chi_{k}\in C_{c}^{0,1}\rbr{\Omega}$, where $0\le k\le m=\ceil{d\rbr{1/q-1/p}}$
(here $\ceil{\wc}$ is the ceiling function), in such a way that $\chi_{0}=\chi$,
$\supp\chi_{k}\subset\supp\chi_{k+1}$ and $\chi_{k+1}=1$ on~$\supp\chi_{k}$.
Let⋅$q_{0}=p$ and⋅$q_{k+1}=\rbr{q_{k}}_{*}$. We want to estimate
the $L_{q_{k}}$\nobreakdash-norm of⋅$D\chi_{k}u$ using the simple
identity
\[
\op A_{\mu}^{\varepsilon}\chi_{k}u=-\rbr{D\wb{\chi}_{k}}^{*}\cdot A^{\varepsilon}Du+D^{*}A^{\varepsilon}\rbr{D\chi_{k}\cdot u}+\chi_{k}\op A_{\mu}^{\varepsilon}u.
\]
The⋅space~$\rbr{W_{\smash[t]{\cramped{p^{\dual}}}}^{1}\rbr{\Omega}^{n}}^{*}$ is naturally
embedded into⋅$\set W_{p}^{-1}\rbr{\Omega;\C^{n}}$ via the projection~$\pi$,
and we see from⋅(\ref{est: Original operator | Norm of (A=001D4B-=0003BC)=00207B=0000B9})
that the inverse of⋅$\op A_{\mu}^{\varepsilon}$ is bounded from⋅$\rbr{W_{\smash[t]{\cramped{p^{\dual}}}}^{1}\rbr{\Omega}^{n}}^{*}$
to⋅$W_{p}^{1}\rbr{\Omega}^{n}$, and hence  from⋅$\rbr{W_{\smash[t]{\cramped{q_{\smash[t]{k}}^{\dual}}}}^{1}\rbr{\Omega}^{n}}^{*}$
to⋅$W_{q_{k}}^{1}\rbr{\Omega}^{n}$ for all⋅$q_{k}\ge q$, by interpolation.
Observe also that $\rbr{q_{k}^{\dual}}^{*}=q_{k+1}^{\dual}$ and
therefore⋅$L_{q_{k+1}}\rbr{\Omega}\subset W_{\smash[t]{\cramped{q_{\smash[t]{k}}^{\dual}}}}^{1}\rbr{\Omega}^{*}$.
Then, by H\"{o}lder's inequality,
\[
\norm{D\chi_{k}u}_{q_{k},\supp\chi_{k}}\lesssim\norm{Du}_{q_{k+1}\join q,\supp\chi_{k}}+\norm u_{p,\supp\chi_{k}}+\norm{\chi_{k}\op A_{\mu}^{\varepsilon}u}_{-1,p,\Omega}^{*}.
\]
Iterating this and using the fact that $\chi_{k+1}=1$ on $\supp\chi_{k}$,
we obtain
\begin{equation}
\norm{D\chi u}_{p,\supp\chi}\lesssim\norm{Du}_{q_{m}\join q,\supp\chi_{m-1}}+\norm u_{p,\supp\chi_{m-1}}+\norm{\chi_{m-1}\op A_{\mu}^{\varepsilon}u}_{-1,p,\Omega}^{*}.\label{est: Original operator | Interior energy estimate, iterating}
\end{equation}
Now note that $q_{m}\le q$, so the hypothesis and H\"{o}lder's
inequality show that
\[
\norm{Du}_{q_{m}\join q,\supp\chi_{m-1}}\lesssim\norm u_{p,\supp\chi_{m}^{\prime}}+\norm{\chi_{m}^{\prime}\op A_{\mu}^{\varepsilon}u}_{-1,p,\Omega}^{*}.
\]
Substituting this to⋅(\ref{est: Original operator | Interior energy estimate, iterating})
 gives the desired estimate with any⋅$\chi^{\prime}\in C_{c}^{0,1}\rbr{\Omega}$
which is~$1$ on~$\supp\chi_{m}^{\prime}$.
\end{proof}

\section{\label{sec: Effective operator}Effective operator}

As⋅usual, the coefficients of the effective operator are described
by the solution of the so-called cell problem.  This problem involves
two variables, the fast and the slow, and so  it is convenient to
introduce the notation $D_{1}$ and⋅$D_{2}$ to denote differentiation
in the first and second variable, respectively.   For⋅any⋅$x\in\Omega$,
define the operator⋅$\op A\rbr x=D_{2}^{*}A\rbr{x,\wc}\hairspace D_{2}\colon\wt W_{p}^{1}\rbr{\cell}^{n}\to\wt W_{p}^{-1}\rbr{\cell}^{n}$.
Then the cell problem is as follows: for each⋅$\xi\in\C^{d\times n}$
and⋅$x\in\Omega$, find⋅$N_{\xi}\rbr{x,\wc}\in\wt W_{p}^{1}\rbr{\cell}^{n}$
with $\int_{\cell}N_{\xi}\rbr{x,y}\dd y=0$, satisfying
\begin{equation}
\op A\rbr x\hairspace N_{\xi}\rbr{x,\wc}=-D_{2}^{*}A\rbr{x,\wc}\hairspace\xi.\label{def: N}
\end{equation}
If⋅such an⋅$N_{\xi}\rbr{x,\wc}$ exists, it is unique in⋅$\wt W_{p}^{1}\rbr{\cell}^{n}\!/\!\ker\op A\rbr x$,
and  we may choose $N_{\xi}\rbr{x,\wc}$ to be linear in⋅$\xi$.
Then the map assigning⋅$N_{\xi}\rbr{x,\wc}$ to each⋅$\xi$ is simply
an operator of multiplication by a function, which we denote by⋅$N\rbr{x,\wc}$,
so that $N\rbr{x,\wc}\xi+\ker\op A\rbr x$ is the solution to~(\ref{def: N}).
We assume that $N\rbr{x,\wc}$ exists for all⋅$x\in\Omega$ and, moreover,
\begin{equation}
N\in C^{0,1}\rbr{\wb{\Omega};\wt W_{p}^{1}\rbr{\cell}}.\label{assu: N is Lipschitz}
\end{equation}
The⋅standard sufficient condition is this:
\begin{lem}
\label{lm: A(x)}Suppose that, for any⋅$x\in\Omega$, $\op A\rbr x$
is an isomorphism of⋅$\wt W_{p}^{1}\rbr{\cell}^{n}\!/\C$ onto⋅$\rbr{\wt W_{p}^{1}\rbr{\cell}^{n}\!/\C}^{*}\simeq\wt W_{p}^{-1}\rbr{\cell}^{n}\cap\C^{\bot}$
with uniformly bounded in⋅$x$ inverse. Then the problem~\textup{(\ref{def: N})}
has a unique solution~$N_{\xi}=N\xi$, where $N$ satisfies⋅$\int_{\cell}N\rbr{x,y}\dd y=0$
and~\textup{(\ref{assu: N is Lipschitz})}.
\end{lem}
\begin{proof}
By⋅assumption,
\[
N_{\xi}\rbr{x,\wc}+\C=-\op A\rbr x^{-1}D_{2}^{*}A\rbr{x,\wc}\hairspace\xi
\]
is⋅the⋅unique solution of⋅(\ref{def: N})⋅and
\[
\norm{D_{2}N_{\xi}\rbr{x,\wc}}_{p,\cell}\lesssim\norm{D_{2}^{*}A\rbr{x,\wc}\hairspace\xi}_{-1,p,\cell},
\]
and⋅therefore
\[
\norm{D_{2}N}_{L_{\infty}\rbr{\Omega;L_{p}\rbr{\cell}}}\lesssim\norm A_{L_{\infty}}.
\]

Let⋅$\op T_{h}$,⋅$h\in\R^{d}$, be the translation operator defined
by⋅$\op T_{h}u\rbr{x,y}=u\rbr{x+h,y}$, where⋅$u\in L_{0}\rbr{\R^{d}\times\R^{d}}$,
and let⋅$\findiff_{h}=\op T_{h}-\op I$. Obviously,
\begin{equation}
\findiff_{h}uv=\findiff_{h}u\cdot v+\op T_{h}u\cdot\findiff_{h}v\label{eq: Delta_h uv}
\end{equation}
for⋅any⋅$u,v\in L_{0}\rbr{\R^{d}\times\R^{d}}$. It follows that if
$x,x+h\in\Omega$, then
\[
\findiff_{h}N_{\xi}\rbr{x,\wc}+\C=-\op A\rbr{x+h}^{-1}D_{2}^{*}\rbr{\findiff_{h}A\rbr{x,\wc}\cdot\rbr{D_{2}N_{\xi}\rbr{x,\wc}+\xi}}.
\]
Hence,
\[
\norm{D_{2}\findiff_{h}N_{\xi}\rbr{x,\wc}}_{p,\cell}\lesssim\norm{D_{2}^{*}\rbr{\findiff_{h}A\rbr{x,\wc}\cdot\rbr{D_{2}N_{\xi}\rbr{x,\wc}+\xi}}}_{-1,p,\cell},
\]
and as a result
\[
\norm{D_{1}D_{2}N}_{L_{\infty}\rbr{\Omega;L_{p}\rbr{\cell}}}\lesssim\norm{D_{1}A}_{L_{\infty}}\norm{I+D_{2}N}_{L_{\infty}\rbr{\Omega;L_{p}\rbr{\cell}}}.
\]

We have verified that $D_{2}N\in C^{0,1}\rbr{\wb{\Omega};\wt L_{p}\rbr{\cell}}$.
It⋅is⋅then immediate from the Poincar\'{e} inequality that⋅$N\in C^{0,1}\rbr{\wb{\Omega};\wt L_{p}\rbr{\cell}}$.
\end{proof}

Now define the effective operator~$\op A^{0}\colon\set W_{p}^{1}\rbr{\Omega;\C^{n}}\to\set W_{p}^{-1}\rbr{\Omega;\C^{n}}$
by setting
\begin{equation}
\op A^{0}=D^{*}A^{0}D,\label{def: A=002070}
\end{equation}
where⋅$A^{0}\colon\wb{\Omega}\to\B\rbr{\C^{d\times n}}$ is given⋅by
\begin{equation}
A^{0}\rbr x=\int_{\cell}A\rbr{x,y}\hairspace\rbr{I+D_{2}N\rbr{x,y}}\dd y.\label{def: Coefficient A=002070}
\end{equation}
Notice that since⋅$A$ and⋅$D_{2}N$ are uniformly continuous in
the first variable, so is~$A^{0}$. In⋅fact, we have⋅$A^{0}\in C^{0,1}\rbr{\wb{\Omega}}$.
Indeed, an easy calculation shows that
\[
\norm{A^{0}}_{L_{\infty}}\le\norm A_{L_{\infty}}\norm{I+D_{2}N}_{L_{\infty}\rbr{\Omega;L_{p}\rbr{\cell}}}
\]
and
\[
\norm{D_{1}A^{0}}_{L_{\infty}}\le\norm A_{L_{\infty}}\norm{D_{1}D_{2}N}_{L_{\infty}\rbr{\Omega;L_{p}\rbr{\cell}}}+\norm{D_{1}A}_{L_{\infty}}\norm{I+D_{2}N}_{L_{\infty}\rbr{\Omega;L_{p}\rbr{\cell}}}.
\]
Thus, by⋅(\ref{assu: N is Lipschitz}), $\norm{A^{0}}_{C^{0,1}\rbr{\wb{\Omega}}}$
is finite. It⋅follows that $\op A^{0}$ is bounded.

We suppose that there is⋅$s\in(0,1]$ such that the operator~$\op A_{\mu}^{0}=\op A^{0}-\mu$,
with the same⋅$\mu$ as in⋅(\ref{est: Original operator | Norm of (A=001D4B-=0003BC)=00207B=0000B9}),
has a continuous inverse from⋅$L_{p}\rbr{\Omega}^{n}$ to⋅$\Lambda_{p}^{1+s}\rbr{\Omega}^{n}$:
\begin{equation}
\abs{\rbr{\op A_{\mu}^{0}}^{-1}f}_{1+s,p,\Omega}\lesssim\norm f_{p,\Omega},\qquad f\in L_{p}\rbr{\Omega}^{n}.\label{est: Norm of (A=002070-=0003BC)=00207B=0000B9}
\end{equation}

\begin{rem}
One usually starts with⋅$\op A_{\mu}^{0}$ being an isomorphism
of⋅$\set W_{p}^{1}\rbr{\Omega;\C^{n}}$ onto⋅$\set W_{p}^{-1}\rbr{\Omega;\C^{n}}$,
while additional regularity as⋅in~(\ref{est: Norm of (A=002070-=0003BC)=00207B=0000B9})
requires that both the domain and the boundary conditions be more
regular as well. For⋅the⋅Dirichlet or the Neumann problems on convex
or uniformly $C^{1,1}$\nobreakdash-regular domains, we have⋅$s=1$
(that is, $\rbr{\op A_{\mu}^{0}}^{-1}$ is bounded from⋅$L_{p}\rbr{\Omega}^{n}$
to⋅$W_{p}^{2}\rbr{\Omega}^{n}$), see, e.g.,⋅\cite[Chapter~3]{Grisvard:2011}
and⋅\cite[Chapter~4]{McL:2000};
the same holds under a weaker assumption that the Jacobi matrix
of each coordinate map~$\omega_{k}$ is a multiplier on the Sobolev
space~$W_{p}^{1}\rbr{W_{k}}$ with  multiplier norm uniformly bounded
in~$k$, see~\cite[Chapter~14]{MSh:2009}.  In⋅the case of  mixed Dirichlet–Neumann problems, one cannot
hope that $s$ will be ``too large'' even for very regular domains
and coefficients, as $s\le2/p-1/2$, $p<4$,  for the Laplacian
on a half-space with mixed boundary conditions, see~\cite{Shamir:1968}.
We refer the reader also to⋅\cite{Jerison+Kenig:1995}, \cite{Savare:1997}
and⋅\cite{Savare:1998} and references therein for more on this matter.
\end{rem}
If⋅(\dualref{assu: N is Lipschitz}) holds, then the operator~$\rbr{\op A^{0}}^{\dual}$
defined by⋅(\dualref{def: A=002070}) coincides with the adjoint
of⋅$\op A^{0}$.  However, the assumptions~(\ref{est: Norm of (A=002070-=0003BC)=00207B=0000B9})
and⋅(\dualref{est: Norm of (A=002070-=0003BC)=00207B=0000B9})
are in general separate. The⋅fact is that the order of smoothness
may depend on the Lebesgue exponent, as for problems with mixed boundary
conditions (see the previous remark). We do not always require that⋅$s^{\dual}>0$,
but if we do, then we  shrink⋅$s^{\dual}$ so~$s^{\dual}\in(0,s]$.

\section{\label{sec: Corrector}Corrector}

Fix an extension operator~$\op E$  mapping⋅$W_{p}^{1}\rbr{\Omega}$
and⋅$\Lambda_{p}^{1+s}\rbr{\Omega}$ continuously into, respectively,
$W_{p}^{1}\rbr{\R^{d}}$ and⋅$\Lambda_{p}^{1+s}\rbr{\R^{d}}$. We
also extend the function~$N$ to⋅$\R^{d}\times\cell$ in such a way
that⋅$N\in C^{0,1}\rbr{\wb{\R}^{d};\wt W_{p}^{1}\rbr{\cell}}$ (e.g.,
by doing a reflection in the boundary). Define the operator⋅$\op K_{\mu}\colon L_{p}\rbr{\Omega}^{n}\to\Lambda_{p}^{s}\rbr{\R^{d};\wt W_{p}^{1}\rbr{\cell}^{n}}$
 to⋅be
\begin{equation}
\op K_{\mu}=ND_{1}\op R_{\mu}^{0},\label{def: K}
\end{equation}
where $\op R_{\mu}^{0}\colon L_{p}\rbr{\Omega}^{n}\to\Lambda_{p}^{1+s}\rbr{\R^{d}}^{n}$
is given⋅by
\begin{equation}
\op R_{\mu}^{0}=\op E\rbr{\op A_{\mu}^{0}}^{-1}.\label{def: R0_mu}
\end{equation}
Both⋅$\op R_{\mu}^{0}$ and⋅$\op K_{\mu}$ are bounded, as we shall
now see.
\begin{lem}
\label{lem: R0_mu and K_mu are bounded}For⋅any⋅$f\in L_{p}\rbr{\Omega}^{n}$,
we have
\begin{align}
\abs{\op R_{\mu}^{0}f}_{1+s,p,\R^{d}} & \lesssim\norm f_{p,\Omega},\label{est: R0_mu is bounded}\\
\abs{D_{2}\op K_{\mu}f}_{s_{1},p,\R^{d}\times\cell}+\abs{\op K_{\mu}f}_{s_{1},p,\R^{d}\times\cell} & \lesssim\norm f_{p,\Omega}.\label{est: K is bounded}
\end{align}
Here the subscript~``$\,1\!$''  in ``$\,s_{1}\!$'' is for
``the first variable''; in other words, $\abs{\wc}_{s_{1},p,\R^{d}\times\cell}$
is the norm on⋅$\Lambda_{p}^{s}\rbr{\R^{d};L_{p}\rbr{\cell}}$.
\end{lem}
\begin{proof}
Since $\op E$ is continuous, (\ref{est: R0_mu is bounded}) follows
from the assumption~(\ref{est: Norm of (A=002070-=0003BC)=00207B=0000B9}).
To⋅prove⋅(\ref{est: K is bounded}), we estimate the seminorm of⋅$D_{2}\op K_{\mu}f$;
the other terms are handled similarly. Using⋅(\ref{eq: Delta_h uv})
and⋅(\ref{def: K}), we write
\[
\findiff_{h}D_{2}\op K_{\mu}f=\findiff_{h}D_{2}N\cdot D_{1}\op R_{\mu}^{0}f+\op T_{h}D_{2}N\cdot\findiff_{h}D_{1}\op R_{\mu}^{0}f.
\]
Then
\[
\begin{aligned}\seminorm{D_{2}\op K_{\mu}f}_{s_{1},p,\R^{d}\times\cell} & =\sup_{t\in\R_{+}}\abs t^{-s}\sup_{\abs h\le t}\norm{\findiff_{h}D_{2}\op K_{\mu}f}_{p,\R^{d}\times\cell}\\
 & \le\norm{D_{2}N}_{C^{0,1}\rbr{\wb{\R}^{d};L_{p}\rbr{\cell}}}\norm{D\op R_{\mu}^{0}f}_{p,\R^{d}}\\
 & \quad+\norm{D_{2}N}_{L_{\infty}\rbr{\R^{d};L_{p}\rbr{\cell}}}\abs{D\op R_{\mu}^{0}f}_{s,p,\R^{d}}.
\end{aligned}
\]
Now (\ref{est: R0_mu is bounded}) gives the desired bound.
\end{proof}
The⋅regularity of the operator~$\op R_{\mu}^{0}$ implies that the
image of⋅$\op K_{\mu}$ is contained in the space~$W_{p}^{1}\rbr{\R^{d};\wt W_{p}^{1}\rbr{\cell}^{n}}$
only if~$s=1$. For⋅the⋅other cases, we will use mollification to
regularize $\op R_{\mu}^{0}$ in~$\op K_{\mu}$.

\subsection{Mollification}

We fix a non-negative function~$J\in C_{c}^{\infty}\rbr{B_{1}\rbr 0}$
satisfying⋅$\int_{\R^{d}}J\rbr x\dd x=1$. For⋅$\delta>0$, let $\op J_{\delta}$
be the standard operator of mollification, that is, $\op J_{\delta}u=J_{\delta}*u$,
where⋅$J_{\delta}\rbr x=\delta^{-d}J\rbr{\delta^{-1}x}$. Obviously,
the operator~$\op J_{\delta}$ maps $\Lambda_{p}^{s}\rbr{\R^{d}}$
into⋅$W_{p}^{1}\rbr{\R^{d}}$, but for $s<1$ its norm blows up as~$\delta\to0$.
It⋅is⋅also known that $\op J_{\delta}$ converges, as⋅$\delta\to0$,
to⋅$\op I$ in the operator norm on~$L_{p}\rbr{\R^{d}}$. The⋅next
two lemmas provide the rates of blow-up and convergence, respectively.

\begin{lem}
\label{lem: J_delta}Let⋅$r\in(0,1]$ and⋅$q\in[1,\infty)$. Then
for any⋅$\delta>0$ and⋅$u\in C_{c}^{\infty}\rbr{\R^{d}}$,
\begin{equation}
\norm{D\op J_{\delta}u}_{q,\R^{d}}\lesssim\delta^{-\rbr{1-r}}\seminorm u_{r,q,\R^{d}}.\label{est: J_delta}
\end{equation}
\end{lem}
\begin{proof}
It⋅is easy to see that
\[
D\op J_{\delta}u\rbr x=-\int_{\R^{d}}DJ_{\delta}\rbr{x-\wh x}\hairspace\findiff_{x-\wh x}u\rbr{\wh x}\dd\wh x,
\]
where the integration is, in fact, running over~$B_{\delta}\rbr x$.
Then, by H\"{o}lder's inequality,
\[
\begin{aligned}\norm{D\op J_{\delta}u}_{q,\R^{d}} & \le\biggrbr{\int_{\R^{d}}\abs{DJ_{\delta}\rbr{\wh x}}^{q^{\dual}}\dd\wh x}^{1/q^{\dual}}\biggrbr{\int_{B_{\delta}\rbr 0}\norm{\findiff_{h}u}_{q,\R^{d}}^{q}\dd h}^{1/q}\\
 & \le\delta^{-1-d/q}\norm{DJ}_{q^{\dual},\R^{d}}\delta^{d/q+r}\abs{B_{1}\rbr 0}^{1/q}\seminorm u_{r,q,\R^{d}},
\end{aligned}
\]
which is what we wished to prove.
\end{proof}
\begin{lem}
\label{lem: J_delta-I}Let⋅$r\in(0,1]$ and⋅$q\in[1,\infty)$. Then
for any⋅$\delta>0$ and⋅$u\in C_{c}^{\infty}\rbr{\R^{d}}$,
\begin{equation}
\norm{\rbr{\op J_{\delta}-\op I}u}_{q,\R^{d}}\lesssim\delta^{r}\seminorm u_{r,q,\R^{d}}.\label{est: J_delta-I}
\end{equation}
\end{lem}
\begin{proof}
We write
\[
\rbr{\op J_{\delta}-\op I}u\rbr x=\int_{\R^{d}}J_{\delta}\rbr{\wh x}\hairspace\findiff_{-\wh x}u\rbr x\dd\wh x
\]
and then repeat the argument in the proof of Lemma~\ref{lem: J_delta}.
\end{proof}

Now for⋅$s\in\rbr{0,1}$, we want to construct mollified versions
of⋅$\op R_{\mu}^{0}$ and⋅$\op K_{\mu}$ and define the operators~$\op R_{\mu}^{0}\rbr{\delta}\colon L_{p}\rbr{\Omega}^{n}\to W_{p}^{2}\rbr{\R^{d}}^{n}$
and⋅$\op K_{\mu}\rbr{\delta}\colon L_{p}\rbr{\Omega}^{n}\to W_{p}^{1}\rbr{\R^{d};\wt W_{p}^{1}\rbr{\cell}^{n}}$⋅by
\begin{align}
\op R_{\mu}^{0}\rbr{\delta} & =\op J_{\delta}\op R_{\mu}^{0},\label{def: R0_mu(delta)}\\
\op K_{\mu}\rbr{\delta} & =ND_{1}\op R_{\mu}^{0}\rbr{\delta}.\label{def: K(delta)}
\end{align}
We agree to set⋅$\op R_{\mu}^{0}\rbr{\delta}=\op R_{\mu}^{0}$ for⋅$s=1$,
in which case also⋅$\op K_{\mu}\rbr{\delta}=\op K_{\mu}$. In⋅the
next two lemmas, we summarize how good these approximations are.
The⋅proofs are immediate from Lemmas~\ref{lem: J_delta} and⋅\ref{lem: J_delta-I},
the assumption~(\ref{assu: N is Lipschitz}) and the bound⋅(\ref{est: R0_mu is bounded})
(cf.~the⋅proof of Lemma~\ref{lem: R0_mu and K_mu are bounded}).
Below, we use $\norm{\wc}_{1_{1},p,\R^{d}\times\cell}$ and⋅$\norm{\wc}_{1_{2},p,\R^{d}\times\cell}$
to denote the norms on⋅$W_{p}^{1}\rbr{\R^{d};L_{p}\rbr{\cell}}$ and⋅$W_{p}^{1}\rbr{\cell;L_{p}\rbr{\R^{d}}}$,
respectively.)
\begin{lem}
\label{lem: R0_mu(delta) and K_mu(delta) are bounded}For⋅any⋅$\delta\in(0,1]$
and⋅$f\in L_{p}\rbr{\Omega}^{n}$, we have
\begin{align}
\delta^{1-s}\norm{\op R_{\mu}^{0}\rbr{\delta}\hairspace f}_{2,p,\R^{d}} & \lesssim\norm f_{p,\Omega},\label{est: R0_mu(delta) is bounded}\\
\delta^{1-s}\norm{D_{1}\op K_{\mu}\rbr{\delta}\hairspace f}_{1_{2},p,\R^{d}\times\cell}+\norm{\op K_{\mu}\rbr{\delta}\hairspace f}_{1_{2},p,\R^{d}\times\cell} & \lesssim\norm f_{p,\Omega}.\label{est: K(delta) is bounded}
\end{align}
\end{lem}
\begin{lem}
\label{lem: R0_mu(delta) and K_mu(delta) convergence rates}For⋅any⋅$\delta>0$
and⋅$f\in L_{p}\rbr{\Omega}^{n}$, we have
\begin{align}
\norm{\rbr{\op R_{\mu}^{0}\rbr{\delta}-\op R_{\mu}^{0}}f}_{1,p,\R^{d}} & \lesssim\delta^{s}\norm f_{p,\Omega},\label{est: R0_mu(delta) convergence rate}\\
\norm{\rbr{\op K_{\mu}\rbr{\delta}-\op K_{\mu}}f}_{1_{2},p,\R^{d}\times\cell} & \lesssim\delta^{s}\norm f_{p,\Omega}.\label{est: K(delta) convergence rate}
\end{align}
\end{lem}
Since we do not impose any extra  assumptions on the coefficients
of⋅$\op A^{\varepsilon}$, the function~$\tau^{\varepsilon}N$ may
fail to be measurable, and therefore the classical corrector~$\tau^{\varepsilon}\op K_{\mu}$~–
and even the mollified one,⋅$\tau^{\varepsilon}\op K_{\mu}\rbr{\delta}$,~–
may not  map⋅$L_{p}\rbr{\Omega}^{n}$  into⋅$L_{0}\rbr{\Omega}^{n}$.
We use the Steklov smoothing operator to further regularize~$\op K_{\mu}\rbr{\delta}$.

\subsection{Smoothing}

Let $\op T^{\varepsilon}\colon L_{0}\rbr{\R^{d}\times\cell}\to L_{0}\rbr{\R^{d}\times\cell;L_{0}\rbr{\cell}}$
be the translation operator
\begin{equation}
\op T^{\varepsilon}u\rbr{x,y,z}=u\rbr{x+\varepsilon z,y},\label{def: T=001D4B}
\end{equation}
where $\rbr{x,y}\in\R^{d}\times\cell$ and~$z\in\cell$.  Obviously,
$\op T^{\varepsilon}\rbr{u+v}=\op T^{\varepsilon}u+\op T^{\varepsilon}v$
and⋅$\op T^{\varepsilon}uv=\op T^{\varepsilon}u\cdot\op T^{\varepsilon}v$,
so $\op T^{\varepsilon}$ is an algebra homomorphism. Next, the formal
adjoint of⋅$\op T^{\varepsilon}$ with respect to the $L_{2}$\nobreakdash-pairing
is given by the formula
\[
\rbr{\op T^{\varepsilon}}^{*}u\rbr{x,y}=\int_{\cell}u\rbr{x-\varepsilon z,y,z}\dd z.
\]
 Then⋅the Steklov smoothing operator~$\op S^{\varepsilon}$ is
the restriction of⋅$\rbr{\op T^{\varepsilon}}^{*}$ to⋅$L_{1}\rbr{\R^{d}\times\cell}+L_{\infty}\rbr{\R^{d}\times\cell}$;
in⋅other words,
\begin{equation}
\op S^{\varepsilon}u\rbr{x,y}=\int_{\cell}\op T^{\varepsilon}u\rbr{x,y,z}\dd z.\label{def: S=001D4B}
\end{equation}
The~operator~$\op S^{\varepsilon}$ thus defined is formally self-adjoint.

Here we collect some well-known facts about $\op T^{\varepsilon}$
and~$\op S^{\varepsilon}$, cf.~\cite[Subsection~2.1]{ZhPas:2016}.
\begin{lem}
\label{lem: =0003C4=001D4BT=001D4B}For⋅any⋅$q\in[1,\infty)$ and⋅$\varepsilon>0$,
$\tau^{\varepsilon}\op T^{\varepsilon}$ is an isometry of⋅$\wt L_{q}\rbr{\R^{d}\times\cell}$
into~$L_{q}\rbr{\R^{d};L_{q}\rbr{\cell}}$.
\end{lem}
\begin{proof}
By⋅change of variable,
\[
\norm{\tau^{\varepsilon}\op T^{\varepsilon}u}_{q,\R^{d}\times\cell}^{q}=\int_{\R^{d}}\int_{\cell}\abs{u\rbr{x+\varepsilon z,\varepsilon^{-1}x}}^{q}\dd x\dd z=\int_{\R^{d}}\int_{\cell}\abs{u\rbr{x,\varepsilon^{-1}x-z}}^{q}\dd x\dd z.
\]
But since⋅$u$ is periodic in the⋅second variable, this equals~$\norm u_{q,\R^{d}\times\cell}^{q}$.
\end{proof}
A⋅related result for⋅$\op S^{\varepsilon}$ is immediate from H\"{o}lder's
inequality and⋅Lemma~\ref{lem: =0003C4=001D4BT=001D4B}.
\begin{lem}
\label{lem: =0003C4=001D4BS=001D4B}For⋅any⋅$q\in[1,\infty)$ and⋅$\varepsilon>0$,
$\tau^{\varepsilon}\op S^{\varepsilon}$ is a bounded operator from⋅$\wt L_{q}\rbr{\R^{d}\times\cell}$
to⋅$L_{q}\rbr{\R^{d}}$ of norm~$1$.
\end{lem}
\begin{rem}
\label{rem: 0-extension for =0003C4=001D4BT=001D4B and =0003C4=001D4BS=001D4B}Notice
that the previous two lemmas are also useful for periodic (in the
second variable) functions on domains of the form~$\Sigma_{\varepsilon}\times\cell$,
because these functions can be extended by zero to~$\R^{d}\times\cell$
and therefore
\begin{align*}
\norm{\tau^{\varepsilon}\op T^{\varepsilon}u}_{q,\Sigma\times\cell} & \le\norm u_{q,\Sigma_{\varepsilon}\times\cell},\\
\norm{\tau^{\varepsilon}\op S^{\varepsilon}u}_{q,\Sigma} & \le\norm u_{q,\Sigma_{\varepsilon}\times\cell}.
\end{align*}
We will use this without further comment.
\end{rem}
Both $\op T^{\varepsilon}$ and⋅$\op S^{\varepsilon}$ converge
to the identity operator in uniform operator topologies whenever
the domain is ``smoother'' than the codomain.
\begin{lem}
\label{lem: T=001D4B-I}Let⋅$\Sigma$ be a domain in⋅$\R^{d}$, and
let⋅$q\in[1,\infty)$. Then for any $\varepsilon>0$ and $u\in C_{c}^{\infty}\rbr{\R^{d}\times\cell}$
we⋅have
\begin{equation}
\norm{\rbr{\op T^{\varepsilon}-\op I}u}_{q,\Sigma\times\cell\times\cell}\lesssim\varepsilon\norm{D_{1}u}_{q,\Sigma_{\varepsilon}\times\cell}.\label{est: T=001D4B-I is of order =0003B5}
\end{equation}
\end{lem}
\begin{proof}
This follows easily from the formula
\[
u\rbr{x+\varepsilon z,y}-u\rbr{x,y}=\varepsilon i\int_{0}^{1}\abr{D_{1}u\rbr{x+\varepsilon tz,y},z}\dd t.\qedhere
\]
\end{proof}

The⋅next lemma comes from the previous one, together with H\"{o}lder's
inequality.
\begin{lem}
\label{lem: S=001D4B-I}Let⋅$\Sigma$ be a domain in⋅$\R^{d}$, and
let⋅$q\in[1,\infty)$. Then for any $\varepsilon>0$ and $u\in C_{c}^{\infty}\rbr{\R^{d}\times\cell}$
we⋅have
\begin{equation}
\norm{\rbr{\op S^{\varepsilon}-\op I}u}_{q,\Sigma\times\cell}\lesssim\varepsilon\norm{D_{1}u}_{q,\Sigma_{\varepsilon}\times\cell}.\label{est: S=001D4B-I is of order =0003B5}
\end{equation}
\end{lem}

\subsection{Corrector}

We define the  corrector~$\op K_{\mu}^{\varepsilon}\colon L_{p}\rbr{\Omega}^{n}\to W_{p}^{1}\rbr{\Omega}^{n}$⋅by
\begin{equation}
\op K_{\mu}^{\varepsilon}=\tau^{\varepsilon}\op S^{\varepsilon}\op K_{\mu}\rbr{\varepsilon}.\label{def: K=001D4B}
\end{equation}
Thanks to the smoothing~$\op S^{\varepsilon}$, it is bounded⋅with
\begin{equation}
\varepsilon\norm{D\op K_{\mu}^{\varepsilon}f}_{p,\Omega}+\norm{\op K_{\mu}^{\varepsilon}f}_{p,\Omega}\lesssim\norm f_{p,\Omega}.\label{est: Norm of K=001D4B}
\end{equation}
Indeed, taking into account that $\varepsilon D\tau^{\varepsilon}\op S^{\varepsilon}=\varepsilon\tau^{\varepsilon}\op S^{\varepsilon}D_{1}+\tau^{\varepsilon}\op S^{\varepsilon}D_{2}$
and using Lemma~\ref{lem: =0003C4=001D4BS=001D4B}, we see that
\[
\varepsilon\norm{D\op K_{\mu}^{\varepsilon}f}_{p,\Omega}+\norm{\op K_{\mu}^{\varepsilon}f}_{p,\Omega}\lesssim\varepsilon\norm{D_{1}\op K_{\mu}\rbr{\varepsilon}\hairspace f}_{p,\Omega_{\varepsilon}\times\cell}+\norm{\op K_{\mu}\rbr{\varepsilon}\hairspace f}_{1_{2},p,\Omega_{\varepsilon}\times\cell}.
\]
The⋅estimate~(\ref{est: Norm of K=001D4B}) then follows from~(\ref{est: K(delta) is bounded}).
 We also notice that Lemmas~\ref{lem: R0_mu(delta) and K_mu(delta) are bounded}
and⋅\ref{lem: =0003C4=001D4BS=001D4B}, together with the Sobolev
embedding theorem, imply that, for each fixed⋅$q\in\sbr{p,p^{*}}$,
\begin{equation}
\varepsilon^{1-s}\norm{\op K_{\mu}^{\varepsilon}f}_{q,\Omega}\lesssim\norm f_{p,\Omega}.\label{est: Norm of K=001D4B with q}
\end{equation}

\begin{rem}
The⋅operator~$\op K_{\mu}^{\varepsilon}$ can be written explicitly⋅as
\[
\op K_{\mu}^{\varepsilon}f\rbr x=\int_{\cell}N\rbr{x+\varepsilon z,\varepsilon^{-1}x}\hairspace\op J_{\varepsilon}\op ED\rbr{\op A_{\mu}^{0}}^{-1}f\rbr{x+\varepsilon z}\dd z.
\]
In⋅this form, it first appeared for⋅$s=1$ (in⋅which case $\op J_{\varepsilon}$
is dropped from~$\op K_{\mu}^{\varepsilon}$) in the paper~\cite{PasT:2007}.
\end{rem}

\section{\label{sec: Main results}Main results}

Now we formulate the main results of the  paper. The⋅first one deals
with approximation under minimal assumptions on the initial problem.
\begin{thm}
\label{thm: Convergence and Approximation with 1st corrector}If⋅\textup{(\ref{est: Original operator | Norm of (A=001D4B-=0003BC)=00207B=0000B9})},
\textup{(\ref{assu: N is Lipschitz})} and⋅\textup{(\ref{est: Norm of (A=002070-=0003BC)=00207B=0000B9})}
hold, then for any⋅$\varepsilon\in\set E_{\mu}$ and⋅$f\in L_{p}\rbr{\Omega}^{n}$
we⋅have
\begin{align}
\norm{\rbr{\op A_{\mu}^{\varepsilon}}^{-1}f-\rbr{\op A_{\mu}^{0}}^{-1}f}_{q,\Omega} & \lesssim\varepsilon^{s/p}\norm f_{p,\Omega},\label{est: Convergence}\\
\norm{D\rbr{\op A_{\mu}^{\varepsilon}}^{-1}f-D\rbr{\op A_{\mu}^{0}}^{-1}f-\varepsilon D\op K_{\mu}^{\varepsilon}f}_{p,\Omega} & \lesssim\varepsilon^{s/p}\norm f_{p,\Omega},\label{est: Approximation with 1st corrector}
\end{align}
where~$q\in\sbr{p,p^{*}}$. The~constants depend only on the parameters~$d$,
$s$, $p$, $q$, $n$, $\mu$, the domain~$\Omega$, the $C^{0,1}$\nobreakdash-norms
of⋅$A$ and⋅$N$ and the constants in the  bounds~\textup{(\ref{est: Original operator | Norm of (A=001D4B-=0003BC)=00207B=0000B9})}
and~\textup{(\ref{est: Norm of (A=002070-=0003BC)=00207B=0000B9})}.
\end{thm}
Notice that the inverse of⋅$\op A_{\mu}^{\varepsilon}$ actually does
converge in the operator norm from $L_{p}$ to $W_{p}^{r}$ with~$r<1$,
yet the rate may be not as good.
\begin{cor}
\label{cor: Corollary 1}Under⋅the⋅hypotheses of Theorem~\textup{\ref{thm: Convergence and Approximation with 1st corrector}},
for any⋅$r\in\rbr{0,1}$, $\varepsilon\in\set E_{\mu}$ and⋅$f\in L_{p}\rbr{\Omega}^{n}$
it holds⋅that
\begin{equation}
\Seminorm{\rbr{\op A_{\mu}^{\varepsilon}}^{-1}f-\rbr{\op A_{\mu}^{0}}^{-1}f}_{r,p,\Omega}\lesssim\varepsilon^{s/p\meet\rbr{1-r}}\norm f_{p,\Omega}.\label{est: Corollary 1}
\end{equation}
The⋅constant depends only on the parameters~$d$, $r$, $s$, $p$,
$n$, $\mu$, the domain~$\Omega$, the $C^{0,1}$\nobreakdash-norms
of⋅$A$ and⋅$N$ and the constants in the  bounds~\textup{(\ref{est: Original operator | Norm of (A=001D4B-=0003BC)=00207B=0000B9})}
and~\textup{(\ref{est: Norm of (A=002070-=0003BC)=00207B=0000B9})}.
\end{cor}
We can improve the estimate~(\ref{est: Convergence}) for⋅$q=p$
provided that the adjoint problem enjoys silimar regularity properties
as for the initial one.
\begin{thm}
\label{thm: Convergence with e}Suppose that \textup{(\ref{est: Original operator | Norm of (A=001D4B-=0003BC)=00207B=0000B9})},
\textup{(\ref{assu: N is Lipschitz})}, \textup{(\ref{est: Norm of (A=002070-=0003BC)=00207B=0000B9})}
and \textup{(\dualref{assu: N is Lipschitz})}, \textup{(\dualref{est: Norm of (A=002070-=0003BC)=00207B=0000B9})}
hold. Then for any⋅$\varepsilon\in\set E_{\mu}$ and⋅$f\in L_{p}\rbr{\Omega}^{n}$
we⋅have
\begin{equation}
\norm{\rbr{\op A_{\mu}^{\varepsilon}}^{-1}f-\rbr{\op A_{\mu}^{0}}^{-1}f}_{p,\Omega}\lesssim\varepsilon^{s/p+s^{\dual}{\mkern-4mu }/p^{\dual}}\norm f_{p,\Omega}.\label{est: Convergence with e}
\end{equation}
The⋅constant depends only on the parameters~$d$, $s$, $s^{\dual}$,
$p$, $n$, $\mu$, the domain~$\Omega$, the $C^{0,1}$\nobreakdash-norms
of⋅$A$, $N$ and⋅$N^{\dual}$ and the constants in the  bounds~\textup{(\ref{est: Original operator | Norm of (A=001D4B-=0003BC)=00207B=0000B9})},
\textup{(\ref{est: Norm of (A=002070-=0003BC)=00207B=0000B9})}
and⋅\textup{(\dualref{est: Norm of (A=002070-=0003BC)=00207B=0000B9})}.
\end{thm}
The⋅other estimate in Theorem~\ref{thm: Convergence and Approximation with 1st corrector}
can be improved as well, but only if restricted to an interior of~$\Omega$.

\begin{thm}
\label{thm: Interior approximation with 1st corrector}Suppose that
\textup{(\ref{est: Original operator | Norm of (A=001D4B-=0003BC)=00207B=0000B9})},
\textup{(\ref{assu: N is Lipschitz})}, \textup{(\ref{est: Norm of (A=002070-=0003BC)=00207B=0000B9})}
and \textup{(\dualref{assu: N is Lipschitz})}, \textup{(\dualref{est: Norm of (A=002070-=0003BC)=00207B=0000B9})}
hold. Suppose further that for a given⋅$\chi\in C^{0,1}\rbr{\wb{\Omega}}$
with⋅$\supp\chi\subset\Omega$ there is $\chi^{\prime}\in C^{0,1}\rbr{\wb{\Omega}}$
with⋅$\supp\chi\subset\supp\chi^{\prime}\subset\Omega$ such that
for all⋅$\varepsilon\in\set E_{\mu}$  the interior energy estimate
\begin{equation}
\norm{D\chi u}_{p,\Omega}\lesssim\norm u_{p,\Omega}+\norm{\chi^{\prime}\op A_{\mu}^{\varepsilon}u}_{-1,p,\Omega}^{*},\qquad u\in\set W_{p}^{1}\rbr{\Omega;\C^{n}},\label{est: Local enery estimate}
\end{equation}
holds. Then for any⋅$\varepsilon\in\set E_{\mu}$ and~$f\in L_{p}\rbr{\Omega}^{n}$
\begin{equation}
\norm{D\chi\rbr{\rbr{\op A_{\mu}^{\varepsilon}}^{-1}f-\rbr{\op A_{\mu}^{0}}^{-1}f-\varepsilon\op K_{\mu}^{\varepsilon}f}}_{p,\Omega}\lesssim\varepsilon^{s/p+s^{\dual}{\mkern-4mu }/p^{\dual}}\norm f_{p,\Omega}.\label{est: Interior approximation with 1st corrector}
\end{equation}
The⋅constant depends only on the parameters~$d$, $s$, $s^{\dual}$,
$p$, $n$, $\mu$, the domain~$\Omega$,  the $C^{0,1}$\nobreakdash-norms
of⋅$A$, $N$, $N^{\dual}$ and⋅$\chi^{\prime}$ and the constants
in the  bounds~\textup{(\ref{est: Original operator | Norm of (A=001D4B-=0003BC)=00207B=0000B9})},
\textup{(\ref{est: Norm of (A=002070-=0003BC)=00207B=0000B9})},
\textup{(\dualref{est: Norm of (A=002070-=0003BC)=00207B=0000B9})}
and⋅\textup{(\ref{est: Local enery estimate})}.
\end{thm}
As⋅a corollary we have:
\begin{cor}
\label{cor: Interior approximation with 1st corrector}Let⋅$\Sigma$
be a subdomain of⋅$\Omega$ with~$\dist\rbr{\Sigma,\partial\Omega}>0$.
Suppose that the hypotheses of Theorem~\textup{\ref{thm: Interior approximation with 1st corrector}}
hold for $\chi$ such that⋅$\chi\restrict{\Sigma}=1$. Then, for any⋅$\varepsilon\in\set E_{\mu}$
and⋅$f\in L_{p}\rbr{\Omega}^{n}$,
\begin{equation}
\norm{D\rbr{\op A_{\mu}^{\varepsilon}}^{-1}f-D\rbr{\op A_{\mu}^{0}}^{-1}f-\varepsilon D\op K_{\mu}^{\varepsilon}f}_{p,\Sigma}\lesssim\varepsilon^{s/p+s^{\dual}{\mkern-4mu }/p^{\dual}}\norm f_{p,\Omega}.\label{est: Interior approximation with 1st corrector-2}
\end{equation}
The⋅constant depends only on the parameters~$d$, $s$, $s^{\dual}$,
$p$, $n$, $\mu$, the domain~$\Omega$, the $C^{0,1}$\nobreakdash-norms
of⋅$A$, $N$, $N^{\dual}$ and⋅$\chi^{\prime}$ and the constants
in the  bounds~\textup{(\ref{est: Original operator | Norm of (A=001D4B-=0003BC)=00207B=0000B9})},
\textup{(\ref{est: Norm of (A=002070-=0003BC)=00207B=0000B9})},
\textup{(\dualref{est: Norm of (A=002070-=0003BC)=00207B=0000B9})}
and⋅\textup{(\ref{est: Local enery estimate})}.
\end{cor}
The⋅next result follows from Corollary~\ref{cor: Interior approximation with 1st corrector}
in the same manner as Corollary~\ref{cor: Corollary 1} comes from
Theorem~\ref{thm: Convergence and Approximation with 1st corrector}.
\begin{cor}
\label{cor: Interior convergence for fractional derivative}Let⋅hypotheses
be as in Corollary~\textup{\ref{cor: Interior approximation with 1st corrector}}.
Then, for any⋅$r\in\rbr{0,1}$, $\varepsilon\in\set E_{\mu}$ and⋅$f\in L_{p}\rbr{\Omega}^{n}$,
\begin{equation}
\Seminorm{\rbr{\op A_{\mu}^{\varepsilon}}^{-1}f-\rbr{\op A_{\mu}^{0}}^{-1}f}_{r,p,\Sigma}\lesssim\varepsilon^{\rbr{s/p+s^{\dual}{\mkern-4mu }/p^{\dual}}\meet\rbr{1-r}}\norm f_{p,\Omega}.\label{est: Interior convergence for fractional derivative}
\end{equation}
The⋅constant depends only on the parameters~$d$, $r$, $s$, $s^{\dual}$,
$p$, $n$, $\mu$, the domain~$\Omega$, the $C^{0,1}$\nobreakdash-norms
of⋅$A$, $N$, $N^{\dual}$ and⋅$\chi^{\prime}$⋅and the constants
in the  bounds~\textup{(\ref{est: Original operator | Norm of (A=001D4B-=0003BC)=00207B=0000B9})},
\textup{(\ref{est: Norm of (A=002070-=0003BC)=00207B=0000B9})},
\textup{(\dualref{est: Norm of (A=002070-=0003BC)=00207B=0000B9})}
and⋅\textup{(\ref{est: Local enery estimate})}.
\end{cor}
\begin{rem}
The⋅corrector~$\varepsilon\op K_{\mu}^{\varepsilon}$ is usually
involved in an approximation for⋅$\rbr{\op A_{\mu}^{\varepsilon}}^{-1}$
in the ``energy'' norm. If⋅we want to approximate⋅$D\rbr{\op A_{\mu}^{\varepsilon}}^{-1}$
only, we may use the operator~$\tau^{\varepsilon}\op S^{\varepsilon}D_{2}\op K_{\mu}\rbr{\varepsilon}$
instead, because
\[
\varepsilon D\op K_{\mu}^{\varepsilon}=\varepsilon\tau^{\varepsilon}\op S^{\varepsilon}D_{1}\op K_{\mu}\rbr{\varepsilon}+\tau^{\varepsilon}\op S^{\varepsilon}D_{2}\op K_{\mu}\rbr{\varepsilon},
\]
where
\[
\norm{\tau^{\varepsilon}\op S^{\varepsilon}D_{1}\op K_{\mu}\rbr{\varepsilon}\hairspace f}_{p,\Omega}\lesssim\norm f_{p,\Omega}
\]
by⋅Lemma~\ref{lem: =0003C4=001D4BS=001D4B} and the estimate~(\ref{est: K(delta) is bounded}),
so this term can be absorbed into the error.
\end{rem}
\begin{rem}
The⋅results of Theorem~\ref{thm: Interior approximation with 1st corrector}
and Corollaries~\ref{cor: Interior approximation with 1st corrector}
and~\ref{cor: Interior convergence for fractional derivative} rely
on an a~priory bound~(\ref{est: Local enery estimate}). In⋅view
of Lemma~\ref{lem: Original operator | Interior energy estimate},
for a compactly supported function~$\chi$ this can be reduced to
a similar bound with a smaller exponent~$q\ge1$, provided that
 the resolvent of⋅$\op A^{\varepsilon}$ is bounded from⋅$\rbr{W_{q^{\dual}}^{1}\rbr{\Omega}^{n}}^{*}$
to~$W_{q}^{1}\rbr{\Omega}^{n}$.
\end{rem}
\begin{rem}
\label{rem: Boundary layer term}A⋅glance at⋅(\ref{est: Approximation with 1st corrector})
and⋅(\ref{est: Interior approximation with 1st corrector-2})  suggests
that the rate of approximation for⋅$D\rbr{\op A_{\mu}^{\varepsilon}}^{-1}$
becomes worse  only near the boundary of~$\Omega$. In⋅fact, one
can introduce a boundary-layer correction term~$\op B_{\mu}^{\varepsilon}$
so that for any⋅$\varepsilon\in\set E_{\mu}$ and⋅$f\in L_{p}\rbr{\Omega}^{n}$
\begin{equation}
\norm{\rbr{\op A_{\mu}^{\varepsilon}}^{-1}f-\rbr{\op A_{\mu}^{0}}^{-1}f-\varepsilon\op K_{\mu}^{\varepsilon}f-\op B_{\mu}^{\varepsilon}f}_{1,p,\Omega}\lesssim\varepsilon^{s}\norm f_{p,\Omega},\label{est: Boundary-layer correction term}
\end{equation}
see⋅(\ref{def: Boundary part}) and Lemmas~\ref{lem: Ie} and~\ref{lem: De}.
For⋅$s=1$ and⋅$p=2$, such a result was the starting point of the
approach suggested in~\cite{ZhPas:2005} (see⋅also~\cite{PasT:2007},
\cite{PSu:2012}, \cite{Su:2013-1} and~\cite{Su:2013-2}). However,
the construction of⋅$\op B_{\mu}^{\varepsilon}$ is no simpler than
the original problem and actually amounts to finding the inverse of~$\op A_{\mu}^{\varepsilon}$.
Thus, that approach required  further analysis of the boundary-layer
correction term to obtain bounds on its norms.  We also note that
if $\Omega=\R^{d}$ (or, more generally, $\Omega$ is a flat manifold
without boundary, such as, e.g.,~$\T^{d}$), then⋅$\op B_{\mu}^{\varepsilon}=0$.
This enables one to improve the rates in⋅(\ref{est: Convergence})–(\ref{est: Approximation with 1st corrector})
to⋅$\varepsilon^{s}$, which, at least for⋅$s=1$, is known to be
sharp.
\end{rem}

\begin{rem}
\label{rem: Main results | Regularity of the effective operator}By⋅inspection
of the proofs, one can see that the estimates in  Theorem~\ref{thm: Convergence and Approximation with 1st corrector}–Corollary~\ref{cor: Interior convergence for fractional derivative}
follow from  inequalities with⋅$\abs{\rbr{\op A_{\mu}^{0}}^{-1}f}_{1+s,p,\Omega}$
in place of⋅$\norm f_{p,\Omega}$ on the right. Thus, if⋅(\ref{est: Norm of (A=002070-=0003BC)=00207B=0000B9})
fails to hold, but $\rbr{\op A_{\mu}^{0}}^{-1}f\in\set W_{p}^{1}\rbr{\Omega;\C^{n}}\cap\Lambda_{p}^{1+s}\rbr{\Omega}^{n}$
for some⋅$f\in L_{p}\rbr{\Omega}^{n}$, then for fixed such~$f$
we still have, e.g., results similar to Theorem~\ref{thm: Convergence and Approximation with 1st corrector}
and⋅Corollary~\ref{cor: Corollary 1}.
\end{rem}

\section{\label{sec: Examples}Examples}

In⋅the examples below we assume that $\Omega$ is a bounded $C^{1,1}$\nobreakdash-do\-main in~$\R^{d}$. For⋅each⋅$1\le k\le n$,
let $\Gamma_{k}$ either be a $C^{1,1}$\nobreakdash-submanifold
of⋅$\partial\Omega$ with boundary (possibly disconnected) or coincide
with⋅$\emptyset$ or~$\partial\Omega$. By⋅$W_{q}^{1}\rbr{\Omega;\Gamma_{k}}$,
we denote the $W_{q}^{1}$\nobreakdash-closure of smooth functions
on⋅$\Omega$ that vanish near~$\Gamma_{k}$. It⋅is⋅known  that
if $1<q_{0}<q_{1}<\infty$, then $W_{q_{\theta}}^{1}\rbr{\Omega;\Gamma_{k}}$
is the complex interpolation space between⋅$W_{q_{0}}^{1}\rbr{\Omega;\Gamma_{k}}$
and⋅$W_{q_{1}}^{1}\rbr{\Omega;\Gamma_{k}}$, where⋅$1/q_{\theta}=\rbr{1-\theta}/q_{0}+\theta/q_{1}$
and⋅$\theta\in\rbr{0,1}$. We set⋅$\set W_{q}^{1}\rbr{\Omega;\C^{n}}=\oplus_{k=1}^{n}W_{q}^{1}\rbr{\Omega;\Gamma_{k}}$.

\subsection{Strongly elliptic operators}

Let⋅$p=2$, and let $\set H^{1}\rbr{\Omega;\C^{n}}=\set W_{2}^{1}\rbr{\Omega;\C^{n}}$.
Suppose that the operator~$\op A^{\varepsilon}$ is weakly coercive
 uniformly in⋅$\varepsilon$ for⋅$\varepsilon$ sufficiently small,
that is, there are $\varepsilon_{0}\in(0,1]$ and $c_{A}>0$ and⋅$C_{A}<\infty$
so that for all⋅$\varepsilon\in\set E=(0,\varepsilon_{0}]$
\begin{equation}
\Re\rbr{\op A^{\varepsilon}u,u}_{\Omega}+C_{A}\norm u_{2,\Omega}^{2}\ge c_{A}\norm{Du}_{2,\Omega}^{2},\qquad u\in\set H^{1}\rbr{\Omega;\C^{n}}.\label{est: A^e is elliptic}
\end{equation}
With⋅this assumption, $\op A^{\varepsilon}$ becomes strongly elliptic,
which means that the function~$A$ satisfies the Legendre–Hadamard
condition
\begin{equation}
\Re\abr{A\rbr{\wc,\wc}\hairspace\xi\otimes\eta,\xi\otimes\eta}\ge c_{A}\abs{\xi}^{2}\abs{\eta}^{2},\qquad\xi\in\R^{d},\eta\in\C^{n}\label{est: Legendre-Hadamard condition}
\end{equation}
(see⋅Lemma~\ref{lem: Strong coercivity of A(x)} below). What is
more,  a simple calculation based on boundedness and coercivity of⋅$\op A^{\varepsilon}$
shows that if $\varepsilon\in\set E$, then $\op A^{\varepsilon}$
is an $m$-sectorial operator with sector
\[
\set S=\bigcbr{z\in\C\colon\abs{\Im z}\le c_{A}^{-1}\norm A_{L_{\infty}}\rbr{\Re z+C_{A}}}
\]
independent of⋅$\varepsilon$, and therefore (\ref{est: Original operator | Norm of (A=001D4B-=0003BC)=00207B=0000B9})
holds for⋅$p=2$ and any⋅$\varepsilon\in\set E$ as long as~$\mu\notin\set S$.
It⋅then follows from Shneiberg's stability theorem~\cite{Shneiberg:1974}
(see also~\cite[Section~17.2]{Agranovich:2013}) that the estimate~(\ref{est: Original operator | Norm of (A=001D4B-=0003BC)=00207B=0000B9})
is valid for any⋅$p$ satisfying
\[
\bigabs{\tan\tfrac{\pi}{2}\bigrbr{\tfrac{1}{2}-\tfrac{1}{p}}}<\inf_{\varepsilon\in\set E}\Bigrbr{\,\inf_{q\in\rbr{1,\infty}}\norm{\op A_{\mu}^{\varepsilon}}_{\set W_{q}^{1}\to\set W_{q}^{-1}}^{-1}\norm{\rbr{\op A_{\mu}^{\varepsilon}}^{-1}}_{\set H^{-1}\to\set H^{1}}^{-1}\!}.
\]
Using⋅(\ref{est: Original operator | A=001D4B is bounded}) and⋅(\ref{est: A^e is elliptic}),
one can easily see that the right-hand side here is  bounded below
by a positive constant~$C_{\mu}$, depending only on⋅$\mu$ and the
ellipticity constants~$c_{A}$, $C_{A}$ and~$\norm A_{C}$. Thus,
we have⋅(\ref{est: Original operator | Norm of (A=001D4B-=0003BC)=00207B=0000B9})
for all⋅$p\in\set P_{\mu}=\rbr{p_{\mu}^{\dual},p_{\mu}}$, where
$p_{\mu}$ solves
\[
\tan\tfrac{\pi}{2}\bigrbr{\tfrac{1}{2}-\tfrac{1}{p_{\mu}}}=C_{\mu}.
\]

Let $\set P=\cup_{\mu\notin\set S}\set P_{\mu}$. We show that
(\ref{est: Original operator | Norm of (A=001D4B-=0003BC)=00207B=0000B9})
holds, in fact, for any⋅$p\in\set P$ and⋅$\mu\notin\set S$ uniformly
in~$\varepsilon\in\set E$. Indeed, suppose that $\mu,\nu\notin\set S$
and choose⋅$p\in\rbr{2,p_{\nu}}$, so that $\rbr{\op A_{\nu}^{\varepsilon}}^{-1}$
maps⋅$\set W_{p}^{-1}\rbr{\Omega;\C^{n}}$ to⋅$\set W_{p}^{1}\rbr{\Omega;\C^{n}}$
continuously. From⋅the Sobolev embedding theorem, we know that $L_{2}\rbr{\Omega}$
is  embedded in⋅$W_{q^{\dual}}^{1}\rbr{\Omega}^{*}$ for⋅$q\in\sbr{2,2^{*}}$,
and in particular in⋅$\set W_{2^{*}\meet p}^{-1}\rbr{\Omega;\C^{n}}$
(see⋅(\ref{est: Norm of pi})). Hence, the first resolvent identity
\begin{equation}
\rbr{\op A_{\mu}^{\varepsilon}}^{-1}=\rbr{\op A_{\nu}^{\varepsilon}}^{-1}+\rbr{\mu-\nu}\rbr{\op A_{\nu}^{\varepsilon}}^{-1}\rbr{\op A_{\mu}^{\varepsilon}}^{-1}\label{eq: 1st resolvent id}
\end{equation}
yields that $\rbr{\op A_{\mu}^{\varepsilon}}^{-1}$ is bounded from⋅$\set W_{2^{*}\meet p}^{-1}\rbr{\Omega;\C^{n}}$
to⋅$\set W_{2^{*}\meet p}^{1}\rbr{\Omega;\C^{n}}$. Repeating this
procedure finitely many times, if need be, we conclude that the operator~$\rbr{\op A_{\mu}^{\varepsilon}}^{-1}$
is bounded from⋅$\set W_{p}^{-1}\rbr{\Omega;\C^{n}}$ to⋅$\set W_{p}^{1}\rbr{\Omega;\C^{n}}$
as well.
\begin{rem}
No⋅necessary and sufficient algebraic condition for⋅$A$ to assure⋅(\ref{est: A^e is elliptic})
is known. A⋅simpler condition not involving⋅$\varepsilon$ and still
implying the weak coercivity  on⋅$\mathring{H}^{1}\rbr{\Omega}^{n}$
is that for some⋅$c>0$ and all~$x\in\Omega$
\begin{equation}
\Re\rbr{A\rbr{x,\wc}\hairspace Du,Du}_{\R^{d}}\ge c\norm{Du}_{2,\R^{d}}^{2},\qquad u\in H^{1}\rbr{\R^{d}}^{n}.\label{est: Sufficient condition for A to be coercive}
\end{equation}
That this hypothesis suffices can be seen by noticing that (\ref{est: Sufficient condition for A to be coercive})
is invariant under dilation and therefore remains true with⋅$A\rbr{x,\varepsilon^{-1}y}$
in place of~$A\rbr{x,y}$. Since $A$ is uniformly continuous in
the first variable, a localization argument then leads to⋅(\ref{est: A^e is elliptic}),
with $c_{A}<c$, $C_{A}>0$ and⋅$\set H^{1}\rbr{\Omega;\C^{n}}=\mathring{H}^{1}\rbr{\Omega}^{n}$.

To⋅give an example of⋅$A$ satisfying the strong coercivity condition
on⋅$\mathring{H}^{1}\rbr{\Omega}^{n}$ (i.e., with⋅$C_{A}=0$), take
a matrix first-order differential operator~$b\rbr D$ with symbol
\[
\xi\mapsto b\rbr{\xi}=\sum_{k=1}^{d}b_{k}\xi_{k},
\]
where~$b_{k}\in\C^{m\times n}$. Suppose that the symbol has the
property that $\rank b\rbr{\xi}=n$ for any⋅$\xi\in\R^{d}\setminus\cbr 0$,
or, equivalently, that there is⋅$c_{b}>0$ such that
\[
b\rbr{\xi}^{*}b\rbr{\xi}\ge c_{b}\abs{\xi}^{2},\qquad\xi\in\R^{d}.
\]
Extending $u\in\mathring{H}^{1}\rbr{\Omega}^{n}$ by zero outside⋅$\Omega$
and applying the Fourier transform, we see that the operator~$b\rbr D^{*}b\rbr D$
is strongly coercive on⋅$\mathring{H}^{1}\rbr{\Omega}^{n}$:
\begin{equation}
\norm{b\rbr D\hairspace u}_{2,\Omega}^{2}\ge c_{b}\norm{Du}_{2,\Omega}^{2},\qquad u\in\mathring{H}^{1}\rbr{\Omega}^{n}.\label{est: Coercivity of b(D)*b(D) on Ho}
\end{equation}
Let $g\in C^{0,1}\rbr{\wb{\Omega};\wt L_{\infty}\rbr{\cell}}^{m\times m}$
with⋅$\Re g$ uniformly positive definite and let⋅$A_{kl}=b_{k}^{*}gb_{l}$.
Then, by⋅(\ref{est: Coercivity of b(D)*b(D) on Ho}),
\[
\begin{aligned}\Re\rbr{A^{\varepsilon}Du,Du}_{\Omega} & =\Re\rbr{g^{\varepsilon}b\rbr D\hairspace u,b\rbr D\hairspace u}_{\Omega}\\
 & \ge c_{b}\norm{\rbr{\Re g}^{-1}}_{L_{\infty}}^{-1}\norm{Du}_{2,\Omega}^{2}
\end{aligned}
\]
for⋅all⋅$u\in\mathring{H}^{1}\rbr{\Omega}^{n}$. Purely periodic
operators of this type were studied, e.g., in⋅\cite{PSu:2012} and⋅\cite{Su:2013-1}.

For⋅coercivity on⋅$H^{1}\rbr{\Omega}^{n}$, one needs to require a
stronger condition on the symbol, namely, that $\rank b\rbr{\xi}=n$
for any⋅$\xi\in\C^{d}\setminus\cbr 0$, not just⋅$\xi\in\R^{d}\setminus\cbr 0$,
implying that
\begin{equation}
\norm{b\rbr D\hairspace u}_{2,\Omega}^{2}\ge c_{b}\norm{Du}_{2,\Omega}^{2}-C_{b}\norm u_{2,\Omega}^{2},\qquad u\in H^{1}\rbr{\Omega}^{n},\label{est: Coercivity of b(D)*b(D) on H}
\end{equation}
with⋅$c_{b}>0$, see~\cite[Section~3.7, Theorem~7.8]{Necas:2012}.
Then, obviously, for any⋅$u\in H^{1}\rbr{\Omega}^{n}$
\[
\begin{aligned}\Re\rbr{A^{\varepsilon}Du,Du}_{\Omega} & =\Re\rbr{g^{\varepsilon}b\rbr D\hairspace u,b\rbr D\hairspace u}_{\Omega}\\
 & \ge\norm{\rbr{\Re g}^{-1}}_{L_{\infty}}^{-1}\bigrbr{c_{b}\norm{Du}_{2,\Omega}^{2}-C_{b}\norm u_{2,\Omega}^{2}},
\end{aligned}
\]
where⋅$A$ and⋅$g$ are as above. Such operators in the purely periodic
setting appeared in~\cite{Su:2013-2}.

\end{rem}
Now we turn to the cell problem and the effective operator. The⋅first
thing that we need to check is that the cell problem~(\ref{def: N})
has a unique solution for which⋅(\ref{assu: N is Lipschitz}) holds.
Lemma~\ref{lm: A(x)} contains a sufficient condition to conclude
these, and we will see in a moment that the operator~$\op A\rbr x$
does indeed meet the hypothesis of that lemma.
\begin{lem}
\label{lem: Strong coercivity of A(x)}Assume that~\textup{(\ref{est: A^e is elliptic})}
holds. Then for any⋅$x\in\Omega$
\begin{equation}
\Re\rbr{\op A\rbr x\hairspace u,u}_{\cell}\ge c_{A}\norm{Du}_{2,\cell}^{2},\qquad u\in\wt H^{1}\rbr{\cell}^{n}.\label{est: Coercivity of D*A(x,=0022C5)D on Q}
\end{equation}

\end{lem}
\begin{proof}
Fix $\smash[t]{u^{\rbr{\varepsilon}}}=\varepsilon u^{\varepsilon}\varphi$
with $u\in\wt C^{1}\rbr{\cell}^{n}$ and~$\varphi\in C_{c}^{\infty}\rbr{\Omega}$.
We substitute $\smash[t]{u^{\rbr{\varepsilon}}}$ into (\ref{est: A^e is elliptic})
and let $\varepsilon$ tend to~0. Then, because $\smash[t]{u^{\rbr{\varepsilon}}}$
and⋅$\smash[t]{Du^{\rbr{\varepsilon}}}-\rbr{Du}^{\varepsilon}\varphi$
converge in⋅$L_{2}$ to~$0$,
\[
\lim_{\varepsilon\to0}\Re\int_{\Omega}\abr{A^{\varepsilon}\rbr x\hairspace\rbr{Du}^{\varepsilon}\rbr x,\rbr{Du}^{\varepsilon}\rbr x}\abs{\varphi\rbr x}^{2}\dd x\ge\lim_{\varepsilon\to0}c_{A}\int_{\Omega}\abs{\rbr{Du}^{\varepsilon}\rbr x}^{2}\abs{\varphi\rbr x}^{2}\dd x.
\]
It⋅is well known  that if⋅$f\in C_{c}\rbr{\R^{d};\wt L_{\infty}\rbr{\cell}}$,
then
\[
\lim_{\varepsilon\to0}\int_{\R^{d}}f^{\varepsilon}\rbr x\dd x=\int_{\R^{d}}\int_{\cell}f\rbr{x,y}\dd x\dd y
\]
(see, e.g.,~\cite[Lemmas~5.5 and~5.6]{Al:1992}). As⋅a⋅result,
\[
\Re\int_{\Omega}\int_{\cell}\abr{A\rbr{x,y}\hairspace Du\rbr y,Du\rbr y}\abs{\varphi\rbr x}^{2}\dd x\dd y\ge c_{A}\int_{\Omega}\int_{\cell}\abs{Du\rbr y}^{2}\abs{\varphi\rbr x}^{2}\dd x\dd y.
\]
But $\varphi$ is an arbitrary function in $C_{c}^{\infty}\rbr{\Omega}$
and $A$ is uniformly continuous in the first variable,⋅so
\[
\Re\int_{\cell}\abr{A\rbr{x,y}\hairspace Du\rbr y,Du\rbr y}\dd y\ge c_{A}\int_{\cell}\abs{Du\rbr y}^{2}\dd y
\]
for⋅all⋅$x\in\Omega$, as required.
\end{proof}
We have shown that, for any⋅$x\in\Omega$, the operator~$\op A\rbr x$
is an isomorphism of⋅$\wt H^{1}\rbr{\cell}^{n}\!/\C$ onto⋅$\wt H^{-1}\rbr{\cell}^{n}$
and that the ellipticity constants of⋅$\op A\rbr x$ are better that
those of⋅$\op A^{\varepsilon}$ (cf.⋅(\ref{est: A^e is elliptic})
with~(\ref{est: Coercivity of D*A(x,=0022C5)D on Q})). Then the
Shneiberg stability theorem yields that $\op A\rbr x$ is  an isomorphism
of⋅$\wt W_{p}^{1}\rbr{\cell}^{n}\!/\C$ onto⋅$\wt W_{p}^{-1}\rbr{\cell}^{n}$
for any~$p\in\set P$. Thus, the hypothesis of Lemma~\ref{lm: A(x)}
 is verified.

As⋅for the effective operator, one can prove that, for any⋅$\mu\notin\set S$,
the inverse for⋅$\op A_{\mu}^{\varepsilon}$ converges in the weak
operator topology and then the limit is an isomorphism of⋅$\set H^{-1}\rbr{\Omega;\C^{n}}$
onto⋅$\set H^{1}\rbr{\Omega;\C^{n}}$, which is, in fact, the inverse
for⋅$\op A_{\mu}^{0}$, see~\cite[Lemma~6.2]{Tartar:2010}. We conclude
that $\rbr{\op A_{\mu}^{0}}^{-1}$ is also an isomorphism as a mapping
from⋅$\set W_{p}^{-1}\rbr{\Omega;\C^{n}}$ to⋅$\set W_{p}^{1}\rbr{\Omega;\C^{n}}$
for any⋅$p\in\set P$, because it is the limit in the weak operator
topology. In⋅the case of the Dirichlet and Neumann problems, this,
and the fact that $A^{0}$ is Lipschitz, yields the assumption~(\ref{est: Norm of (A=002070-=0003BC)=00207B=0000B9})
with⋅$s=1$ and any⋅$p\in\set P$ (see, e.g.,~\cite[Chapter~4]{McL:2000}).
In⋅the case of the mixed Dirichlet–Neumann problem, we assume additionally
that $A^{0}$ is self-adjoint; then it follows from~\cite{Savare:1997}
that (\ref{est: Norm of (A=002070-=0003BC)=00207B=0000B9}) holds
for⋅$s=1/2$ and⋅$p=2$.

Of⋅course, all these results are true for the dual counterparts with
the same range of⋅$p$, because $c_{A^{\dual}}=c_{A}$, $C_{A^{\dual}}=C_{A}$
and~$\norm{A^{\dual}}_{L_{\infty}}=\norm A_{L_{\infty}}$.

It⋅remains to discuss the interior energy estimate~(\ref{est: Local enery estimate}).
Let⋅$p\in\set P\cap[2,\infty)$. Applying the functional~$\op A_{\mu}^{\varepsilon}u$
to⋅$\abs{\chi}^{2}u$, where⋅$\chi\in C_{c}^{0,1}\rbr{\Omega}$, and
using⋅(\ref{est: A^e is elliptic}), we arrive at the well-known Caccioppoli
inequality:
\[
\norm{\chi Du}_{2,\supp\chi}\lesssim\norm u_{2,\supp\chi}+\norm{\chi\op A_{\mu}^{\varepsilon}u}_{-1,2,\Omega}^{*},\qquad u\in\set H^{1}\rbr{\Omega;\C^{n}}.
\]
Therefore, Lemma~\ref{lem: Original operator | Interior energy estimate}
with⋅$q=2$ yields~(\ref{est: Local enery estimate}).

To⋅summarize, assume that $\op A^{\varepsilon}$ satisfies the coercivity
condition~(\ref{est: A^e is elliptic}). Then for Dirichlet or Neumann
boundary conditions the global results (see Theorem~\ref{thm: Convergence and Approximation with 1st corrector}–Theorem~\ref{thm: Convergence with e}) are valid with⋅$s=s^{\dual}=1$ and⋅$p\in\set P$ and the local
results (see Theorem~\ref{thm: Interior approximation with 1st corrector}–Corollary~\ref{cor: Interior convergence for fractional derivative})
are valid with⋅$s=s^{\dual}=1$ and~$p\in\set P\cap[2,\infty)$;
for mixed Dirichlet–Neumann boundary conditions all these results
hold true with⋅$s=s^{\dual}=1/2$ and⋅$p=2$, provided that $A^{0}$
is self-adjoint.
\begin{rem}
The⋅constants~$p_{\mu}$, and hence the interval~$\set P$, can
be expressed explicitly. We note that generally one would not expect
$\set P$ to be too wide. In⋅fact,  it must shrink to⋅$\cbr 2$
as the ellipticity of the family~$\op A^{\varepsilon}$ becomes ``bad''
(that is, the ratio~$c_{A}^{-1}\norm A_{L_{\infty}}$ grows), see~\cite{Mey:1963}.
In⋅the⋅next subsection we provide an example where $p$ may be chosen
 arbitrary large.
\end{rem}

\subsection{Strongly elliptic operators with $\boldsymbol{\protect\VMO}$\protect\nobreakdash-coefficients}

Throughout this subsection, we restrict our attention to the Dirichlet
and Neumann problems. Let $\op A^{\varepsilon}$ be as in the previous
subsection. Assume further that $A\in L_{\infty}\rbr{\Omega;\VMO\rbr{\R^{d}}}$,
meaning that $\sup_{x\in\Omega}\eta_{A\rbr{x,\wc}}\rbr r\to0$ as⋅$r\to0$.

Using the reflection technique, we extend $A$ to be a function belonging
to both $C_{c}^{0,1}\rbr{\R^{d};\wt L_{\infty}\rbr{\cell}}$ and~$L_{\infty}\rbr{\R^{d};\VMO\rbr{\R^{d}}}$.
Notice that $A^{\varepsilon}$ is then a $\VMO$\nobreakdash-function.
Indeed, $A^{\varepsilon}$ obviously belongs to the space~$\BMO\rbr{\R^{d}}$,
with $\norm{A^{\varepsilon}}_{\BMO}\le2\norm A_{L_{\infty}}$. Next,
after dilation, we may suppose that~$\varepsilon=1$. Given an⋅$\epsilon>0$
small, there is $r>0$ such that $\osc_{A}\rbr r<\epsilon/3$ and
$\eta_{A\rbr{x,\wc}}\rbr r<\epsilon/3$. Then, since
\[
\begin{aligned}\dashint_{B_{R}\rbr{x_{0}}}\abs{A^{1}\rbr x-m_{B_{R}\rbr{x_{0}}}\rbr{A^{1}}}\dd x & \le\dashint_{B_{R}\rbr{x_{0}}}\abs{A\rbr{x_{0},x}-m_{B_{R}\rbr{x_{0}}}\rbr{A\rbr{x_{0},\wc}}}\dd x\\
 & \quad+2\dashint_{B_{R}\rbr{x_{0}}}\abs{A\rbr{x,x}-A\rbr{x_{0},x}}\dd x,
\end{aligned}
\]
we have
\[
\eta_{A^{1}}\rbr r\le\sup_{x_{0}\in\R^{d}}\eta_{A\rbr{x_{0},\wc}}\rbr r+2\osc_{A}\rbr r<\epsilon,
\]
and the claim follows.

As⋅a⋅result, if $\varepsilon\in\set E$ is fixed and⋅$\mu\notin\set S$,
the inverse of⋅$\op A_{\mu}^{\varepsilon}$ is a continuous map from⋅$\set W_{p}^{-1}\rbr{\Omega;\C^{n}}$
to⋅$\set W_{p}^{1}\rbr{\Omega;\C^{n}}$ for each⋅$p\in\rbr{1,\infty}$,
see~\cite{Shen:2018}. Hence, in order to prove⋅(\ref{est: Original operator | Norm of (A=001D4B-=0003BC)=00207B=0000B9}),
we need only show that its norm is uniformly bounded in~$\varepsilon$.
We do this by  treating⋅$\op A^{\varepsilon}$ as a local perturbation
of a purely periodic operator and then applying results for purely
periodic operators with rapidly oscillating coefficients.

First observe that if $B_{R}$ is a ball with center in⋅$\wb{\Omega}$
and radius~$R$, then, by⋅(\ref{est: A^e is elliptic}),
\[
\Re\rbr{A^{\varepsilon}\rbr{x_{0},\wc}\hairspace Dv,Dv}_{B_{R}\cap\Omega}\ge\rbr{c_{A}-\osc_{A}\rbr R}\hairspace\norm{Dv}_{2,B_{R}\cap\Omega}^{2}-C_{A}\norm v_{2,B_{R}\cap\Omega}^{2}
\]
for⋅all⋅$v$ in⋅$\set H^{1}\rbr{B_{R}\cap\Omega;\C^{n}}$, the space
of functions whose zero extensions to⋅$\Omega$ belong to⋅$\set H^{1}\rbr{\Omega;\C^{n}}$.
 It⋅follows that for⋅$R$ small enough, the operator~$D^{*}A^{\varepsilon}\rbr{x_{0},\wc}\hairspace D$
from⋅$\set H^{1}\rbr{B_{R}\cap\Omega;\C^{n}}$ to the dual space~$\set H^{-1}\rbr{B_{R}\cap\Omega;\C^{n}}$
is $m$\nobreakdash-sectorial, with sector
\[
\set S_{R}=\bigcbr{z\in\C\colon\abs{\Im z}\le\rbr{c_{A}-\osc_{A}\rbr R}^{-1}\norm A_{L_{\infty}}\rbr{\Re z+C_{A}}}
\]
converging pointwise to $\set S$ as⋅$R\to0$, that is, $\dist\rbr{z,\set S_{R}}\to\dist\rbr{z,\set S}$
for~$z\in\C$.

Now,  fix $\mu\not\in\set S$ and find $R_{0}>0$ such that $\osc_{A}\rbr R\le c_{A}/2$
and~$\mu\notin\set S_{R}$ as long as~$R\le R_{0}$. Let⋅$F\in C_{c}^{\infty}\rbr{\Omega}^{dn}$
and $u_{\varepsilon}=\rbr{\op A_{\mu}^{\varepsilon}}^{-1}D^{*}F$.
Take⋅$\chi\in C_{c}^{\infty}\rbr{B_{R}}$, $R\le R_{0}$, with the
properties that $0\le\chi\rbr x\le1$ and⋅$\chi=1$ on~$1/2B_{R}$.
Then $v_{\varepsilon}=\chi u_{\varepsilon}$  obviously satisfies
\[
D^{*}A^{\varepsilon}\rbr{x_{0},\wc}\hairspace Dv_{\varepsilon}-\mu v_{\varepsilon}=D^{*}\rbr{A^{\varepsilon}\rbr{x_{0},\wc}-A^{\varepsilon}}Dv_{\varepsilon}+g
\]
in⋅the sense of functionals on⋅$\set H^{1}\rbr{B_{R_{0}}\cap\Omega;\C^{n}}$,
where ${g=\chi D^{*}F+D^{*}\rbr{A^{\varepsilon}D\chi\cdot u_{\varepsilon}}}-\rbr{D\chi}^{*}\cdot A^{\varepsilon}Du_{\varepsilon}$.
This is a purely periodic problem, for which we know that the operator~$D^{*}A^{\varepsilon}\rbr{x_{0},\wc}\hairspace D-\mu$
is an isomorphism of⋅$\set W_{q}^{1}\rbr{B_{R_{0}}\cap\Omega;\C^{n}}$
onto⋅$\set W_{q}^{-1}\rbr{B_{R_{0}}\cap\Omega;\C^{n}}$ for any⋅$q\in\rbr{1,\infty}$,
with uniformly  bounded inverse, see~\cite{Shen:2018}. Assuming
 that⋅$p>2$ (the⋅other case will follow by duality), we immediately
find that
\[
\norm{Dv_{\varepsilon}}_{2^{*}\meet p,B_{R_{0}}\cap\Omega}\lesssim\osc_{A}\rbr R\hairspace\norm{Dv_{\varepsilon}}_{2^{*}\meet p,B_{R_{0}}\cap\Omega}+\norm g_{-1,2^{*}\meet p,B_{R_{0}}\cap\Omega}^{*},
\]
the⋅constant not depending on~$R$. Choosing $R$ sufficiently
small, we may absorb the first term on the right into the left-hand
side. Since
\[
\norm g_{-1,2^{*}\meet p,B_{R_{0}}\cap\Omega}^{*}\lesssim\norm F_{2^{*}\meet p,B_{R}\cap\Omega}+\norm{u_{\varepsilon}}_{1,2,B_{R}\cap\Omega}
\]
(we have used the Sobolev embedding theorem to estimate the⋅$L_{2^{*}\meet p}$\nobreakdash-norm
of⋅$u_{\varepsilon}$ and the⋅$\set W_{2^{*}\meet p}^{-1}$\nobreakdash-norm
of~$Du_{\varepsilon}$), it follows that
\[
\norm{Du_{\varepsilon}}_{2^{*}\meet p,1/2B_{R}\cap\Omega}\lesssim\norm F_{2^{*}\meet p,B_{R}\cap\Omega}+\norm{u_{\varepsilon}}_{1,2,B_{R}\cap\Omega}.
\]
Now, cover⋅$\Omega$ with balls of radius~$R$ to obtain
\[
\norm{Du_{\varepsilon}}_{2^{*}\meet p,\Omega}\lesssim\norm F_{2^{*}\meet p,\Omega}+\norm{u_{\varepsilon}}_{1,2,\Omega}\lesssim\norm F_{2^{*}\meet p,\Omega}.
\]
After⋅a⋅finite number of repetitions, if need be, we get
\[
\norm{Du_{\varepsilon}}_{p,\Omega}\lesssim\norm F_{p,\Omega}.
\]

Next, the hypothesis of Lemma~\ref{lm: A(x)} is satisfied, because
of Lemma~\ref{lem: Strong coercivity of A(x)} and the fact that
$A\rbr{x,\wc}\in\VMO\rbr{\R^{d}}$ with $\VMO$\nobreakdash-modulus
bounded uniformly in~$x$. Finally,  (\ref{est: Norm of (A=002070-=0003BC)=00207B=0000B9})
 and⋅(\ref{est: Local enery estimate}) hold for, respectively, $s=1$
and any⋅$p\in\rbr{1,\infty}$ and any⋅$p\in[2,\infty)$, as indicated
previously.

Summarizing, if $A\in L_{\infty}\rbr{\Omega;\VMO\rbr{\R^{d}}}$ satisfies
the coercivity condition~(\ref{est: A^e is elliptic}), then the
global results (see Theorem~\ref{thm: Convergence and Approximation with 1st corrector}–Theorem~\ref{thm: Convergence with e})
are valid with⋅$s=s^{\dual}=1$ and⋅$p\in\rbr{1,\infty}$ and the
local results (see Theorem~\ref{thm: Interior approximation with 1st corrector}–Corollary~\ref{cor: Interior convergence for fractional derivative})
are valid with⋅$s=s^{\dual}=1$ and~$p\in[2,\infty)$.

\section{\label{sec: Proof of the main results}Proof of the  main results}

We start with a resolvent identity involving⋅$\rbr{\op A_{\mu}^{\varepsilon}}^{-1}$,
$\rbr{\op A_{\mu}^{0}}^{-1}$ and⋅$\op K_{\mu}^{\varepsilon}$ that
will play a central role in the proof.

\subsection{The resolvent identity}

Fix⋅$f\in L_{p}\rbr{\Omega}^{n}$ and⋅$f^{\dual}\in\rbr{W_{p}^{1}\rbr{\Omega}^{n}}^{*}$.
For⋅$\delta=\varepsilon$, we set $u_{0}=\op R_{\mu}^{0}f$, $u_{0,\delta}=\op R_{\mu}^{0}\rbr{\delta}\hairspace f$,
$U=\op K_{\mu}f$, $U_{\delta}=\op K_{\mu}\rbr{\delta}\hairspace f$
and⋅$U_{\varepsilon,\delta}=\tau^{\varepsilon}\op S^{\varepsilon}U_{\delta}=\op K_{\mu}^{\varepsilon}f$.
Then we have
\[
\rbr{\rbr{\op A_{\mu}^{\varepsilon}}^{-1}f-\rbr{\op A_{\mu}^{0}}^{-1}f-\varepsilon\op K_{\mu}^{\varepsilon}f,f^{\dual}}_{\Omega}=\rbr{f,u_{\varepsilon}^{\dual}}_{\Omega}-\rbr{u_{0},f^{\dual}}_{\Omega}-\varepsilon\rbr{U_{\varepsilon,\delta},f^{\dual}}_{\Omega},
\]
where, as usual, the ``$\dual$'' labels the dual counterparts,
e.g.,⋅$u_{\varepsilon}^{\dual}=\rbr{\rbr{\op A_{\mu}^{\varepsilon}}^{\dual}}^{-1}f^{\dual}$.
By⋅definition of⋅$u_{0}$ and~$u_{\varepsilon}^{\dual}$,
\[
\rbr{f,u_{\varepsilon}^{\dual}}_{\Omega}-\rbr{u_{0},f^{\dual}}_{\Omega}=\rbr{A^{0}Du_{0},Du_{\varepsilon}^{\dual}}_{\Omega}-\rbr{A^{\varepsilon}Du_{0},Du_{\varepsilon}^{\dual}}_{\Omega}.
\]
Choose a function~$\rho_{\varepsilon}\in C^{0,1}\rbr{\wb{\Omega}}$
with support in the closure of⋅$\rbr{\partial\Omega}_{3\varepsilon}\cap\Omega$
and values in⋅$\sbr{0,1}$ such that $\rho_{\varepsilon}\restrict{\rbr{\partial\Omega}_{2\varepsilon}\cap\Omega}=1$
and $\norm{D\rho_{\varepsilon}}_{\infty,\Omega}\lesssim\varepsilon^{-1}$.
For⋅example, we may take $\rho_{\varepsilon}\rbr x=3-\rbr{\cellradius\varepsilon}^{-1}\dist\rbr{x,\partial\Omega}$
for⋅$x\in\Omega\cap\rbr{\partial\Omega}_{3\varepsilon}\!\setminus\rbr{\partial\Omega}_{2\varepsilon}$.
If⋅$\chi_{\varepsilon}=1-\rho_{\varepsilon}$, then $\chi_{\varepsilon}U_{\varepsilon,\delta}\in\set W_{p}^{1}\rbr{\Omega;\C^{n}}$,
and we immediately conclude that
\[
\rbr{\chi_{\varepsilon}U_{\varepsilon,\delta},f^{\dual}}_{\Omega}=\rbr{A^{\varepsilon}D\chi_{\varepsilon}U_{\varepsilon,\delta},Du_{\varepsilon}^{\dual}}_{\Omega}-\mu\rbr{\chi_{\varepsilon}U_{\varepsilon,\delta},u_{\varepsilon}^{\dual}}_{\Omega}.
\]
As⋅a⋅result,
\begin{equation}
\begin{aligned}\hspace{2em} & \hspace{-2em}\rbr{\rbr{\op A_{\mu}^{\varepsilon}}^{-1}f-\rbr{\op A_{\mu}^{0}}^{-1}f-\varepsilon\op K_{\mu}^{\varepsilon}f,f^{\dual}}_{\Omega}\\
 & =\rbr{\chi_{\varepsilon}A^{0}Du_{0},Du_{\varepsilon}^{\dual}}_{\Omega}-\rbr{\chi_{\varepsilon}A^{\varepsilon}D\rbr{u_{0}+\varepsilon U_{\varepsilon,\delta}},Du_{\varepsilon}^{\dual}}_{\Omega}+\varepsilon\mu\rbr{\chi_{\varepsilon}U_{\varepsilon,\delta},u_{\varepsilon}^{\dual}}_{\Omega}\\
 & \quad+\rbr{\rho_{\varepsilon}\rbr{A^{0}-A^{\varepsilon}}Du_{0},Du_{\varepsilon}^{\dual}}_{\Omega}+\varepsilon\rbr{A^{\varepsilon}D\rho_{\varepsilon}\cdot U_{\varepsilon,\delta},Du_{\varepsilon}^{\dual}}_{\Omega}-\varepsilon\rbr{\rho_{\varepsilon}U_{\varepsilon,\delta},f^{\dual}}_{\Omega}.
\end{aligned}
\label{eq: First step}
\end{equation}

Let us focus on the first two terms on the right-hand side. The⋅first
one can be written, using⋅(\ref{def: Coefficient A=002070}),⋅as
\begin{equation}
\begin{aligned}\rbr{\chi_{\varepsilon}A^{0}Du_{0},Du_{\varepsilon}^{\dual}}_{\Omega} & =\rbr{\chi_{\varepsilon}A\rbr{D_{1}u_{0}+D_{2}U},D_{1}u_{\varepsilon}^{\dual}}_{\Omega\times\cell}\\
 & =\rbr{\chi_{\varepsilon}A\rbr{D_{1}u_{0,\delta}+D_{2}U_{\delta}},D_{1}u_{\varepsilon}^{\dual}}_{\Omega\times\cell}\\
 & \quad+\rbr{\chi_{\varepsilon}A\rbr{D_{1}\rbr{u_{0}-u_{0,\delta}}+D_{2}\rbr{U-U_{\delta}}},D_{1}u_{\varepsilon}^{\dual}}_{\Omega\times\cell}.
\end{aligned}
\label{eq: A=002070}
\end{equation}
For⋅the⋅second, observe that $\varepsilon DU_{\varepsilon,\delta}=\varepsilon\tau^{\varepsilon}\op S^{\varepsilon}D_{1}U_{\delta}+\tau^{\varepsilon}\op S^{\varepsilon}D_{2}U_{\delta}$,
and hence
\begin{equation}
\begin{aligned}\rbr{\chi_{\varepsilon}A^{\varepsilon}D\rbr{u_{0}+\varepsilon U_{\varepsilon,\delta}},Du_{\varepsilon}^{\dual}}_{\Omega} & =\rbr{\tau^{\varepsilon}\chi_{\varepsilon}A\op T^{\varepsilon}\rbr{D_{1}u_{0,\delta}+D_{2}U_{\delta}},D_{1}u_{\varepsilon}^{\dual}}_{\Omega\times\cell}\\
 & \quad+\rbr{\tau^{\varepsilon}\chi_{\varepsilon}A\op T^{\varepsilon}D_{1}\rbr{u_{0}-u_{0,\delta}},D_{1}u_{\varepsilon}^{\dual}}_{\Omega\times\cell}\\
 & \quad+\varepsilon\rbr{\tau^{\varepsilon}\chi_{\varepsilon}A\op T^{\varepsilon}D_{1}U_{\delta},D_{1}u_{\varepsilon}^{\dual}}_{\Omega\times\cell}\\
 & \quad+\rbr{\chi_{\varepsilon}A^{\varepsilon}\rbr{\op I-\op S^{\varepsilon}}Du_{0},Du_{\varepsilon}^{\dual}}_{\Omega}.
\end{aligned}
\label{eq: A=001D4B}
\end{equation}
We commute⋅$\op T^{\varepsilon}$ past⋅$\chi_{\varepsilon}A$ in
the first term on the right,
\begin{equation}
\begin{aligned}\rbr{\tau^{\varepsilon}\chi_{\varepsilon}A\op T^{\varepsilon}\rbr{D_{1}u_{0,\delta}+D_{2}U_{\delta}},D_{1}u_{\varepsilon}^{\dual}}_{\Omega\times\cell} & =\rbr{\tau^{\varepsilon}\op T^{\varepsilon}\chi_{\varepsilon}A\rbr{D_{1}u_{0,\delta}+D_{2}U_{\delta}},D_{1}u_{\varepsilon}^{\dual}}_{\Omega\times\cell}\\
 & \quad+\rbr{\tau^{\varepsilon}\chi_{\varepsilon}\sbr{A,\op T^{\varepsilon}}\rbr{D_{1}u_{0,\delta}+D_{2}U_{\delta}},D_{1}u_{\varepsilon}^{\dual}}_{\Omega\times\cell}\\
 & \quad-\rbr{\tau^{\varepsilon}\sbr{\rho_{\varepsilon},\op T^{\varepsilon}}A\rbr{D_{1}u_{0,\delta}+D_{2}U_{\delta}},D_{1}u_{\varepsilon}^{\dual}}_{\Omega\times\cell},
\end{aligned}
\label{eq: =00005BxA,T=00005D}
\end{equation}
and then examine the difference
\begin{equation}
\rbr{\chi_{\varepsilon}A\rbr{D_{1}u_{0,\delta}+D_{2}U_{\delta}},D_{1}u_{\varepsilon}^{\dual}}_{\Omega\times\cell}-\rbr{\tau^{\varepsilon}\op T^{\varepsilon}\chi_{\varepsilon}A\rbr{D_{1}u_{0,\delta}+D_{2}U_{\delta}},D_{1}u_{\varepsilon}^{\dual}}_{\Omega\times\cell}.\label{eq: @Proofs | The difference}
\end{equation}
Noticing that $\chi_{\varepsilon}$ vanishes near the boundary and,
moreover, so does $\op T^{\varepsilon}\chi_{\varepsilon}$, and using
Lemma~\ref{lem: =0003C4=001D4BT=001D4B}, we obtain
\begin{equation}
\begin{aligned}\rbr{\chi_{\varepsilon}A\rbr{D_{1}u_{0,\delta}+D_{2}U_{\delta}},D_{1}u_{\varepsilon}^{\dual}}_{\Omega\times\cell} & =\rbr{D_{1}^{*}\chi_{\varepsilon}A\rbr{D_{1}u_{0,\delta}+D_{2}U_{\delta}},u_{\varepsilon}^{\dual}}_{\Omega\times\cell}\\
 & =\rbr{\tau^{\varepsilon}\op T^{\varepsilon}D_{1}^{*}\chi_{\varepsilon}A\rbr{D_{1}u_{0,\delta}+D_{2}U_{\delta}},\op T^{\varepsilon}u_{\varepsilon}^{\dual}}_{\Omega\times\cell}.
\end{aligned}
\label{eq: @Proofs | The difference | 1st term}
\end{equation}
A⋅similar result for the other term in⋅(\ref{eq: @Proofs | The difference})
requires a technical lemma.
\begin{lem}
\label{lem: Differentiation of =0003C4=0003B5T=0003B5F}Fix~$\varepsilon>0$.
Let $F\in C^{0,1}\rbr{\wb{\Omega};\wt L_{p}\rbr{\cell}}^{d}$ be such
that $F\rbr{x,\wc}=0$ for⋅$x\in\rbr{\partial\Omega}_{\varepsilon}$
and $D_{2}^{*}F\rbr{x,\wc}=0$ for each⋅$x\in\Omega$ as a functional
in~$\wt W_{p}^{-1}\rbr{\cell}$. Then $D_{1}^{*}\tau^{\varepsilon}\op T^{\varepsilon}F=\tau^{\varepsilon}\op T^{\varepsilon}D_{1}^{*}F$
on⋅$C_{c}^{1}\rbr{\Omega}$, viewed as a subspace of~$C_{c}\rbr{\Omega\times\cell}$.
\end{lem}
\begin{rem}
We point out that the statement of the lemma is trivial if $F$ is
either  smooth or purely periodic. In⋅the⋅former case, the divergence
of the function~$x\mapsto F\rbr{x+\varepsilon z,\varepsilon^{-1}x}$
equals $\Div_{1}F\rbr{x+\varepsilon z,\varepsilon^{-1}x}$ whenever
$\Div_{2}F=0$ (in the strong sense). In⋅the⋅latter case, we have
that $F$ is divergence-free on the torus~$\T^{d}$ and hence on~$\R^{d}$
(see~\cite[Section~1.1]{ZhKO:1993}). Then the function~$x\mapsto\tau^{\varepsilon}\op T^{\varepsilon}F\rbr{x,z}=F\rbr{\varepsilon^{-1}x}$
is divergence-free on⋅$\R^{d}$ as well, which is exactly what the
lemma says.
\end{rem}
\begin{proof}
Let⋅$\varphi$ be a function in⋅$C_{c}^{1}\rbr{\Omega}^{n}$, extended
by zero to all of~$\R^{d}$. After a change of variables, we must
show that
\begin{equation}
\begin{aligned}\hspace{2em} & \hspace{-2em}\int_{\Omega}\int_{\cell}\abr{F\rbr{x,\varepsilon^{-1}x+y},D_{1}\varphi\rbr{x+\varepsilon y}}\dd x\dd y\\
 & =\int_{\Omega}\int_{\cell}\abr{D_{1}^{*}F\rbr{x,\varepsilon^{-1}x+y},\varphi\rbr{x+\varepsilon y}}\dd x\dd y.
\end{aligned}
\label{eq: Differentiation of =0003C4=0003B5T=0003B5F | Change of variables}
\end{equation}
Were $F\rbr{x,\wc}$  smooth, this would be nothing but the usual
integration by parts formula. But we can find a sequence of  smooth
functions~$F_{K}$ with⋅$D_{2}^{*}F_{K}=0$ that converges, in a
suitable sense, to the function~$F$, and that will complete the
proof.

If⋅$e_{k}\rbr y=e^{2\pi i\abr{y,k}}$, where⋅$k\in\Z^{d}$, then
we let $F_{K}\rbr{x,\wc}$ denote the square partial sum of the Fourier
series for~$F\rbr{x,\wc}$:
\[
F_{K}\rbr{x,\wc}=\sum_{\abs{k_{j}}\le K}\hat{F}_{k}\rbr x\hairspace e_{k}.
\]
By⋅hypothesis, $D_{2}^{*}F\rbr{x,\wc}=0$ on~$\wt W_{\smash[t]{\cramped{p^{\dual}}}}^{1}\rbr{\cell}^{n}$,~so
\[
\abr{\hat{F}_{k}\rbr x,k}=\rbr{2\pi}^{-1}\int_{\cell}\abr{F\rbr{x,y},De_{k}\rbr y}\dd y=0
\]
for⋅each~$k\in\Z^{d}$. Also notice that $D^{*}\hat{F}_{k}\rbr x$
are  the Fourier coefficients of⋅$D_{1}^{*}F\rbr{x,\wc}$. An~integration
by⋅parts then gives
\begin{equation}
\begin{aligned}\hspace{2em} & \hspace{-2em}\int_{\Omega}\int_{\cell}\abr{F_{K}\rbr{x,\varepsilon^{-1}x+y},D_{1}\varphi\rbr{x+\varepsilon y}}\dd x\dd y\\
 & =\int_{\Omega}\int_{\cell}\abr{\rbr{D_{1}^{*}F}_{K}\rbr{x,\varepsilon^{-1}x+y},\varphi\rbr{x+\varepsilon y}}\dd x\dd y.
\end{aligned}
\label{eq: Differentiation of =0003C4=0003B5T=0003B5F | Identity for Fk}
\end{equation}
Here $\rbr{D_{1}^{*}F}_{K}\rbr{x,\wc}$ is the square partial sum
of the Fourier series for~$D_{1}^{*}F\rbr{x,\wc}$.

We now show that (\ref{eq: Differentiation of =0003C4=0003B5T=0003B5F | Identity for Fk})
implies~(\ref{eq: Differentiation of =0003C4=0003B5T=0003B5F | Change of variables}).
Let⋅$G$ be a function in⋅$L_{\infty}\rbr{\R^{d};\wt L_{p}\rbr{\cell}}$,
and let⋅$G_{K}\rbr{x,\wc}$ be the square partial sum of the Fourier⋅series
for~$G\rbr{x,\wc}$. We claim that $G_{K}\to G$ in the weak\nobreakdash-$*$~topology
on⋅$C_{c}\rbr{\R^{d}\times\cell}^{*}$ as~$K\to\infty$. Indeed,
given any⋅$\psi\in C_{c}\rbr{\R^{d}\times\cell}$, the sequence of
 functions~$x\mapsto\rbr{G_{K}\rbr{x,\wc},\psi\rbr{x,\wc}}_{\cell}$
converges pointwise to the function~$x\mapsto\rbr{G\rbr{x,\wc},\psi\rbr{x,\wc}}_{\cell}$,
because $G_{K}\rbr{x,\wc}\to G\rbr{x,\wc}$ in~$L_{p}\rbr{\cell}$
(see~\cite[Theorem~4.1.8]{Gra:2014-1}). In⋅addition, all the functions
in the sequence are supported in a single compact set and are uniformly
bounded, since
\[
\begin{aligned}\abs{\rbr{G_{K}\rbr{x,\wc},\psi\rbr{x,\wc}}_{\cell}} & \lesssim\norm{G\rbr{x,\wc}}_{p,\cell}\norm{\psi\rbr{x,\wc}}_{p^{\dual},\cell}\\
 & \le\norm G_{L_{\infty}\rbr{\R^{d};L_{p}\rbr{\cell}}}\norm{\psi}_{C},
\end{aligned}
\]
where we have used the fact that $\sup_{K\in\N}\norm{G_{K}\rbr{x,\wc}}_{p,\cell}\lesssim\norm{G\rbr{x,\wc}}_{p,\cell}$
(see~\cite[Corollary~4.1.3]{Gra:2014-1}). Then $\rbr{G_{K},\psi}_{\R^{d}\times\cell}\to\rbr{G,\psi}_{\R^{d}\times\cell}$
by the Lebesgue dominated convergence theorem, and the claim follows.
 Applying this to the functions~$\rbr{x,y}\mapsto\chf_{\Omega}\rbr x\hairspace F\rbr{x,\varepsilon^{-1}x+y}$
and~$\rbr{x,y}\mapsto\chf_{\Omega}\rbr x\hairspace D_{1}^{*}F\rbr{x,\varepsilon^{-1}x+y}$
($\chf_{\Omega}$ is the characteristic function of~$\Omega$), which
obviously belong  to⋅$L_{\infty}\rbr{\R^{d};\wt L_{p}\rbr{\cell}}$,
we immediately obtain~(\ref{eq: Differentiation of =0003C4=0003B5T=0003B5F | Change of variables}).
\end{proof}

Choose a cutoff function~$\eta_{\varepsilon}\in C^{0,1}\rbr{\wb{\Omega}}$
satisfying⋅$\eta_{\varepsilon}\restrict{\rbr{\supp\chi_{\varepsilon}}_{\varepsilon}}=1$.
By⋅definition of⋅$U_{\delta}$, the second term in⋅(\ref{eq: @Proofs | The difference})⋅is
\[
\rbr{\tau^{\varepsilon}\op T^{\varepsilon}\chi_{\varepsilon}A\rbr{I+D_{2}N}D_{1}u_{0,\delta},D_{1}\eta_{\varepsilon}u_{\varepsilon}^{\dual}}_{\Omega\times\cell}.
\]
Assume for the moment that  $\eta_{\varepsilon}u_{\varepsilon}^{\dual}\in C_{c}^{1}\rbr{\Omega}^{n}$
and recall from⋅(\ref{def: N}) that, for each fixed⋅$x\in\Omega$,
$D_{2}^{*}A\rbr{x,\wc}\hairspace\rbr{I+D_{2}N\rbr{x,\wc}}\hairspace Du_{0,\delta}\rbr x=0$
on~$\wt W_{\smash[t]{\cramped{p^{\dual}}}}^{1}\rbr{\cell}^{n}$. Then Lemma~\ref{lem: Differentiation of =0003C4=0003B5T=0003B5F}
 tells us that
\[
\rbr{\tau^{\varepsilon}\op T^{\varepsilon}\chi_{\varepsilon}A\rbr{I+D_{2}N}D_{1}u_{0,\delta},D_{1}\eta_{\varepsilon}u_{\varepsilon}^{\dual}}_{\Omega\times\cell}=\rbr{\tau^{\varepsilon}\op T^{\varepsilon}D_{1}^{*}\chi_{\varepsilon}A\rbr{I+D_{2}N}D_{1}u_{0,\delta},\eta_{\varepsilon}u_{\varepsilon}^{\dual}}_{\Omega\times\cell}.
\]
But the form
\[
\eta_{\varepsilon}u_{\varepsilon}^{\dual}\mapsto\rbr{\tau^{\varepsilon}\op T^{\varepsilon}\chi_{\varepsilon}A\rbr{I+D_{2}N}D_{1}u_{0,\delta},D_{1}\eta_{\varepsilon}u_{\varepsilon}^{\dual}}_{\Omega\times\cell}
\]
is⋅continuous on⋅$\mathring{W}_{p^{\dual}}^{1}\rbr{\Omega}^{n}$
and the form
\[
\eta_{\varepsilon}u_{\varepsilon}^{\dual}\mapsto\rbr{\tau^{\varepsilon}\op T^{\varepsilon}D_{1}^{*}\chi_{\varepsilon}A\rbr{I+D_{2}N}D_{1}u_{0,\delta},\eta_{\varepsilon}u_{\varepsilon}^{\dual}}_{\Omega\times\cell}
\]
is⋅continuous on⋅$L_{p^{\dual}}\rbr{\Omega}^{n}$ (by⋅Lemma~\ref{lem: =0003C4=001D4BT=001D4B}
and the hypothesis~(\ref{assu: N is Lipschitz})), so the equality
\begin{equation}
\rbr{\tau^{\varepsilon}\op T^{\varepsilon}\chi_{\varepsilon}A\rbr{D_{1}u_{0,\delta}+D_{2}U_{\delta}},D_{1}\eta_{\varepsilon}u_{\varepsilon}^{\dual}}_{\Omega\times\cell}=\rbr{\tau^{\varepsilon}\op T^{\varepsilon}D_{1}^{*}\chi_{\varepsilon}A\rbr{D_{1}u_{0,\delta}+D_{2}U_{\delta}},\eta_{\varepsilon}u_{\varepsilon}^{\dual}}_{\Omega\times\cell}\label{eq: @Proofs | The difference | 2nd term}
\end{equation}
holds, in fact, for any~$u_{\varepsilon}^{\dual}\in W_{\smash[t]{\cramped{p^{\dual}}}}^{1}\rbr{\Omega}^{n}$.
Recalling that $\eta_{\varepsilon}=1$ on $\rbr{\supp\chi_{\varepsilon}}_{\varepsilon}$
and combining (\ref{eq: @Proofs | The difference | 1st term}) with⋅(\ref{eq: @Proofs | The difference | 2nd term}),
we see that
\begin{equation}
\begin{aligned}\hspace{2em} & \hspace{-2em}\rbr{\chi_{\varepsilon}A\rbr{D_{1}u_{0,\delta}+D_{2}U_{\delta}},D_{1}u_{\varepsilon}^{\dual}}_{\Omega\times\cell}-\rbr{\tau^{\varepsilon}\op T^{\varepsilon}\chi_{\varepsilon}A\rbr{D_{1}u_{0,\delta}+D_{2}U_{\delta}},D_{1}u_{\varepsilon}^{\dual}}_{\Omega\times\cell}\\
 & =\rbr{\tau^{\varepsilon}\op T^{\varepsilon}D_{1}^{*}\chi_{\varepsilon}A\rbr{D_{1}u_{0,\delta}+D_{2}U_{\delta}},\rbr{\op T^{\varepsilon}-\op I}u_{\varepsilon}^{\dual}}_{\Omega\times\cell}.
\end{aligned}
\label{eq: @Proofs | The difference after lemma}
\end{equation}

Putting  together⋅(\ref{eq: First step})–(\ref{eq: =00005BxA,T=00005D})
and~(\ref{eq: @Proofs | The difference after lemma}), we arrive
at the operator identity
\begin{equation}
\rbr{\op A_{\mu}^{\varepsilon}}^{-1}-\rbr{\op A_{\mu}^{0}}^{-1}-\varepsilon\op K_{\mu}^{\varepsilon}\restrict{L_{p}\rbr{\Omega}^{n}}=\op I_{\mu}^{\varepsilon}+\op D_{\mu}^{\varepsilon}+\op B_{\mu}^{\varepsilon}\label{eq: The identity}
\end{equation}
that effectively splits the problem into the interior parts, given⋅by
\begin{equation}
\begin{aligned}\rbr{\op I_{\mu}^{\varepsilon}f,f^{\dual}}_{\Omega} & =\rbr{\tau^{\varepsilon}\op T^{\varepsilon}\chi_{\varepsilon}D_{1}^{*}A\rbr{D_{1}u_{0,\delta}+D_{2}U_{\delta}},\rbr{\op T^{\varepsilon}-\op I}u_{\varepsilon}^{\dual}}_{\Omega\times\cell}\\
 & \quad-\rbr{\tau^{\varepsilon}\chi_{\varepsilon}\sbr{A,\op T^{\varepsilon}}\rbr{D_{1}u_{0,\delta}+D_{2}U_{\delta}},D_{1}u_{\varepsilon}^{\dual}}_{\Omega\times\cell}\\
 & \quad-\varepsilon\rbr{\tau^{\varepsilon}\chi_{\varepsilon}A\op T^{\varepsilon}D_{1}U_{\delta},D_{1}u_{\varepsilon}^{\dual}}_{\Omega\times\cell}\\
 & \quad-\rbr{\chi_{\varepsilon}A^{\varepsilon}\rbr{\op I-\op S^{\varepsilon}}Du_{0},Du_{\varepsilon}^{\dual}}_{\Omega}\\
 & \quad+\varepsilon\mu\rbr{\chi_{\varepsilon}U_{\varepsilon,\delta},u_{\varepsilon}^{\dual}}_{\Omega}
\end{aligned}
\label{def: Interior part}
\end{equation}
and
\begin{equation}
\begin{aligned}\rbr{\op D_{\mu}^{\varepsilon}f,f^{\dual}}_{\Omega} & =\rbr{\chi_{\varepsilon}A\rbr{D_{1}\rbr{u_{0}-u_{0,\delta}}+D_{2}\rbr{U-U_{\delta}}},D_{1}u_{\varepsilon}^{\dual}}_{\Omega\times\cell}\\
 & \quad-\rbr{\tau^{\varepsilon}\chi_{\varepsilon}A\op T^{\varepsilon}D_{1}\rbr{u_{0}-u_{0,\delta}},D_{1}u_{\varepsilon}^{\dual}}_{\Omega\times\cell},
\end{aligned}
\label{def: Interior part, discrepancy}
\end{equation}
and the boundary part, given⋅by
\begin{equation}
\begin{aligned}\rbr{\op B_{\mu}^{\varepsilon}f,f^{\dual}}_{\Omega} & =\rbr{\op T^{\varepsilon}\rbr{D_{1}\rho_{\varepsilon}}^{*}\cdot\tau^{\varepsilon}\op T^{\varepsilon}A\rbr{D_{1}u_{0,\delta}+D_{2}U_{\delta}},\rbr{\op T^{\varepsilon}-\op I}u_{\varepsilon}^{\dual}}_{\Omega\times\cell}\\
 & \quad+\rbr{\tau^{\varepsilon}\sbr{\rho_{\varepsilon},\op T^{\varepsilon}}A\rbr{D_{1}u_{0,\delta}+D_{2}U_{\delta}},D_{1}u_{\varepsilon}^{\dual}}_{\Omega\times\cell}\\
 & \quad+\rbr{\rho_{\varepsilon}\rbr{A^{0}-A^{\varepsilon}}Du_{0},Du_{\varepsilon}^{\dual}}_{\Omega}\\
 & \quad+\varepsilon\rbr{A^{\varepsilon}D\rho_{\varepsilon}\cdot U_{\varepsilon,\delta},Du_{\varepsilon}^{\dual}}_{\Omega}\\
 & \quad-\varepsilon\rbr{\rho_{\varepsilon}U_{\varepsilon,\delta},f^{\dual}}_{\Omega}.
\end{aligned}
\label{def: Boundary part}
\end{equation}
This is the resolvent identity that we seek.

\subsection{Proof of the main results}

We now estimate each operator in the identity and begin with the
interior part. In⋅what follows, we will frequently use  the fact
that $\chi_{\varepsilon}$ vanishes on⋅$\rbr{\partial\Omega}_{2\varepsilon}\cap\Omega$
and⋅$\rho_{\varepsilon}$ is supported in~$\rbr{\partial\Omega}_{3\varepsilon}\cap\Omega$.
\begin{lem}
\label{lem: Ie}For⋅any⋅$\varepsilon\in\set E_{\mu}$ and⋅$f\in L_{p}\rbr{\Omega}^{n}$,
\begin{equation}
\norm{\op I_{\mu}^{\varepsilon}f}_{1,p,\Omega}\lesssim\varepsilon^{s}\norm f_{p,\Omega}.\label{est: Ie}
\end{equation}
\end{lem}
\begin{proof}
By⋅Lemmas~\ref{lem: =0003C4=001D4BT=001D4B} and~\ref{lem: T=001D4B-I},
\[
\begin{aligned}\hspace{2em} & \hspace{-2em}\bigabs{\rbr{\tau^{\varepsilon}\op T^{\varepsilon}\chi_{\varepsilon}D_{1}^{*}A\rbr{D_{1}u_{0,\delta}+D_{2}U_{\delta}},\rbr{\op T^{\varepsilon}-\op I}u_{\varepsilon}^{\dual}}_{\Omega\times\cell}}\\
 & \le\norm{\tau^{\varepsilon}\op T^{\varepsilon}\chi_{\varepsilon}D_{1}^{*}A\rbr{D_{1}u_{0,\delta}+D_{2}U_{\delta}}}_{p,\rbr{\supp\chi_{\varepsilon}}_{\varepsilon}\times\cell}\norm{\rbr{\op T^{\varepsilon}-\op I}u_{\varepsilon}^{\dual}}_{p^{\dual},\rbr{\supp\chi_{\varepsilon}}_{\varepsilon}\times\cell}\\
 & \lesssim\varepsilon\bigrbr{\norm{Du_{0,\delta}}_{1,p,\Omega}+\norm{D_{1}D_{2}U_{\delta}}_{p,\Omega\times\cell}+\norm{D_{2}U_{\delta}}_{p,\Omega\times\cell}}\norm{Du_{\varepsilon}^{\dual}}_{p^{\dual},\Omega}.
\end{aligned}
\]
Next, observe that⋅$\tau^{\varepsilon}\sbr{A,\op T^{\varepsilon}}=\tau^{\varepsilon}\rbr{\op I-\op T^{\varepsilon}}A\cdot\tau^{\varepsilon}\op T^{\varepsilon}$,
which, with Lemma~\ref{lem: =0003C4=001D4BT=001D4B}, gives
\[
\begin{aligned}\hspace{2em} & \hspace{-2em}\bigabs{\rbr{\tau^{\varepsilon}\chi_{\varepsilon}\sbr{A,\op T^{\varepsilon}}\rbr{D_{1}u_{0,\delta}+D_{2}U_{\delta}},D_{1}u_{\varepsilon}^{\dual}}_{\Omega\times\cell}}\\
 & \le\norm{\rbr{\op I-\op T^{\varepsilon}}A}_{L_{\infty}}\norm{\tau^{\varepsilon}\op T^{\varepsilon}\rbr{D_{1}u_{0,\delta}+D_{2}U_{\delta}}}_{p,\supp\chi_{\varepsilon}\times\cell}\norm{D_{1}u_{\varepsilon}^{\dual}}_{p^{\dual},\Omega\times\cell}\\
 & \lesssim\varepsilon\bigrbr{\norm{Du_{0,\delta}}_{p,\Omega}+\norm{D_{2}U_{\delta}}_{p,\Omega\times\cell}}\norm{Du_{\varepsilon}^{\dual}}_{p^{\dual},\Omega}.
\end{aligned}
\]
By⋅Lemma~\ref{lem: =0003C4=001D4BT=001D4B} again, we see that
\[
\begin{aligned}\varepsilon\bigabs{\rbr{\tau^{\varepsilon}\chi_{\varepsilon}A\op T^{\varepsilon}D_{1}U_{\delta},D_{1}u_{\varepsilon}^{\dual}}_{\Omega\times\cell}} & \le\varepsilon\norm A_{L_{\infty}}\norm{\tau^{\varepsilon}\op T^{\varepsilon}D_{1}U_{\delta}}_{p,\supp\chi_{\varepsilon}\times\cell}\norm{D_{1}u_{\varepsilon}^{\dual}}_{p^{\dual},\Omega\times\cell}\\
 & \lesssim\varepsilon\norm{D_{1}U_{\delta}}_{p,\Omega\times\cell}\norm{Du_{\varepsilon}^{\dual}}_{p^{\dual},\Omega},
\end{aligned}
\]
while Lemmas~\ref{lem: =0003C4=001D4BS=001D4B} and~\ref{lem: S=001D4B-I}
show that
\[
\begin{aligned}\bigabs{\rbr{\chi_{\varepsilon}A^{\varepsilon}\rbr{\op I-\op S^{\varepsilon}}Du_{0},Du_{\varepsilon}^{\dual}}_{\Omega}} & \le\norm A_{L_{\infty}}\norm{\rbr{\op I-\op S^{\varepsilon}}Du_{0}}_{p,\supp\chi_{\varepsilon}}\norm{Du_{\varepsilon}^{\dual}}_{p^{\dual},\Omega}\\
 & \lesssim\bigrbr{\varepsilon\norm{DDu_{0,\delta}}_{p,\Omega}+\norm{D\rbr{u_{0}-u_{0,\delta}}}_{p,\Omega}}\norm{Du_{\varepsilon}^{\dual}}_{p^{\dual},\Omega}
\end{aligned}
\]
and
\[
\varepsilon\bigabs{\rbr{\chi_{\varepsilon}U_{\varepsilon,\delta},u_{\varepsilon}^{\dual}}_{\Omega}}\le\varepsilon\norm{U_{\varepsilon,\delta}}_{p,\supp\chi_{\varepsilon}}\norm{u_{\varepsilon}^{\dual}}_{p^{\dual},\Omega}\lesssim\varepsilon\norm{U_{\delta}}_{p,\Omega\times\cell}\norm{u_{\varepsilon}^{\dual}}_{p^{\dual},\Omega}.
\]
Combining these estimates with Lemmas~\ref{lem: R0_mu(delta) and K_mu(delta) are bounded}
and⋅\ref{lem: R0_mu(delta) and K_mu(delta) convergence rates} and
the inequality~(\dualref{est: Original operator | Norm of (A=001D4B-=0003BC)=00207B=0000B9})
yields
\[
\bigabs{\rbr{\op I_{\mu}^{\varepsilon}f,f^{\dual}}_{\Omega}}\lesssim\varepsilon^{s}\norm f_{p,\Omega}\norm{f^{\dual}}_{-1,p^{\dual},\Omega}^{*}.
\]
This is what we wished to prove.
\end{proof}
\begin{lem}
\label{lem: De}For⋅any⋅$\varepsilon\in\set E_{\mu}$ and⋅$f\in L_{p}\rbr{\Omega}^{n}$,
\begin{equation}
\norm{\op D_{\mu}^{\varepsilon}f}_{1,p,\Omega}\lesssim\varepsilon^{s}\norm f_{p,\Omega}.\label{est: De}
\end{equation}
\end{lem}
\begin{proof}
We have that 
\[
\begin{aligned}\hspace{2em} & \hspace{-2em}\bigabs{\rbr{\chi_{\varepsilon}A\rbr{D_{1}\rbr{u_{0}-u_{0,\delta}}+D_{2}\rbr{U-U_{\delta}}},D_{1}u_{\varepsilon}^{\dual}}_{\Omega\times\cell}}\le\\
 & \le\norm A_{L_{\infty}}\bigrbr{\norm{D_{1}\rbr{u_{0}-u_{0,\delta}}}_{p,\Omega\times\cell}+\norm{D_{2}\rbr{U-U_{\delta}}}_{p,\Omega\times\cell}}\norm{D_{1}u_{\varepsilon}^{\dual}}_{p^{\dual},\Omega\times\cell}\\
 & \lesssim\bigrbr{\norm{D\rbr{u_{0}-u_{0,\delta}}}_{p,\Omega}+\norm{D_{2}\rbr{U-U_{\delta}}}_{p,\Omega\times\cell}}\norm{Du_{\varepsilon}^{\dual}}_{p^{\dual},\Omega}
\end{aligned}
\]
and
\[
\begin{aligned}\hspace{2em} & \hspace{-2em}\bigabs{\rbr{\tau^{\varepsilon}\chi_{\varepsilon}A\op T^{\varepsilon}D_{1}\rbr{u_{0}-u_{0,\delta}},D_{1}u_{\varepsilon}^{\dual}}_{\Omega\times\cell}}\\
 & \le\norm A_{L_{\infty}}\norm{\tau^{\varepsilon}\op T^{\varepsilon}D_{1}\rbr{u_{0}-u_{0,\delta}}}_{p,\supp\chi_{\varepsilon}\times\cell}\norm{D_{1}u_{\varepsilon}^{\dual}}_{p^{\dual},\Omega\times\cell}\\
 & \lesssim\norm{D\rbr{u_{0}-u_{0,\delta}}}_{p,\Omega}\norm{Du_{\varepsilon}^{\dual}}_{p^{\dual},\Omega}
\end{aligned}
\]
(according to Lemma~\ref{lem: =0003C4=001D4BT=001D4B}). Hence,
by Lemma~\ref{lem: R0_mu(delta) and K_mu(delta) convergence rates}
and the inequality~(\dualref{est: Original operator | Norm of (A=001D4B-=0003BC)=00207B=0000B9}),
\[
\bigabs{\rbr{\op D_{\mu}^{\varepsilon}f,f^{\dual}}_{\Omega}}\lesssim\varepsilon^{s}\norm f_{p,\Omega}\norm{f^{\dual}}_{-1,p^{\dual},\Omega}^{*}.\qedhere
\]
\end{proof}
It⋅remains to estimate the boundary operator~$\op B_{\mu}^{\varepsilon}$,
for which we need the following.
\begin{lem}
\label{lem: Be wo last term}For⋅any⋅$\varepsilon\in\set E_{\mu}$,
$f\in L_{p}\rbr{\Omega}^{n}$ and⋅$f^{\dual}\in\rbr{W_{p}^{1}\rbr{\Omega}^{n}}^{*}$,
\begin{equation}
\bigabs{\rbr{\op B_{\mu}^{\varepsilon}f+\varepsilon\rho_{\varepsilon}U_{\varepsilon,\delta},f^{\dual}}_{\Omega}}\lesssim\varepsilon^{s/p}\norm f_{p,\Omega}\norm{Du_{\varepsilon}^{\dual}}_{p^{\dual},\rbr{\partial\Omega}_{5\varepsilon}\cap\Omega}.\label{est: Be wo last term}
\end{equation}
\end{lem}
\begin{proof}
Arguing as in the proof of Lemma~\ref{lem: Ie}, we easily find
that
\[
\begin{aligned}\hspace{2em} & \hspace{-2em}\bigabs{\rbr{\op T^{\varepsilon}\rbr{D_{1}\rho_{\varepsilon}}^{*}\cdot\tau^{\varepsilon}\op T^{\varepsilon}A\rbr{D_{1}u_{0,\delta}+D_{2}U_{\delta}},\rbr{\op T^{\varepsilon}-\op I}u_{\varepsilon}^{\dual}}_{\Omega\times\cell}}\\
 & \lesssim\bigrbr{\norm{Du_{0,\delta}}_{p,\supp D\rho_{\varepsilon}}+\norm{D_{2}U_{\delta}}_{p,\supp D\rho_{\varepsilon}\times\cell}}\norm{Du_{\varepsilon}^{\dual}}_{p^{\dual},\rbr{\supp D\rho_{\varepsilon}}_{2\varepsilon}}
\end{aligned}
\]
and
\[
\begin{aligned}\hspace{2em} & \hspace{-2em}\bigabs{\rbr{\tau^{\varepsilon}\sbr{\rho_{\varepsilon},\op T^{\varepsilon}}A\rbr{D_{1}u_{0,\delta}+D_{2}U_{\delta}},D_{1}u_{\varepsilon}^{\dual}}_{\Omega\times\cell}}\\
 & \lesssim\bigrbr{\norm{Du_{0,\delta}}_{p,\rbr{\supp D\rho_{\varepsilon}}_{2\varepsilon}}+\norm{D_{2}U_{\delta}}_{p,\rbr{\supp D\rho_{\varepsilon}}_{2\varepsilon}\times\cell}}\norm{Du_{\varepsilon}^{\dual}}_{p^{\dual},\rbr{\supp D\rho_{\varepsilon}}_{\varepsilon}}.
\end{aligned}
\]
Likewise,
\[
\begin{aligned}\hspace{2em} & \hspace{-2em}\bigabs{\rbr{\rho_{\varepsilon}\rbr{A^{0}-A^{\varepsilon}}Du_{0},Du_{\varepsilon}^{\dual}}_{\Omega}}\\
 & \lesssim\bigrbr{\norm{D\rbr{u_{0}-u_{0,\delta}}}_{p,\Omega}+\norm{Du_{0,\delta}}_{p,\supp\rho_{\varepsilon}}}\norm{Du_{\varepsilon}^{\dual}}_{p^{\dual},\supp\rho_{\varepsilon}}
\end{aligned}
\]
and
\[
\varepsilon\bigabs{\rbr{A^{\varepsilon}D\rho_{\varepsilon}\cdot U_{\varepsilon,\delta},Du_{\varepsilon}^{\dual}}_{\Omega}}\lesssim\norm{U_{\delta}}_{p,\rbr{\supp D\rho_{\varepsilon}}_{\varepsilon}\times\cell}\norm{Du_{\varepsilon}^{\dual}}_{p^{\dual},\supp D\rho_{\varepsilon}}.
\]
Thus
\[
\begin{aligned}\bigabs{\rbr{\op B_{\mu}^{\varepsilon}f+\varepsilon\rho_{\varepsilon}U_{\varepsilon,\delta},f^{\dual}}_{\Omega}} & \lesssim\bigl(\norm{Du_{0,\delta}}_{p,\rbr{\partial\Omega}_{5\varepsilon}\cap\Omega}+\norm{D\rbr{u_{0}-u_{0,\delta}}}_{p,\Omega}\\
 & \hphantom{{}\lesssim\bigl(\norm{Du_{0,\delta}}_{p,\rbr{\partial\Omega}_{5\varepsilon}\cap\Omega}+\norm{D\rbr{u_{0}-u_{0,\delta}}}_{p,\Omega}}\mathllap{{}+\norm{U_{\delta}}_{1_{2},p,\rbr{\partial\Omega}_{5\varepsilon}\cap\Omega\times\cell}}\bigr)\norm{Du_{\varepsilon}^{\dual}}_{p^{\dual},\rbr{\partial\Omega}_{5\varepsilon}\cap\Omega}.
\end{aligned}
\]
Now notice that, for all⋅$v\in W_{p}^{1}\rbr{\Omega}$ and⋅$\varepsilon\le1$,
\begin{equation}
\norm v_{p,\rbr{\partial\Omega}_{5\varepsilon}\cap\Omega}\lesssim\varepsilon^{1/p}\norm v_{1,p,\Omega}^{1/p}\norm v_{p,\Omega}^{1/p^{\dual}}\lesssim\varepsilon^{s/p}\bigrbr{\varepsilon^{1-s}\norm{Dv}_{p,\Omega}+\norm v_{p,\Omega}}\label{est: Norm in the neighborhood, with s}
\end{equation}
by⋅Lemma~\ref{lem: Norm in the neighborhood} and Young's inequality.
This together with Lemmas~\ref{lem: R0_mu(delta) and K_mu(delta) are bounded}
and⋅\ref{lem: R0_mu(delta) and K_mu(delta) convergence rates} complete
the proof.
\end{proof}
Now we are ready to obtain a bound on~$\op B_{\mu}^{\varepsilon}$.
\begin{lem}
\label{lem: Be in Wp1}For⋅any⋅$\varepsilon\in\set E_{\mu}$ and⋅$f\in L_{p}\rbr{\Omega}^{n}$,
\[
\norm{\op B_{\mu}^{\varepsilon}f}_{1,p,\Omega}\lesssim\varepsilon^{s/p}\norm f_{p,\Omega}.
\]
\end{lem}
\begin{proof}
It⋅follows from⋅(\dualref{est: Original operator | Norm of (A=001D4B-=0003BC)=00207B=0000B9})
and⋅Lemma~\ref{lem: Be wo last term} that
\[
\norm{\op B_{\mu}^{\varepsilon}f}_{1,p,\Omega}\lesssim\varepsilon\norm{\rho_{\varepsilon}U_{\varepsilon,\delta}}_{1,p,\Omega}+\varepsilon^{s/p}\norm f_{p,\Omega}.
\]
Next, $\varepsilon DU_{\varepsilon,\delta}=\varepsilon\tau^{\varepsilon}\op S^{\varepsilon}D_{1}U_{\delta}+\tau^{\varepsilon}\op S^{\varepsilon}D_{2}U_{\delta}$,
and then, by⋅Lemma~\ref{lem: =0003C4=001D4BS=001D4B},
\[
\begin{aligned}\varepsilon\norm{\rho_{\varepsilon}U_{\varepsilon,\delta}}_{1,p,\Omega} & \lesssim\varepsilon\norm{U_{\varepsilon,\delta}}_{1,p,\supp\rho_{\varepsilon}}+\norm{U_{\varepsilon,\delta}}_{p,\supp D\rho_{\varepsilon}}\\
 & \lesssim\varepsilon\norm{U_{\delta}}_{1_{1},p,\Omega_{1}\times\cell}+\norm{U_{\delta}}_{1_{2},p,\rbr{\supp\rho_{\varepsilon}}_{\varepsilon}\times\cell}
\end{aligned}
\]
(recall that $U_{\delta}$ is extended to all of⋅$\R^{d}$ and thus
is well-defined on~$\Omega_{1}$). We decompose⋅$\rbr{\supp\rho_{\varepsilon}}_{\varepsilon}$
into $\rbr{\supp\rho_{\varepsilon}}_{\varepsilon}\cap\Omega$ and⋅$\rbr{\supp\rho_{\varepsilon}}_{\varepsilon}\!\setminus\Omega$
and apply Lemma~\ref{lem: Norm in the neighborhood} with, respectively,⋅$\Sigma=\Omega$
and⋅$\Sigma=\Omega_{1}\!\setminus\wb{\Omega}$ (cf.~(\ref{est: Norm in the neighborhood, with s}))
 to get
\begin{equation}
\varepsilon\norm{\rho_{\varepsilon}U_{\varepsilon,\delta}}_{1,p,\Omega}\lesssim\varepsilon^{s/p}\bigrbr{\varepsilon^{1-s}\norm{D_{1}U_{\delta}}_{1_{2},p,\Omega_{1}\times\cell}+\norm{U_{\delta}}_{1_{2},p,\Omega_{1}\times\cell}}.\label{est: Norm of DUe near boundary}
\end{equation}
Hence
\[
\varepsilon\norm{\rho_{\varepsilon}U_{\varepsilon,\delta}}_{1,p,\Omega}\lesssim\varepsilon^{s/p}\norm f_{p,\Omega}
\]
by⋅Lemma~\ref{lem: R0_mu(delta) and K_mu(delta) are bounded}, and
the proof is complete.
\end{proof}
We are finally in a position  to prove Theorem~\ref{thm: Convergence and Approximation with 1st corrector}.
\begin{proof}[Proof of Theorem~\textup{\ref{thm: Convergence and Approximation with 1st corrector}}]
The⋅identity~(\ref{eq: The identity}), together with Lemmas~\ref{lem: Ie},
\ref{lem: De} and~\ref{lem: Be in Wp1} immediately implies that
\begin{equation}
\norm{\rbr{\op A_{\mu}^{\varepsilon}}^{-1}f-\rbr{\op A_{\mu}^{0}}^{-1}f-\varepsilon\op K_{\mu}^{\varepsilon}f}_{1,p,\Omega}\lesssim\varepsilon^{s/p}\norm f_{p,\Omega}.\label{est: Wp1-bound}
\end{equation}
In⋅particular,⋅(\ref{est: Approximation with 1st corrector}) holds.
The⋅$L_{q}$\nobreakdash-bound~(\ref{est: Convergence}) comes
from⋅(\ref{est: Wp1-bound}) as well, since
\[
\norm{\rbr{\op A_{\mu}^{\varepsilon}}^{-1}f-\rbr{\op A_{\mu}^{0}}^{-1}f}_{q,\Omega}\lesssim\norm{\rbr{\op A_{\mu}^{\varepsilon}}^{-1}f-\rbr{\op A_{\mu}^{0}}^{-1}f-\varepsilon\op K_{\mu}^{\varepsilon}f}_{1,p,\Omega}+\varepsilon\norm{\op K_{\mu}^{\varepsilon}f}_{q,\Omega}
\]
according to the Sobolev⋅embedding theorem, and the terms on the right
are estimated by using⋅(\ref{est: Norm of K=001D4B with q}) and~(\ref{est: Wp1-bound}).
\end{proof}
\begin{proof}[Proof of Corollary~\textup{\ref{cor: Corollary 1}}]
From⋅(\ref{est: Wp1-bound}) and the fact that $W_{p}^{1}\rbr{\Omega}^{n}$
is continuously embedded in⋅$W_{p}^{r}\rbr{\Omega}^{n}$, we conclude
that
\[
\Seminorm{\rbr{\op A_{\mu}^{\varepsilon}}^{-1}f-\rbr{\op A_{\mu}^{0}}^{-1}f-\varepsilon\op K_{\mu}^{\varepsilon}f}_{r,p,\Omega}\lesssim\varepsilon^{s/p}\norm f_{p,\Omega}.
\]
On⋅the⋅other hand, interpolation between the $W_{p}^{1}$\nobreakdash-
and⋅$L_{p}$\nobreakdash-bounds  in~(\ref{est: Norm of K=001D4B})
gives
\[
\varepsilon^{r}\Seminorm{\op K_{\mu}^{\varepsilon}f}_{r,p,\Omega}\lesssim\norm f_{p,\Omega},
\]
so⋅(\ref{est: Corollary 1}) follows.
\end{proof}
The⋅estimate for the boundary part, the operator~$\op B_{\mu}^{\varepsilon}$,
was the worst one. Knowing that $\rbr{\op A_{\mu}^{\varepsilon}}^{\dual}$
satisfies the hypotheses of Theorem~\dualref{thm: Convergence and Approximation with 1st corrector},
we can get a better  bound on the operator norm of⋅$\op B_{\mu}^{\varepsilon}$
on⋅$L_{p}\rbr{\Omega}^{n}$ by estimating⋅$Du_{\varepsilon}^{\dual}$
near the boundary more carefully (this idea is due to Griso, see~\cite[Lemma~3.1]{Gri:2006}).
\begin{lem}
\label{lem: Be in Lp}Suppose that the hypotheses of Theorem~\textup{\dualref{thm: Convergence and Approximation with 1st corrector}}
hold. Then for any⋅$\varepsilon\in\set E_{\mu}$ and⋅$f\in L_{p}\rbr{\Omega}^{n}$
\[
\norm{\op B_{\mu}^{\varepsilon}f}_{p,\Omega}\lesssim\varepsilon^{s/p+s^{\dual}{\mkern-4mu }/p^{\dual}}\norm f_{p,\Omega}.
\]
\end{lem}
\begin{proof}
We know from Lemma~\ref{lem: Be wo last term} and the estimate~(\ref{est: Norm of K=001D4B})
that
\[
\bigabs{\rbr{\op B_{\mu}^{\varepsilon}f,f^{\dual}}_{\Omega}}\lesssim\varepsilon^{s/p+s^{\dual}{\mkern-4mu }/p^{\dual}}\norm f_{p,\Omega}\bigrbr{\varepsilon^{-s^{\dual}{\mkern-4mu }/p^{\dual}}\norm{Du_{\varepsilon}^{\dual}}_{p^{\dual},\rbr{\partial\Omega}_{5\varepsilon}\cap\Omega}+\norm{f^{\dual}}_{p^{\dual},\Omega}}.
\]
Theorem~\dualref{thm: Convergence and Approximation with 1st corrector}
tells us that $Du_{\varepsilon}^{\dual}$ can be approximated by⋅$D\rbr{u_{0}^{\dual}+\varepsilon U_{\varepsilon,\delta}^{\dual}}$,
yielding
\[
\norm{Du_{\varepsilon}^{\dual}}_{p^{\dual},\rbr{\partial\Omega}_{5\varepsilon}\cap\Omega}\lesssim\norm{Du_{0}^{\dual}}_{p^{\dual},\rbr{\partial\Omega}_{5\varepsilon}\cap\Omega}+\varepsilon\norm{DU_{\varepsilon,\delta}^{\dual}}_{p^{\dual},\rbr{\partial\Omega}_{5\varepsilon}\cap\Omega}+\varepsilon^{s^{\dual}{\mkern-4mu }/p^{\dual}}\norm{f^{\dual}}_{p^{\dual},\Omega}.
\]
Next,
\[
\norm{Du_{0}^{\dual}}_{p^{\dual},\rbr{\partial\Omega}_{5\varepsilon}\cap\Omega}\lesssim\norm{D\rbr{u_{0}^{\dual}-u_{0,\delta}^{\dual}}}_{p^{\dual},\Omega}+\norm{Du_{0,\delta}^{\dual}}_{p^{\dual},\rbr{\partial\Omega}_{5\varepsilon}\cap\Omega}
\]
and
\[
\varepsilon\norm{DU_{\varepsilon,\delta}^{\dual}}_{p^{\dual},\rbr{\partial\Omega}_{5\varepsilon}\cap\Omega}\lesssim\varepsilon^{s^{\dual}{\mkern-4mu }/p^{\dual}}\bigrbr{\varepsilon^{1-s^{\dual}}\norm{D_{1}U_{\delta}^{\dual}}_{1_{2},p,\Omega_{1}\times\cell}+\norm{D_{2}U_{\delta}^{\dual}}_{p,\Omega_{1}\times\cell}}
\]
(cf.~(\ref{est: Norm of DUe near boundary})), whence
\[
\norm{Du_{\varepsilon}^{\dual}}_{p^{\dual},\rbr{\partial\Omega}_{5\varepsilon}\cap\Omega}\lesssim\varepsilon^{s^{\dual}{\mkern-4mu }/p^{\dual}}\norm{f^{\dual}}_{p^{\dual},\Omega}
\]
by⋅Lemmas~\dualref{lem: R0_mu(delta) and K_mu(delta) are bounded},⋅\dualref{lem: R0_mu(delta) and K_mu(delta) convergence rates}
and⋅\ref{lem: Norm in the neighborhood}, as required.
\end{proof}
\begin{proof}[Proof of Theorem~\textup{\ref{thm: Convergence with e}}]
This is immediate from the identity~(\ref{eq: The identity})
and the estimate~(\ref{est: Norm of K=001D4B}) and Lemmas~\ref{lem: Ie},
\ref{lem: De} and~\ref{lem: Be in Lp}.
\end{proof}
As⋅we have seen, the interior terms in⋅(\ref{eq: The identity})
are of order~$\varepsilon^{s}$ even in the⋅$W_{p}^{1}$\nobreakdash-norm,
unlike the boundary term.  To~go further, we establish an ``interior''
operator identity, which is  similar to⋅(\ref{eq: The identity})
but involves no boundary terms.

So⋅let $\chi^{\prime}\in C^{0,1}\rbr{\wb{\Omega}}$ with $\chi^{\prime}=0$
in⋅$\rbr{\partial\Omega}_{\sigma}$ for some~$\sigma>0$. Define
the linear operator~$\op P^{\varepsilon}\colon W_{p}^{1}\rbr{\Omega}^{n}\to\rbr{C_{c}^{\infty}\rbr{\Omega}^{n}}^{*}$
associated with the form~$\rbr{u,v}\mapsto\rbr{A^{\varepsilon}Du,Dv}_{\Omega}$
and set~$\op P_{\mu}^{\varepsilon}=\op P^{\varepsilon}-\mu$. If⋅$u_{\varepsilon}=\rbr{\op A_{\mu}^{\varepsilon}}^{-1}f$,
then we have
\[
\rbr{\chi^{\prime}\op P_{\mu}^{\varepsilon}u_{\varepsilon},u_{\varepsilon}^{\dual}}_{\Omega}=\rbr{f,\chi^{\prime}u_{\varepsilon}^{\dual}}_{\Omega}=\rbr{\op A_{\mu}^{0}u_{0},\chi^{\prime}u_{\varepsilon}^{\dual}}_{\Omega}.
\]
Thus,
\[
\begin{aligned}\rbr{\chi^{\prime}\op P_{\mu}^{\varepsilon}\rbr{u_{\varepsilon}-u_{0}-\varepsilon U_{\varepsilon,\delta}},u_{\varepsilon}^{\dual}}_{\Omega} & =\rbr{A^{0}Du_{0},D\chi^{\prime}u_{\varepsilon}^{\dual}}_{\Omega}-\rbr{A^{\varepsilon}D\rbr{u_{0}+\varepsilon U_{\varepsilon,\delta}},D\chi^{\prime}u_{\varepsilon}^{\dual}}_{\Omega}\\
 & \quad+\varepsilon\mu\rbr{U_{\varepsilon,\delta},\chi^{\prime}u_{\varepsilon}^{\dual}}_{\Omega}.
\end{aligned}
\]
The⋅first two terms on the right-hand side are similar to those in⋅(\ref{eq: First step}),
with $u_{\varepsilon}^{\dual}$ replaced by⋅$\chi^{\prime}u_{\varepsilon}^{\dual}$,
in which case $\chi_{\varepsilon}\restrict{\supp\chi^{\prime}}=1$
for $5\varepsilon\le\sigma$, so the previous calculations go over
without change to yield, for such~$\varepsilon$,
\begin{equation}
\rbr{\op A_{\mu}^{\varepsilon}}^{-1}\chi^{\prime}\op P_{\mu}^{\varepsilon}\rbr{\rbr{\op A_{\mu}^{\varepsilon}}^{-1}-\rbr{\op A_{\mu}^{0}}^{-1}-\varepsilon\op K_{\mu}^{\varepsilon}}\restrict{L_{p}\rbr{\Omega}^{n}}=\mathring{\op I}_{\mu}^{\varepsilon}+\mathring{\op D}_{\mu}^{\varepsilon},\label{eq: The interior identity}
\end{equation}
where
\[
\begin{aligned}\rbr{\mathring{\op I}_{\mu}^{\varepsilon}f,f^{\dual}}_{\Omega} & =\rbr{\tau^{\varepsilon}\op T^{\varepsilon}D_{1}^{*}A\rbr{D_{1}u_{0,\delta}+D_{2}U_{\delta}},\rbr{\op T^{\varepsilon}-\op I}\chi^{\prime}u_{\varepsilon}^{\dual}}_{\Omega\times\cell}\\
 & \quad-\rbr{\tau^{\varepsilon}\sbr{A,\op T^{\varepsilon}}\rbr{D_{1}u_{0,\delta}+D_{2}U_{\delta}},D_{1}\chi^{\prime}u_{\varepsilon}^{\dual}}_{\Omega\times\cell}\\
 & \quad-\varepsilon\rbr{\tau^{\varepsilon}A\op T^{\varepsilon}D_{1}U_{\delta},D_{1}\chi^{\prime}u_{\varepsilon}^{\dual}}_{\Omega\times\cell}\\
 & \quad-\rbr{A^{\varepsilon}\rbr{\op I-\op S^{\varepsilon}}Du_{0},D\chi^{\prime}u_{\varepsilon}^{\dual}}_{\Omega}\\
 & \quad+\varepsilon\mu\rbr{U_{\varepsilon,\delta},\chi^{\prime}u_{\varepsilon}^{\dual}}_{\Omega}
\end{aligned}
\]
and
\[
\begin{aligned}\rbr{\mathring{\op D}_{\mu}^{\varepsilon}f,f^{\dual}}_{\Omega} & =\rbr{A\rbr{D_{1}\rbr{u_{0}-u_{0,\delta}}+D_{2}\rbr{U-U_{\delta}}},D_{1}\chi^{\prime}u_{\varepsilon}^{\dual}}_{\Omega\times\cell}\\
 & \quad-\rbr{\tau^{\varepsilon}A\op T^{\varepsilon}D_{1}\rbr{u_{0}-u_{0,\delta}},D_{1}\chi^{\prime}u_{\varepsilon}^{\dual}}_{\Omega\times\cell}.
\end{aligned}
\]
This is the interior operator identity that we seek.
\begin{proof}[Proof of Theorem~\textup{\ref{thm: Interior approximation with 1st corrector}}]
Set⋅$v_{\varepsilon}=u_{\varepsilon}-u_{0}-\varepsilon U_{\varepsilon,\delta}$
and⋅$f_{\varepsilon}=\chi^{\prime}\op P_{\mu}^{\varepsilon}v_{\varepsilon}$.
If⋅$\eta$ is a smooth cutoff function which is supported in⋅$\Omega$
and is identically~$1$ on⋅$\supp\chi^{\prime}$, then $\eta v_{\varepsilon}$
belongs to⋅$\set W_{p}^{1}\rbr{\Omega;\C^{n}}$ and therefore $f_{\varepsilon}=\chi^{\prime}\op A_{\mu}^{\varepsilon}\eta v_{\varepsilon}$
belongs to~$\set W_{p}^{-1}\rbr{\Omega;\C^{n}}$. To⋅estimate the
norm of⋅$f_{\varepsilon}$, we use the identity~(\ref{eq: The interior identity}):
\[
\begin{aligned}\bigabs{\rbr{f_{\varepsilon},u_{\varepsilon}^{\dual}}_{\Omega}} & \lesssim\varepsilon^{s}\bigl(\varepsilon^{1-s}\norm{Du_{0,\delta}}_{1,p,\Omega}+\varepsilon^{1-s}\norm{D_{1}U_{\delta}}_{1_{2},p,\Omega\times\cell}+\norm{U_{\delta}}_{1_{2},p,\Omega\times\cell}\\
 & \hphantom{{}\lesssim\varepsilon^{s}\bigl(\varepsilon^{1-s}\norm{Du_{0,\delta}}_{1,p,\Omega}+\varepsilon^{1-s}\norm{D_{1}U_{\delta}}_{1_{2},p,\Omega\times\cell}+\norm{U_{\delta}}_{1_{2},p,\Omega\times\cell}}\mathllap{{}+\varepsilon^{-s}\norm{D\rbr{u_{0}-u_{0,\delta}}}_{p,\Omega}+\varepsilon^{-s}\norm{D_{2}\rbr{U-U_{\delta}}}_{p,\Omega\times\cell}}\bigr)\norm{u_{\varepsilon}^{\dual}}_{1,p^{\dual},\Omega}
\end{aligned}
\]
(cf.⋅the proofs of Lemmas~\ref{lem: Ie} and⋅\ref{lem: De}). Taking
the supremum over all⋅$f^{\dual}\in\rbr{W_{p}^{1}\rbr{\Omega}^{n}}^{*}$,
or, equivalently, over all⋅$u_{\varepsilon}^{\dual}\in\set W_{\smash[t]{\cramped{p^{\dual}}}}^{1}\rbr{\Omega;\C^{n}}$
(recall that the quotient map~(\dualref{def: Quotient map})
is an epimorphism), and  applying Lemmas~\ref{lem: R0_mu(delta) and K_mu(delta) are bounded}
and⋅\ref{lem: R0_mu(delta) and K_mu(delta) convergence rates} shows
that
\begin{equation}
\Norm{f_{\varepsilon}}_{-1,p,\Omega}\lesssim\varepsilon^{s}\norm f_{p,\Omega}.\label{est: Norm of fe}
\end{equation}
On⋅the⋅other hand, according to~(\ref{est: Local enery estimate}),
\[
\norm{D\chi v_{\varepsilon}}_{p,\Omega}\lesssim\norm{v_{\varepsilon}}_{p,\Omega}+\Norm{f_{\varepsilon}}_{-1,p,\Omega},
\]
because $\chi^{\prime}\op A_{\mu}^{\varepsilon}\eta v_{\varepsilon}=f_{\varepsilon}$
and⋅$\eta=1$ on~$\supp\chi^{\prime}$. The⋅result now follows from⋅(\ref{est: Norm of K=001D4B}),
(\ref{est: Norm of fe}) and Theorem~\ref{thm: Convergence with e}.
\end{proof}

\appendix

\section{An estimate for integrals over a neighborhood of the boundary}

The following lemma is a slight modification of~\cite[Lemma~5.1]{PSu:2012}.
\begin{lem}
\label{lem: Norm in the neighborhood}Let⋅$\Sigma$ be a uniformly
weakly Lipschitz domain in~$\R^{d}$. Then for each⋅$q\in[1,\infty)$
and any~$\varepsilon>0$
\begin{equation}
\norm u_{q,\rbr{\partial\Sigma}_{\varepsilon}\cap\Sigma}\lesssim\varepsilon^{1/q}\norm u_{1,q,\Sigma}^{1/q}\norm u_{q,\Sigma}^{1/q^{\dual}},\qquad u\in C_{c}^{\infty}\rbr{\wb{\Sigma}}.\label{est: Norm in the neighborhood}
\end{equation}
The⋅constant in the inequality depends only on⋅$q$,⋅$d$ and~$\Sigma$.
\end{lem}
\begin{proof}
Recall that $B$ denotes the open unit ball centered at the origin
and $B_{+}$ denotes the open unit half-ball with~$x_{d}\in\rbr{0,1}$.
Let $S_{t}$ be the cross-section of⋅$B$ at⋅$x_{d}=t$ and $P_{t}$
be the piece of⋅$B_{+}$ with~$x_{d}\in\rbr{0,t}$. If⋅$\rbr{W_{k},\omega_{k}}$
are local boundary coordinate patches, then $\omega_{k}\rbr{W_{k}\cap\Sigma}=B_{+}$
and⋅$\omega_{k}\rbr{W_{k}\cap\partial\Sigma}=S_{0}$, and for any⋅$y\in\omega_{k}\rbr{W_{k}\cap\Sigma}$
\[
\dist\rbr{y,S_{0}}\le L_{\Sigma}\dist\rbr{x,W_{k}\cap\partial\Sigma},
\]
where~$x=\omega_{k}^{-1}\rbr y$ and~$L_{\Sigma}=\sup_{k}\seminorm{\omega_{k}}_{C^{0,1}}$.
It⋅follows that $\omega_{k}\rbr{W_{k}\cap\rbr{\partial\Sigma}_{\varepsilon}\cap\Sigma}\subset P_{\varepsilon/\varepsilon_{1}}$
with~$\varepsilon_{1}\cellradius=L_{\Sigma}^{-1}$. On⋅the⋅other
hand, we know that the cover is sufficiently tight in the sense that
the union of⋅$\omega_{k}^{-1}\rbr{B_{+}}$ contains $\rbr{\partial\Sigma}_{\delta}\cap\Sigma$
for some~$\delta>0$. Therefore, taking $\varepsilon_{0}=\varepsilon_{1}\meet\delta$,
we can insure that $\rbr{\partial\Sigma}_{\varepsilon}\cap\Sigma$
is covered by⋅$\cbr{W_{k}}$ for any~$\varepsilon\le\varepsilon_{0}$.

Now, using a partition of unity~$\cbr{\varphi_{k}}$ subordinate
to⋅$\cbr{W_{k}}$ (see⋅Section~\ref{sec: Notation}) and making a
change of variables to flatten out the boundary, we reduce⋅(\ref{est: Norm in the neighborhood})
to proving that, for any $\varepsilon\le\varepsilon_{0}$ and any
smooth function~$u$ on⋅$B_{+}$ vanishing near the boundary of⋅$B$,
it holds that 
\begin{equation}
\norm u_{q,P_{\varepsilon/\varepsilon_{0}}}\lesssim\varepsilon^{1/q}\norm u_{1,q,B_{+}}^{1/q}\norm u_{q,B_{+}}^{1-1/q}.\label{est: Norm in the neighborhood for half-ball}
\end{equation}
By⋅the⋅divergence theorem, for any⋅$t\in\rbr{0,1}$ we have
\[
\int_{S_{t}}\abs{u\rbr{x^{\prime},t}}^{q}\dd x^{\prime}=-\int_{B_{+}\!\setminus P_{t}}\partial_{x_{d}}\abs{u\rbr x}^{q}\dd x,
\]
and hence
\[
\begin{aligned}\int_{S_{t}}\abs{u\rbr{x^{\prime},t}}^{q}\dd x^{\prime} & \le q\int_{B_{+}\!\setminus P_{t}}\abs{\partial_{x_{d}}u\rbr x}\abs{u\rbr x}^{q-1}\dd x\\
 & \le q\biggrbr{\int_{B_{+}}\abs{\partial_{x_{d}}u\rbr x}^{q}\dd x}^{1/q}\biggrbr{\int_{B_{+}}\abs{u\rbr x}^{q}\dd x}^{1-1/q}.
\end{aligned}
\]
Integrating in $t$ from $0$ to $\varepsilon/\varepsilon_{0}$ now
gives~(\ref{est: Norm in the neighborhood for half-ball}).
\end{proof}

\section*{Acknowledgment}

The~author is grateful to T.\thinglue A.~Suslina for helpful discussions.

\bibliographystyle{amsalpha}
\bibliography{bibliography}

\begin{thebibliography}{ZhKO94}

\bibitem[AF03]{AF:2003}
\textsc{R.~Adams and J.~Fournier},
\newblock \textit{Sobolev Spaces}, 2nd~ed,
\newblock Academic Press, Amsterdam, 2003.

\bibitem[Agr13]{Agranovich:2013}
\textsc{M.~S.~Agranovich},
\newblock \textit{Sobolev Spaces, Their Generalizations and Elliptic Problems
                  in Smooth and Lipschitz Domains},
\newblock Moscow Center for Continuous Mathematical Education, Moscow, 2013 (in Russian);
\newblock Springer International, 2015 (in English).

\bibitem[A92]{Al:1992}
\textsc{G.~Allaire},
\newblock \textit{Homogenization and two-scale convergence},
\newblock SIAM J. Math. Anal., 23 (1992), pp.~1482--1518.

\bibitem[AC98]{AlaireConca:1998}
\textsc{G.~Allaire and C.~Conca},
\newblock \textit{Bloch wave homogenization and spectral asymptotic analysis},
\newblock J. Math. Pures. Appl., 77 (1998), pp.~153--208.

\bibitem[BP84]{BP:1984}
\textsc{N.~Bakhvalov and G.~Panasenko},
\newblock \textit{Homogenisation: Averaging Processes in Periodic Media:
                  Mathematical Problems in the Mechanics of Composite Materials},
\newblock Nauka, Moscow, 1984 (in Russian);
\newblock Kluwer Academic, Dordrecht, 1989 (in English).

\bibitem[BLP78]{BLP:1978}
\textsc{A.~Bensoussan, J.-L.~Lions and G.~Papanicolaou},
\newblock \textit{Asymptotic Analysis for Periodic Structures},
\newblock North-Holland, Amsterdam, 1978.

\bibitem[BSu01]{BSu:2001}
\textsc{M.\thinglue Sh.~Birman and T.\thinglue A.~Suslina},
\newblock \textit{Threshold effects near the lower edge of the spectrum
                  for periodic differential operators of mathematical physics},
\newblock in Systems, Approximation, Singular Integral Operators, and Related Topics,
          A.\thinglue A.~Borichev and N.\thinglue K.~Nikolski, eds.,
\newblock Birkh\"auser, Basel, 2001, pp.~71--107.

\bibitem[BSu03]{BSu:2003}
\bysame,
\newblock \textit{Second order periodic differential operators. Threshold properties and homogenization},
\newblock Algebra i~Analiz, 15 (2003), no.~5, pp.~1--108 (in Russian);
\newblock St.~Petersburg Math.~J., 15 (2004), pp.~639--714 (in English).

\bibitem[BSu05]{BSu:2005}
\bysame,
\newblock \textit{Homogenization with corrector term for periodic elliptic differential operators},
\newblock Algebra i~Analiz, 17 (2005), no.~6, pp.~1--104 (in Russian);
\newblock St.~Petersburg Math.~J., 17 (2006), pp.~897--973 (in English).

\bibitem[B08]{Bor:2008}
\textsc{D.\thinglue I.~Borisov},
\newblock \textit{Asymptotics for the solutions of elliptic systems with rapidly oscillating coefficients},
\newblock Algebra i~Analiz, 20 (2008), no.~2, pp.~19--42 (in Russian);
\newblock St.~Petersburg Math.~J., 20 (2009), pp.~175--191 (in English).

\bibitem[ChC16]{ChC:2016}
\textsc{K.\thinglue D.~Cherednichenko and S.~Cooper},
\newblock \textit{Resolvent estimates for high-contrast elliptic problems with periodic coefficients},
\newblock Arch. Ration. Mech. Anal., 219 (2016), pp.~1061--1086.

\bibitem[CDG02]{CDG:2002}
\textsc{D.~Cioranescu, A.~Damlamian and G.~Griso},
\newblock \textit{Periodic unfolding and homogenization},
\newblock C.~R.~Acad. Sci. Paris, Ser.~I, 335 (2002), pp.~99--104.

\bibitem[Gar07]{Gar:2007}
\textsc{J.~Garnett},
\newblock \textit{Bounded Analytic Functions},
\newblock Springer, New~York, 2007.

\bibitem[GiM79]{GiM:1979}
\textsc{M.~Giaquinta and G.~Modica},
\newblock \textit{Regularity results for some classes of higher order non linear elliptic systems},
\newblock J.~Reine~u.~angew. Math., 311/312 (1979), pp.~145--169.

\bibitem[Gia83]{Gia:1983}
\textsc{M.~Giaquinta},
\newblock \textit{Multiple Integrals in the Calculus of Variations and Nonlinear Elliptic Systems},
\newblock Princeton University Press, New~Jersey, 1983.

\bibitem[Gra14\textsubscript{1}]{Gra:2014-1}
\textsc{L.~Grafakos},
\newblock \textit{Classical Fourier Analysis},
\newblock 3rd~ed., Springer, New~York, 2014.

\bibitem[Gra14\textsubscript{2}]{Gra:2014-2}
\bysame,
\newblock \textit{Modern Fourier Analysis},
\newblock 3rd~ed., Springer, New~York, 2014.

\bibitem[Gri04]{Gri:2004}
\textsc{G.~Griso},
\newblock \textit{Error estimate and unfolding for periodic homogenization},
\newblock Asymptot. Anal., 40 (2004), pp.~269--286.

\bibitem[Gri06]{Gri:2006}
\bysame,
\newblock \textit{Interior error estimate for periodic homogenization},
\newblock Anal. Appl., 4 (2006), pp.~61--79.

\bibitem[Grv11]{Grisvard:2011}
\textsc{P.~Grisvard},
\newblock \textit{Elliptic Problems in Nonsmooth Domains},
\newblock 2nd~ed., SIAM, Philadelphia, 2011.

\bibitem[JK95]{Jerison+Kenig:1995}
\textsc{D.~Jerison and C.\thinglue E.~Kenig},
\newblock \textit{The inhomogeneous Dirichlet problem in Lipschitz domains},
\newblock J. Funct. Anal., 130 (1995), pp.~161--219.

\bibitem[KLS12]{KLS:2012}
\textsc{C.\thinglue E.~Kenig, F.~Lin and Z.~Shen},
\newblock \textit{Convergence rates in $L_{2}$ for elliptic homogenization problems},
\newblock Arch. Ration. Mech. Anal., 203 (2012), pp.~1009--1036.

\bibitem[MSh09]{MSh:2009}
\textsc{V.~G.~Maz'ya and T.~O.~Shaposhnikova},
\newblock \textit{Theory of Sobolev Multipliers: With Applications to Differential and Integral Operators},
\newblock Springer, Berlin, 2009.

\bibitem[McL00]{McL:2000}
\textsc{W.~McLean},
\newblock \textit{Strongly Elliptic Systems and Boundary Integral Equations},
\newblock Cambridge University Press, Cambridge, 2000.

\bibitem[Mey63]{Mey:1963}
\textsc{N.~G.~Meyers},
\newblock \textit{An $L^p$-estimate for the gradient of solutions of second order elliptic divergence equations},
\newblock Ann. Sc. Norm. Super. Pisa Cl. Sci., 17 (1963), pp.~189--206.

\bibitem[MT97]{MT:1997}
\textsc{F.~Murat and L.~Tartar},
\newblock \textit{$H$-Convergence},
\newblock in Topics in the Mathematical Modelling of Composite Materials,
          A.~Cherkaev and R.~Kohn, eds.,
\newblock Birkh\"auser, Boston, 1997, pp.~21--43.

\bibitem[N06]{Nazarov:2006}
\textsc{S.~A.~Nazarov},
\newblock \textit{Homogenization of elliptic systems with periodic coefficients:
                  Weighted $L^p$ and $L^\infty$ estimates for asymptotic remainders},
\newblock Algebra i~Analiz, 18 (2006), no.~2, pp.~117--166 (in Russian);
\newblock St.~Petersburg Math.~J., 18 (2006), pp.~269--304 (in English).

\bibitem[Ne12]{Necas:2012}
\textsc{J.~Ne\v{c}as},
\newblock \textit{Direct Methods in the Theory of Elliptic Equations},
\newblock Springer, Berlin, 2012.

\bibitem[PSu12]{PSu:2012}
\textsc{M.~A.~Pakhnin and T.~A.~Suslina},
\newblock \textit{Operator error estimates for homogenization
                  of the elliptic Dirichlet problem in a bounded domain},
\newblock Algebra i~Analiz, 24 (2012), no.~6, pp.~139--177 (in Russian);
\newblock St.~Petersburg Math.~J., 24 (2013), pp.~949--976 (in English).

\bibitem[PT07]{PasT:2007}
\textsc{S.\thinglue E.~Pastukhova and R.\thinglue N.~Tikhomirov},
\newblock \textit{Operator estimates in reiterated and locally periodic homogenization},
\newblock Dokl. Acad. Nauk, 415 (2007), pp.~304--309 (in Russian);
\newblock Dokl. Math., 76 (2007), pp.~548--553 (in English).

\bibitem[Sav97]{Savare:1997}
\textsc{G.~Savar\'{e}},
\newblock \textit{Regularity and perturbation results for mixed second order
                  elliptic problems},
\newblock Comm. Partial Differential Equations, 22 (1997), pp.~869--899.

\bibitem[Sav98]{Savare:1998}
\bysame,
\newblock \textit{Regularity results for elliptic equations in Lipschitz domains},
\newblock J. Funct. Anal., 152 (1998), pp.~176--201.

\bibitem[Se17\textsubscript{1}]{Se:2017-1}
\textsc{N.~N.~Senik},
\newblock \textit{Homogenization for non-self-adjoint locally periodic
                  elliptic operators},
\newblock May 2017, \url{https://arxiv.org/abs/1703.02023}.

\bibitem[Se17\textsubscript{2}]{Se:2017-2}
\bysame,
\newblock \textit{Homogenization for non-self-adjoint periodic elliptic
                  operators on an infinite cylinder},
\newblock SIAM J. Math.~Anal., 49 (2017), pp.~874--898.

\bibitem[Se17\textsubscript{3}]{Se:2017-3}
\bysame,
\newblock \textit{On homogenization for non-self-adjoint locally periodic
                  elliptic operators},
\newblock Funktsional. Anal. i Prilozhen., 51 (2017), no.~2, pp.~92--96 (in Russian);
\newblock Funct. Anal. Appl., 51 (2017), pp.~152--156 (in English).

\bibitem[Se20]{Se:2020}
\bysame,
\newblock \textit{On homogenization of locally periodic elliptic and parabolic operators},
\newblock Funktsional. Anal. i Prilozhen., 54 (2020), no.~1, pp.~87--92 (in Russian).

\bibitem[Sha68]{Shamir:1968}
\textsc{E.~Shamir},
\newblock \textit{Regularization of mixed second-order elliptic problems},
\newblock Israel J. Math., 6 (1968), pp.~150--168.

\bibitem[She18]{Shen:2018}
\textsc{Zh.~Shen},
\newblock \textit{Periodic Homogenization of Elliptic Systems},
\newblock Birkh\"auser, Cham, 2018.

\bibitem[Shn74]{Shneiberg:1974}
\textsc{I.~Ya.~Shneiberg},
\newblock \textit{Spectral properties of linear operators in interpolation
                  families of Banach spaces},
\newblock Mat. Issled., 9 (1974), no.~2, pp.~214--227 (in Russian).

\bibitem[Ste70]{St:1970}
\textsc{E.~Stein},
\newblock \textit{Singular Integrals and Differentiability Properties of Functions},
\newblock Princeton University Press, New~Jersey, 1970.

\bibitem[Su13\textsubscript{1}]{Su:2013-1}
\textsc{T.\thinglue A.~Suslina},
\newblock \textit{Homogenization of the Dirichlet problem for elliptic systems:
                  $L^2$\nobreakdash-operator error estimates},
\newblock Mathematika, 59 (2013), pp.~463--476.

\bibitem[Su13\textsubscript{2}]{Su:2013-2}
\bysame,
\newblock \textit{Homogenization of the Neumann problem for elliptic systems
                  with periodic coefficients},
\newblock SIAM J. Math. Anal., 45 (2013), pp.~3453--3493.

\bibitem[Tar10]{Tartar:2010}
\textsc{L.~Tartar},
\newblock \textit{The General Theory of Homogenization},
\newblock Springer, Berlin, 2010.

\bibitem[ZhKO94]{ZhKO:1993}
\textsc{V.\thinglue V.~Jikov, S.\thinglue M.~Kozlov and O.\thinglue A.~Oleinik},
\newblock \textit{Homogenization of Differential Operators and Integral Functionals},
\newblock Springer, Berlin, 1994.

\bibitem[Zh05]{Zh:2005}
\textsc{V.\thinglue V.~Zhikov},
\newblock \textit{On operator estimates in homogenization theory},
\newblock Dokl. Acad. Nauk, 403 (2005), pp.~305--308 (in Russian);
\newblock Dokl. Math., 72 (2005), pp.~535--538 (in English).

\bibitem[ZhP05]{ZhPas:2005}
\textsc{V.\thinglue V.~Zhikov and S.\thinglue E.~Pastukhova},
\newblock \textit{On operator estimates for some problems in homogenization theory},
\newblock Russ.~J. Math. Phys., 12 (2005), pp.~515--524.

\bibitem[ZhP16]{ZhPas:2016}
\bysame,
\newblock \textit{Operator estimates in homogenization theory},
\newblock Uspekhi Mat. Nauk, 71 (2016), no.~3, pp.~27--122 (in Russian);
\newblock Russian Math. Surveys, 71 (2016), pp.~417--511 (in English).

\end{thebibliography}

\end{document}